\documentclass[12pt,leqno]{article}
\usepackage{amsfonts}
\usepackage{amsmath, amsthm, amsfonts, amssymb, color}
\usepackage{mathrsfs}
\usepackage{url}
\usepackage{color}
\theoremstyle{definition}

\newcommand{\scr}[1]{\mathscr #1}
\definecolor{wco}{rgb}{0.5,0.2,0.3}

\numberwithin{equation}{section} \theoremstyle{remark}

\newcommand{\ua}{\uparrow}

\title{{\bf   Distribution Dependent Reflecting Stochastic Differential Equations}\footnote{Supported in
 part by   NNSFC (11771326, 11831014, 11921001).}}
\author{
{\bf    Feng-Yu Wang  }\\
\footnotesize{Center for Applied Mathematics, Tianjin University, Tianjin 300072, China } }

\begin{document}
\allowdisplaybreaks
\def\R{\mathbb R}  \def\ff{\frac} \def\ss{\sqrt} \def\B{\mathbf
B}\def\TO{\mathbb T}
\def\I{\mathbb I_{\pp M}}\def\p<{\preceq}
\def\N{\mathbb N} \def\kk{\kappa} \def\m{{\bf m}}
\def\ee{\varepsilon}\def\ddd{D^*}
\def\dd{\delta} \def\DD{\Delta} \def\vv{\varepsilon} \def\rr{\rho}
\def\<{\langle} \def\>{\rangle} \def\GG{\Gamma} \def\gg{\gamma}
  \def\nn{\nabla} \def\pp{\partial} \def\E{\mathbb E}
\def\d{\text{\rm{d}}} \def\bb{\beta} \def\aa{\alpha} \def\D{\scr D}
  \def\si{\sigma} \def\ess{\text{\rm{ess}}}
\def\beg{\begin} \def\beq{\begin{equation}}  \def\F{\scr F}
\def\Ric{{\rm Ric}} \def\Hess{\text{\rm{Hess}}}
\def\e{\text{\rm{e}}} \def\ua{\underline a} \def\OO{\Omega}  \def\oo{\omega}
 \def\tt{\tilde}
\def\cut{\text{\rm{cut}}} \def\P{\mathbb P} \def\ifn{I_n(f^{\bigotimes n})}
\def\C{\scr C}      \def\aaa{\mathbf{r}}     \def\r{r}
\def\gap{\text{\rm{gap}}} \def\prr{\pi_{{\bf m},\varrho}}  \def\r{\mathbf r}
\def\Z{\mathbb Z} \def\vrr{\varrho} \def\ll{\lambda}
\def\L{\scr L}\def\Tt{\tt} \def\TT{\tt}\def\II{\mathbb I}
\def\i{{\rm in}}\def\Sect{{\rm Sect}}  \def\H{\mathbb H}
\def\M{\scr M}\def\Q{\mathbb Q} \def\texto{\text{o}} \def\LL{\Lambda}
\def\Rank{{\rm Rank}} \def\B{\scr B} \def\i{{\rm i}} \def\HR{\hat{\R}^d}
\def\to{\rightarrow}\def\l{\ell}\def\iint{\int}
\def\EE{\scr E}\def\Cut{{\rm Cut}}\def\W{\mathbb W}
\def\A{\scr A} \def\Lip{{\rm Lip}}\def\S{\scr S}
\def\BB{\scr B}\def\Ent{{\rm Ent}} \def\i{{\rm i}}\def\itparallel{{\it\parallel}}
\def\g{{\mathbf g}}\def\Sect{{\mathcal Sec}}\def\T{\mathcal T}\def\V{{\mathbb V}}
\def\PP{{\bf P}}\def\HL{{\bf L}}\def\Id{{\rm Id}}\def\f{{\bf f}}\def\cut{{\rm cut}}
\def\n{{\mathbf n}}
\def\sm{\preceq}
\def\BL{\scr A} \def\I{{\mathbf I}}

\maketitle

\begin{abstract} To characterize the Neumann problem for nonlinear Fokker-Planck equations, we investigate  distribution dependent reflecting  SDEs (DDRSDEs)  in a domain.
We  first  prove the   well-posedness and  establish  functional inequalities  for reflecting   SDEs with singular drifts, then extend  these results  to   DDRSDEs with singular or monotone coefficients, for which  a general criterion deducing the well-posedness of DDRSDEs from that of reflecting SDEs is established.  

\end{abstract} \noindent
 \vskip 2cm

\section{Introduction }

Because of  intrinsic links to nonlinear Fokker-Planck equations/mean-field particle systems  and many other applications, distribution dependent (McKean-Vlasov)
 SDEs have been intensively investigated,    see  for instances  the monograph/surveys  \cite{CD, 20HRW, SN} among many other references.
To characterize the Neumann problem for nonlinear Fokker-Planck equations in a domain,
we aim to  develop a counterpart theory   for   distribution dependent reflecting    SDEs  (DDRSDEs for short). 

The only reference we know on this topic is  \cite{Reis}, where DDRSDEs  are studied   in a convex domain for coefficients satisfying the $\W_2$-Lipschitz   condition in the distribution  variable and the semi-Lipschitz   condition in the space variable.   We will work on a general framework where $D$ may be non-convex and the coefficients could be singular   in both space and distribution variables.

We first state the fundamental assumption on the domain   in the study of reflecting SDEs, then introduce the link of DDRSDEs and nonlinear Neumann problems,     and finally summarize the main results derived in the paper with an example  of  (singular) granular media equation with Neumann boundary.

\subsection{Assumption on the domain}  Let $D\subset \R^d$ be a connected  open domain with boundary $\pp D$. For any $x\in\pp D$ and $r>0$, let
$$\scr N_{x,r}:=\big\{\n\in \R^d: |\n|=1, B(x-r\n, r)\cap D=\emptyset\big\},$$ where $B(x,r):=\{y\in \R^d: |x-y|<r\}.$ Since  $\scr N_{x,r}$ is decreasing in $r>0$, we have
 $$\scr N_x:=\cup_{r>0} \scr N_{x,r}=\lim_{r\downarrow 0} \scr N_{x,r}, \ \ x\in \pp D.$$ We call $\scr N_x$ the set of inward unit normal vectors of $\pp D$ at point $x$.
When $\pp D$ is differentiable at $x$, $\scr N_x$ is a singleton set.  Otherwise $\scr N_x$ may be empty or contain more than one    vectors.  For instance, letting $D$ be the interior of a triangle in $\R^2$,
at each vertex $x$ the set $\scr N_x$ contains infinite many vectors, whereas for $D$ being the exterior of the triangle $\scr N_x$ is empty at each vertex point $x$.

Following \cite{LSZ,SS}, throughout the paper we make the following assumption on $D$, which automatically holds for $D=\R^d$ where $\pp D=\emptyset$.

\beg{enumerate} \item[{\bf (D)}]   Either $D$ is convex, or there exists a constant $r_0>0$ such that $\scr N_x=\scr N_{x, r_0}\ne\emptyset$   and
\beq\label{OI}  \sup_{v\in \R^d, |v|=1} \inf\big\{ \<v, \n(y)\>:    y \in B(x,r_0)\cap \pp D, \n(y)\in \scr N_y\big\}\ge r_0,\ \ x\in    \pp D.\end{equation}
\end{enumerate}

 \paragraph{Remark 2.1.}   We present below some facts on assumption  {\bf (D)}.
 \beg{enumerate} \item[(1)] According to \cite[Remark 1.1]{SS},  for any $x\in\pp D$ and $r>0$, $\n\in \scr N_{x,r}$   if and only if
 $\<y-x, \n\>\ge -\ff {|y-x|^2}  {2 r}$ for $\in \bar D,$  so that the condition
  $\scr N_x=\scr N_{x, r_0}$ in {\bf (D)} implies
 \beq\label{LCC} \<y-x, \n(x)\>\ge -\ff {|y-x|^2}  {2 r_0},\ \ y\in \bar D, x\in\pp D,  \n(x)\in \scr N_x.\end{equation}
When $D$ is convex,   {\bf (D)} holds for any $r_0>0$ so that
  \beq\label{CVX} \<y-x, \n(x)\>\ge 0,\ \ y\in\bar D, x\in\pp D, \n(x)\in \scr N_x,\end{equation} and  \eqref{OI} holds if $d=2$ or $D$ is bounded,   see \cite{Tanaka}.
 \item[(2)]   When $\pp D$ is $C^{1}$-smooth, for each $x\in\pp D$ the set $\scr N_x$ is  singleton.
 If $\n(x)\in \scr N_x$ is uniformly continuous in $x\in \pp D$, then \eqref{OI} holds for small $r_0>0$.
 In particular,  {\bf (D)} holds when  $\pp D\in C_b^2$ in the following sense. \end{enumerate}

\beg{defn}\label{Cb2} For any $r>0$, let
 \beg{align*} &\pp_r D:=\big\{x\in\bar D:  {\rm dist}(x,\pp D)\le r\big\},\ \ \
 \pp_{-r} D:=\big\{x\in D^c:  {\rm dist}(x,\pp D)\le r\big\},\\
 &\pp_{\pm r} D:= (\pp_r D)\cup \pp_{-r}D,\ \ \ D_r:= D\cup (\pp_{-r}D).\end{align*}
For any  $k\in\mathbb N$, we write $\pp D\in C^k_{b}$  if there exists a constant $r_0>0$ such that  the polar coordinate map
$$ I: \pp D\times [-r_0,r_0] \ni (\theta, \rr)\mapsto   \theta+ \rr\n(\theta) \in  \pp_{\pm r_0}D$$
is a $C^k$-diffeomorphism, such that   $(\theta(x),\rr(x)):= I^{-1}(x)$  having  bounded and continuous derivatives in $x\in  \pp_{\pm r_0}D$ up to the $k$-th order, where $\theta(x)$  is   the projection  of $x$ to $\pp D$ and
\beq\label{RR} \rr(x)=  {\rm dist}(x,\pp D)1_{\{\pp_{r_0}D\}}(x)- {\rm dist}(x,\pp D)1_{\{\pp_{-r_0}D\}}(x),\ \  x\in \pp_{\pm r_0}D. \end{equation}
Moreover, for $\vv\in (0,1)$, we denote $\pp D\in C^{k+\vv}_b$ if  it is in $C_b^k$  with $ \nn^k\rr$ and $\nn^k\theta $ being  $\vv$-H\"older continuous
on $\pp_{r_0}D$.
Finally, we write $\pp D\in C_b^{k,L}$ if it is $C_b^k$ with $\nn^k\rr$ being Lipschitiz continuous on $\pp_{r_0}D.$
\end{defn}

Note that $\pp D\in C_b^k$ does not imply the boundedness of $D$ or $\pp D$, but any bounded $C^k$ domain  satisfies $\pp D\in C_b^k$.  

\subsection{DDRSDE and nonlinear Neumann problem}

Let $\scr P(\bar D)$ be the space of all probability measures on the closure $\bar D$ of $D$,   equipped with the weak topology.
 Consider the following DDRSDE   on $\bar D\subset \R^d$:
\beq\label{E1} \d X_t= b_t(X_t,\L_{X_t}) \d t+ \si_t(X_t,\L_{X_t})\d W_t + \n(X_t) \d l_t,\ \ t\ge 0,\end{equation}
 where   $(W_t)_{t\ge 0}$ is an $m$-dimensional Brownian motion on a complete filtration probability space $(\OO,\{\F_t\}_{t\ge 0},\P)$, $\L_{X_t}$ is the distribution of $X_t$,
 $\n(x)\in\scr N_x$ for $x\in \pp D$, $l_t$ is an adapted continuous increasing process which increases only when $X_t\in \pp D$, and
$$b: [0,\infty)\times D\times \scr P(\bar D) \to \R^d,\ \  \si: [0,\infty)\times D \times \scr P(\bar D)\to \R^d\otimes\R^m$$ are measurable.
When different probability measures are considered, we denote by $\L_{X|\P}$ the distribution of a random variable $X$ under the probability $\P$.

\beg{defn}  (1) A pair $(X_t,l_t)_{t\ge 0}$ is called a solution of $\eqref{E1}$, if $X_t$ is an adapted continuous process on $\bar D$, $l_t$ is an adapted continuous increasing process with   $\d l_t$ supported on $\{t\ge 0: X_t\in \pp D\}$, such that $\P$-a.s.
$$  \int_0^t \big\{|b_r(X_r, \L_{X_r})| + \| \si_r(X_r,\L_{X_r})\|^2\} \d r <\infty,\ \ t\ge 0,$$ and  for some measurable map  $\pp D\ni x\mapsto \n(x)\in \scr N_x,$ $\P$-a.s.
$$ X_t= X_0 +\int_0^t b_r(X_r, \L_{X_r})\d r + \int_0^t \si_r(X_r,\L_{X_r})\d W_r+\int_0^t \n(X_r)\d l_r,\ \ t\ge 0.$$
In this case, $l_t$ is called the local time of $X_t$ on $\pp D$. We call \eqref{E1}  strongly well-posed for distributions in a subspace $\hat{\scr P}\subset \scr P(\bar D),$ if for any  $\F_0$-measurable
variable $X_0$ with $\L_{X_0}\in \hat {\scr P}$,  the equation has a unique solution   with $\L_{X_t}\in \hat {\scr P}$ for $t\ge 0$; if this is true for $\hat {\scr P}=\scr P(\bar D)$, we called it strongly well-posed.

(2) A triple $(X_t,l_t, W_t)_{t\ge 0}$ is called a weak solution of \eqref{E1}, if $W_t$ is an $m$-dimensional Brownian motion under a probability space and $(X_t,l_t)_{t\ge 0}$ solves \eqref{E1}.
\eqref{E1} is called  weakly unique (resp. jointly weakly unique), if for any   two weak solutions
  $(X_t,l_t, W_t)_{t\ge 0}$ under probability $\P$ and $(\tt X_t,\tt l_t, \tt W_t)_{t\ge 0}$ under probability $\tt \P$,
$\L_{X_0|\P}=\L_{\tt X_0|\tt\P}$ implies $\L_{(X_t,l_t)_{t\ge 0}|\P} = \L_{(\tt X_t,\tt l_t)_{t\ge 0}|\tt\P}$ (resp. $\L_{(X_t,l_t,W_t)_{t\ge 0}|\P} = \L_{(\tt X_t,\tt l_t,\tt W_t)_{t\ge 0}|\tt\P}).$ We call \eqref{E1}   weakly well-posed for distributions in $\hat {\scr P}\subset \scr P(\bar D)$, if   it has a unique weak solution for initial distributions in $\hat{\scr P}$ and the   distribution of the solution at any time is in $\hat {\scr P}$; it is called weakly well-posed  if moreover  $\hat {\scr P}=\scr P(\bar D).$

(3) We call  \eqref{E1} well-posed (for distributions in $\hat{\scr P})$, if it is both strongly and weakly well-posed (for distributions in $\hat{\scr P})$. \end{defn}

To characterize the  nonlinear Fokker-Planck equation  associated with \eqref{E1},  consider the following time-distribution dependent second order differential operator:
\beq\label{LMU} L_{t,\mu}:= \ff 1 2  {\rm tr} \big\{(\si_t\si_t^*)(\cdot,\mu) \nn^2\big\} +\nn_{b_t(\cdot,\mu)},\ \  t\ge 0,\mu\in \scr P(\bar D),\end{equation}
where $\nn$ and $\nn^2$ are the gradient and Hessian operators in $\R^d$ respectively, and $\nn_v$ is  the directional derivative along $v\in \R^d$.
Assume that for any $\mu\in C([0,\infty); \scr P(\bar D)),$
\beq\label{SBB} \si^\mu_t(x):=\si_t(x,\mu_t),\ \ b_t^\mu(x):=b_t(x,\mu_t)\end{equation}
satisfy $\|\si^\mu\|^2+|b^\mu|\in L^1_{loc} ([0,\infty)\times \bar D;\d t\,\mu_t(\d x)).$

Let  $C_N^2(\bar D)$ be the class of $C^2$-functions on $\bar D$ with compact support   satisfying the Neumann boundary condition
  $\nn_\n f|_{\pp D}=0$.
 By It\^o's formula,  for any (weak) solution $X_t$ to \eqref{E1}, $\mu_t:=\L_{X_t}$ solves the nonlinear Fokker-Planck equation
 \beq\label{NFP} \pp_t \mu_t= L_{t,\mu_t}^*\mu_t\ \ \text{with\ respect\ to\ }C_N^2(\bar D),\ \ t\ge 0   \end{equation} for probability measures on $\bar D,$ in the sense that $\mu_\cdot\in C([0,\infty);\scr P(\bar D))$ and
 \beq\label{NFP2}
  \mu_t(f):=\int_{\bar D}f\d\mu_t= \mu_0(f)+\int_0^t \mu_s(L_{s,\mu_s} f) \d s,\ \ t\ge 0, f\in C_N^2(\bar D). \end{equation}
On the  other hand,  by establishing the $``$superposition principle" as in  \cite{RWBR18,RWBR} based on \cite{TR},  under reasonable conditions we may prove that a solution to \eqref{NFP}
also provides a weak solution to \eqref{E1}. We leave this to a future study.

To understand \eqref{NFP} as a nonlinear  Neumann problem on $D$, let   $L_{t,\mu_t}^*$  be the adjoint operator of  $L_{t,\mu_t}$:  for any $g\in  L_{loc}^1(D, (\|\si_t(x,\mu_t)\|^2+|b_t(x,\mu_t)|)\d x)$, $L_{t,\mu_t}^*g$ is the linear functional on $C_0^2(D)$ (the class of $C^2$-functions on $D$ with compact support) given by
\beq\label{INTT} C_0^2(D) \ni f\mapsto \int_{D} \{f L_{t,\mu_t}^*g\}(x)  \d x := \int_{D} \{g L_{t,\mu_t}f\}(x) \d x.\end{equation}
Assume that  $\L_{X_t}$ has a density function $\rr_t$, i.e. $\mu_t:=\L_{X_t}=\rr_t(x)\d x.$ It is the case under a general non-degenerate or
 H\"ormander condition (see for instance \cite{BKRS}), and it  follows from   Krylov's estimates \eqref{KR} or \eqref{KR-2} below. When $\pp D\in C^2$,    \eqref{NFP} implies that $\rr_t$ solves the following nonlinear Neumann problem on $\bar D$:
\beq\label{NFP3} \pp_t \rr_t = L_{t,\rr_t}^*\rr_t,\ \ \nn_{t,\n}\rr_t|_{\pp D}=0,\    t\ge 0 \end{equation} in the weak sense,
  where $L_{t,\rr_t}:= L_{t,\rr_t(x)\d x}$, and for a function $g$ on $\pp D$
$$\nn_{t,\n}g:= \nn_{\si_t\si_t^*\n} g+ {\rm div}_{\pp D} (g\pi \si_t\si_t^*\n)$$
for   the  divergence ${\rm div}_{\pp D}$ on $\pp D$ and the projection $\pi$ to the tangent space of $\pp D$:
$$\pi_x v:=v- \<v,\n(x)\>\n(x),\ \ v\in\R^d, x\in\pp D. $$
If in particular $\si\si^*\n=\ll\n$ holds on $[0,\infty)\times \pp D$ for a function $\ll\ne 0$ a.e., $\nn_{t,\n}\rr_t|_{\pp D}=0$ is equivalent to
the standard Neumann boundary condition $\nn_\n \rr_t|_{\pp D}=0.$

We now deduce  \eqref{NFP3} from  \eqref{NFP2}.
Firstly,  by \eqref{INTT},  \eqref{NFP2}  implies
 $$\int_D (f\rr_t)(x)\d x=  \int_D (f\rr_0)(x)\d x+\int_0^t \d s \int_D (fL_{s,\rr_s}^*\rr_s)(x)\d x,\ \ f\in C_0^2(D), t\ge 0,$$ so that
 $\pp_t\rr_t= L_{t,\rr_t}^*\rr_t$. Next,  by  the integration by parts formula, \eqref{NFP2} implies
\beg{align*} &\int_D (f\rr_t)(x)\d x=  \int_D (f\rr_0)(x)\d x+\int_0^t  \d s \int_D (\rr_s L_{s,\rr_s}f)(x)\d x\\
&=  \int_D (f\rr_0)(x)\d x+\int_0^t  \bigg( \int_D (fL_{s,\rr_s}^*\rr_s)(x)\d x+\int_{\pp D} \big\{f \nn_{\si_s\si_s^*\n} \rr_s-  \rr_s\nn_{\si_s\si_s^*\n} f \big\}(x)\d x\bigg)\d s\\
&= \int_D (f\rr_0)(x)\d x+\int_0^t  \bigg( \int_D (f\pp_s\rr_s)(x)\d x+\int_{\pp D} \big\{f\nn_{\si_s\si_s^*\n}\rr_s +f {\rm div}_{\pp D} (\rr_s \pi \si_s\si_s^*\n)  \big\}(x)\d x\bigg)\d s\\
&= \int_D (f\rr_t)(x)\d x + \int_0^t   \d s  \int_{\pp D}   \big\{f(\nn_{t,\n} \rr_t)\big\} (x)\d x,\ \ f\in C_N^2(\bar D), t\ge 0.\end{align*}
Thus, $\nn_{t,\n} \rr_t|_{\pp D}=0.$

\subsection{Summary of main results}

 Theorems \ref{T2.1}-\ref{T2.3} provide sufficient conditions for the well-posedness and functional inequalities of reflecting SDEs with singular drifts. These results generalize    the corresponding ones derived in recent years for singular SDEs without reflection, and improve some existing results for reflecting SDEs. The essential difficulty in the study of singular reflecting SDEs is explained in the beginning of Section 2.

Theorems \ref{T1}-\ref{T3} present the weak and strong well-posedness of the DDRSDE \eqref{E1} under different conditions, where the first  result  applies to  locally integrable   drifts with the distribution dependence bounded by $\|\cdot\|_{k,var}+\W_k$ (see Section 2 for definitions of  probability distances), the second result includes a general criterion deducing the well-posedness of \eqref{E1} from that of reflecting SDEs,  and the last two results  work for   the monotone case with the   dependence on distribution given by   $\W_k (k>1)$ or more general     $\W_\psi$ induced by a cost function $\psi$.

 Theorems \ref{T4} and \ref{T5} establish the log-Harnack inequality for  solutions to \eqref{E1} with respect to the initial distributions, which in particular implies the gradient estimate
and entropy-cost inequality  for the distributions of  the solutions.  The first result applies to  the singular case and the other works for the monotone case.

 \ \newline To conclude this section,    we consider an example   of \eqref{NFP3} arising from kinetic mechanics. For simplicity, we only consider  bounded domain, but  our general results also work for unbounded domains. See \cite{W21c} for the study of exponential ergodicity.

\paragraph{Example 1.1 (Granular media equation with Neumann boundary).} Let $D$ be a bounded domain with $\pp D\in C_b^{2,L}.$ For a potential $V:\bar D\to\R$ and  an interaction functional $W: \R^d\to\R,$ consider
the following nonlinear PDE for probability density functions on $\bar D$:
$$  \pp_t  \varrho_t= \DD  \varrho_t +{\rm div} \big\{ \varrho_t\nn V +  \varrho_t \nn (W* \varrho_t)\big\}, \ \ \nn_\n \varrho_t|_{\pp D}=0,$$
 where $(W* \varrho_t)(x):=\int_{\R^d} W(x-z)  \varrho_t(z)\d z.$
It is easy to see that this equation is covered by \eqref{NFP3}   with
 $$b(x,\mu)= -\nn V(x) -\nn (W*\mu)(x) ,\ \ \si(x,\mu)= \ss 2 {\bf I}_d,$$ where ${\bf I}_d$ is the $d\times d$ identity matrix, and $(W*\mu)(x):= \int_{\R^d} W(x-z)\mu(\d z).$
 
   If $V$ and $W$ are  weakly differentiable with $\|\nn W\|_\infty<\infty$ and $|\nn V|\in L^p(\bar D)$ for some $p>d\lor 2$, then Theorem \ref{T1} with $k=0$ implies that the associated SDE \eqref{E1} is well-posed, and Theorem \ref{T5} provides some  functional   inequalities for the solution. These results apply to   $W(x):=|x|^3$ which is of special interest from  physics \cite{BCP97}.
 
\section{Reflecting SDE with singular drift}

Let $\si_t(x,\mu)=\si_t(x)$ and $b_t(x,\mu)=b_t(x)$ do not depend on $\mu$, so that \eqref{E1} reduces to   the following reflecting SDE on $\bar D$:
 \beq\label{E01} \d X_t= b_t(X_t)\d t +\si_t(X_t) \d W_t+\n(X_t)\d l_t,\ \ t\in [0,T],\end{equation}
 where $T>0$ is a fixed time.
The associated time dependent generator   reads
 \beq\label{LLT} L_t:= \ff 1 2 {\rm tr}\big\{\si_t\si_t^*\nn^2\big\} +\nn_{b_t},\ \ t\in [0,T].\end{equation}
 The problem of confining a stochastic process to a domain goes back to Skorokhod \cite{SK1, SK2}, and has been well developed under monotone (or locally semi-Lipschitz) conditions, see 
 the recent work \cite{Hino} and references within.   In this section, we  solve \eqref{E01} with a singular (unbounded on bounded sets) drift.

 SDEs with singular coefficients have already  been well investigated by using Zvokin's transform, see for instances \cite{KR, XXZZ, XZ, Z2} and references within. However, the corresponding study for singular reflecting   SDEs is very limited. With great effort overcoming difficulty induced by the local time, 
in the recent work   \cite{YZ}  Yang and Zhang were  able to prove  the well-posedness of \eqref{E01} for bounded $C^3$ domain, bounded $b$ and $\si= {\bf I}_d.$  So, the general  setup we discussed here is new in the literature.

 Before  moving on, let us explain the main difficulty of the study by considering   the following simple reflecting SDE on $\bar D$:
\beq\label{E*4} \d X_t = b_t(X_t)\d t +\ss 2 \d W_t + \n(X_t)\d l_t,\ \ t\in [0,T],\end{equation}
where $W_t$ is the $d$-dimensional Brownian motion and $ \int_0^T \|b_t\|_{L^p(\R^d)}^q\d t  <\infty$ for some $p,q>2$ with $\ff d p+\ff 2 q<1.$  When $\ll>0$ is large enough,
the unique solution of the PDE
$$(\pp_t +\DD+\nn_{b_t})u_t=\ll u_t-b_t,\ \ t\in [0,T], u_T=0$$
 satisfies
$$\|u\|_\infty+\|\nn u\|_\infty\le \ff 1 2,\ \ \|u\|_{L_q^p}:= \bigg(\int_0^T\|\nn^2 u_t\|_{L^p(\R^d)}^q\d t\bigg)^{\ff 1 q} <\infty,$$ see \cite{KR,Z2}.
Thus, for any $t\in [0,T]$, $\Theta_t:=id + u_t$ ($id$ is the identity map) is a homeomorphism on $\R^d$, and by It\^o's formula,    $Y_t:= \Theta_t(X_t)$ solves
$$\d Y_t= \ll \{u_t\circ \Theta_t^{-1}\}(Y_t)\d t + \d W_t  + \{(\nn u_t)\circ\Theta_t^{-1}\}(Y_t)\d W_t + \{\n(X_t)+\nn_{\n}u_t(X_t)\}\d l_t.$$
When $D=\R^d$, we have $l_t=0$ so that this SDE is regular enough to have well-posedness, which implies the same property of \eqref{E*4} since $\Theta_t$ is a homeomorphism, see \cite{KR}. When $D\ne \R^d$, to prove the pathwise uniqueness of $Y_t$ by applying It\^o's formula to $|Y_t-\tt Y_t|^2$, where $\tt Y_t:=\Theta_t(\tt X_t)$ for another solution $\tt X_t$ of \eqref{E*4} with local time $\tt l_t$, one needs to
find a constant $c>0$ such that
\beq\label{ERT} \beg{split} &\<\Theta_t(X_t)-\Theta_t(\tt X_t), (\n+\nn_\n u_t)(X_t) \>\d l_t+ \<\Theta_t(\tt X_t)-\Theta_t(X_t), (\n+\nn_\n u_t) (\tt X_t))\>\d \tt l_t\\
&\le c|X_t-\tt X_t|^2 (\d l_t+\d\tt l_t).\end{split} \end{equation}
This is not implied by \eqref{LCC}  except   for   $d=1,$  since only in this case the vectors $\Theta_t(x)-\Theta_t(y)$ and $(\n+\nn_\n u_t)(x)$ are in the same directions of $x-y$ and $\n(x)$ respectively for large   $\ll>0$.

 To overcome this difficulty, we will construct  a  Zvokin's transform by solving the associated Neumann problem on $\bar D$, for which
 $\nn_\n u_t|_{\pp D}=0$. Even in this case, $\Theta_t$ may also  map a point from $\bar D$ to $\bar D^c$ such that \eqref{LCC} does not apply.
 To this end, we will construct a modified process  of $|X_t-\tt X_t|^2$    by using a function   from \cite{DI}. Our construction  simplifies that   in \cite{YZ} and enables us to work in a more general framework.  
 
\subsection{Conditions and main results}

We first  recall some functional spaces used in the study of singular SDEs, see for instance  \cite{XXZZ}. For any $p\ge 1$, $L^p(\R^d)$ is the class of   measurable  functions $f$ on $\R^d$ such that
 $$\|f\|_{L^p(\R^d)}:=\bigg(\int_{\R^d}|f(x)|^p\d x\bigg)^{\ff 1 p}<\infty.$$
For any $\epsilon >0$ and $p\ge 1$, let  $H^{\epsilon,p}(\R^d):=(1-\DD)^{-\ff\epsilon 2} L^p(\R^d)$ with
 $$\|f\|_{H^{\epsilon,p}(\R^d)}:= \|(1-\DD)^{\ff\epsilon 2} f\|_{L^p(\R^d)}<\infty,\ \ f\in  H^{\epsilon,p}(\R^d).$$

 For any $z\in\R^d$ and $r>0$,   let $B(z,r):=\{x\in\R^d: |x-z|< r\}$ be the open ball centered at $z$ with radius $r$.
For any $p,q>1$ and $t_0<t_1$, let $\tt L_q^p(t_0,t_1)$ denote the class of measurable functions $f$ on $[t_0,t_1]\times\R^d$ such that
$$\|f\|_{\tt L_q^p(t_0,t_1)}:= \sup_{z\in \R^d}\bigg( \int_{t_0}^{t_1} \|1_{B(z,1)}f_t\|_{L^p(\R^d)}^q\d t\bigg)^{\ff 1 q}<\infty.$$
For any $\epsilon>0$, let $\tt H_{q}^{\epsilon,p}(t_0,t_1)$ be the space of $f\in \tt L_q^p(t_0,t_1)$ with
$$\|f\|_{\tt H_q^{\epsilon,p}(t_0,t_1)}:=  \sup_{z\in \R^d}\bigg( \int_{t_0}^{t_1}\|g(z+\cdot)f_t\|_{\H^{\epsilon,p} (\R^d)}^q\d t\bigg)^{\ff 1 q}<\infty$$ for some $g\in C_0^\infty(\R^d)$
satisfying $g|_{B(0,1)}=1,$ where $C_0^\infty(\R^d)$ is the class of $C^\infty$ functions on $\R^d$ with compact support.
We remark that the space $\tt H_{q}^{\epsilon,p}(t_0,t_1)$ does not depend on the choice of $g$.
When $t_0=0$, we simply denote
$$\tt L_q^p(t_1):=\tt L_q^p(0,t_1),\ \  \tt H_q^{\epsilon,p}(t_1):=\tt H_q^{\epsilon,p}(0,t_1), \  \ \ t_1>0.$$
For a   domain  $D\subset \R^d$, we denote $f\in \tt L_q^p(t_0,t_1,D) (=: \tt L_q^p(t_1,D)$ for $t_0=0)$, if $f$ is a measurable function on $[t_0,t_1]\times \bar D$  such that
$$\|f\|_{\tt L_q^p(t_0,t_1,D)}:= \|1_D f\|_{\tt L_q^p(t_0,t_1)}<\infty.$$
A vector or matrix valued function   is said in one of the above introduced spaces, if so are its components.

We will  take $(p,q)$ from the class
$$\scr K:=\Big\{(p,q): p,q\in  (1,\infty),\  \ff d p+\ff 2 q<1\Big\},$$
and use the following  assumptions on  the coefficients $b$ and $\si$. Let $\|\cdot\|_\infty$ denote the uniform norm for real (or vector/matrix) valued functions.

\emph{\beg{enumerate}
\item[{$(A_0^{\si,b})$}] {\bf (D)} holds,   $a:= \si\si^*$ and $b$   are extended    to  measurable functions on $[0,T]\times\R^d$,    $b$ has decomposition    $b=b^{(0)}+b^{(1)}$ with  $b^{(0)}_t|_{\bar D^c}=0,$
such that the following conditions hold:
\item[$(1)$]     $a_t$ is invertible with $\|a\|_\infty+\|a^{-1}\|_\infty<\infty$,   and
 \beq\label{*UN}  \lim_{\vv\to 0} \sup_{|x-y|\le \vv, t\in [0,T]} \|a_t(x)-a_t(y)\|=0.\end{equation}
\item[$(2)$]  There exists  $(p_0,q_0)\in \scr K$  such that  $|b^{(0)}|\in \tt L_{q_0}^{p_0}(T).$
Moreover, $b^{(1)}$ is locally bounded on $[0,T]\times \R^d$,   and
    there exist  a constant   $L >1$ and a function $\tt\rr\in C_b^2(\bar D)$  such that
\beq\label{LIP} \|\nn b^{(1)}\|_{\infty }  :=\sup_{t\in [0,T], x\ne y} \ff{|b_t^{(1)}(x)-b_t^{(1)}(y)|}{|x-y|}  \le L,\end{equation}
\beq\label{GB1} \<b_t^{(1)}, \nn \tt \rr\>|_{\bar D} \ge - L, \  \ \<\nn \tt \rr,\n\>|_{\pp D}\ge 1, \ \ t\in [0,T].\end{equation}
 \item[{$(A_1^{\si,b})$}] $(A_0^{\si,b})$ holds, and there exist $\{(p_i,q_i)\}_{0\le i\le l} \subset \scr K$  and $0\le f_i\in\tt L_{q_i}^{p_i}(T), 1\le i\le l,$ such that
 $$|b^{(0)}|^2\in \tt L_{q_0}^{p_0}(T),\ \ \|\nn\si\|^2\le \sum_{i=1}^l f_i.$$
\end{enumerate}}

\paragraph{Remark 2.1.} Each of the following two  conditions implies the existence of $\tt \rr$ in \eqref{GB1}:
\beg{enumerate} \item[(a)] $\pp D\in C_b^2$ and   there exists a constant $K>0$ such that
$\<b_t^{(1)},\n\>|_{\pp D} \ge -K$ for $t\in [0,T];$
\item[(b)]  $D$ is bounded and  there exist   $\vv\in (0,1)$ and $x_0\in D$ such that
\beq\label{*GB}\<x_0-x, \n(x)\> \ge \vv |x-x_0|,\ \ x\in\pp D.\end{equation}
 \end{enumerate}
 Indeed, if (a) holds then there exists $r_0>0$ such that $\rr\in C_b^2(\pp_{r_0} D)$. Let $h\in C^\infty([0,\infty))$ with $h(r)=r$ for $r\in [0, r_0/4]$ and $h(r)=r_0/2$ for $r\ge r_0/2.$ By taking $\tt\rr=h\circ\rr$ we have $\tt\rr\in C_b^2(\bar D)$, $\<\nn \tt\rr,\n\>|_{\pp D}=1$, and for any $x\in  D$ letting $\bar x\in \pp D$ such that $|x-\bar x|=\rr(x),$ we deduce from \eqref{LIP} that
\beg{align*}  \<b_t^{(1)}(x), \nn \tt\rr(x)\>= h'(\rr(x)) \big\{\<b_t^{(1)} (\bar x), \n(\bar x)\>+ \<b_t^{(1)}(x)-b_t^{(1)}(\bar x),\n(\bar x)\>\big\} \ge
  -(1+r_0)L \|h'\|_\infty.\end{align*}
Therefore,  \eqref{GB1} holds for some  (different) constant $L$.
Next, if (b) holds, by \eqref{*GB} we may take $\tt\rr(x)=N\ss{1+|x-x_0|^2} $ for large enough $N\ge 1$ such that $\<\nn \tt\rr, \n\>|_{\pp D}\ge 1$. So, by the boundedness of $D$ and $b^{(1)}\in C([0,T]\times\R^d),$
\eqref{GB1} holds for some constant $L>0$.

\

Assumption   $(A^{\si,b}_0)$    will be used to establish   Krylov's estimate for functions $f\in \cap_{(p,q)\in\scr K}\tt L_q^p(T)$, which is crucial to solve singular SDEs, see Lemma \ref{L1} below.
To improve this estimate for $(p,q)$ satisfying   $\ff {d}p+\ff 2 q<2$ as in the case without reflecting (see \cite{XXZZ}),
we introduce one  more assumption.

Consider the following differential operators on $\bar D$:
\beq\label{BARL}    L_t^{\si,b^{(1)}} := \ff 1 2  {\rm tr}\big( \si_t\si_t^* \nn^2\big) +     \nn_{b_t^{(1)}},\ \ t\in [0,T].\end{equation}
Let
$\{P_{s,t}^{\si, b^{(1)}}\}_{T\ge t_1\ge t\ge s\ge 0}$ be the Neumann semigroup on $\bar D$ generated by
$L_t^{\si,b^{(1)}},$   that is, for any $\phi\in C_b^2(\bar D)$,  and any  $t\in (0,T]$, $(P_{s,t}^{\si,b^{(1)}}\phi)_{s\in [0,t]}$ is the unique solution of the PDE
\beq\label{NMM} \pp_s u_s= - L_s^{\si,b^{(1)}} u_s,\ \ \nn_{\n}u_s|_{\pp D}=0\ \text{for}\ s\in [0,t), u_t=\phi.\end{equation}
For any $t>0,$ let $C_b^{1,2}([0,t]\times\bar D)$ be the set of  functions  $f\in C_b([0,t]\times\bar D)$  with bounded and continuous derivatives $\pp_t f, \nn f$ and $\nn^2 f$.

\emph{\beg{enumerate}  \item[$(A_2^{\si,b})$]    $\pp D\in C_b^{2,L}$ and the following conditions hold for $\si$ and $b$   on $[0,T]\times \bar D$:
\item[$(1)$]  $a_t:=\si_t\si_t^*$ is invertible,    $\eqref{*UN}$ holds for $x,y\in \bar D$ and there exist $\{(p_i,q_i)\}_{0\le i\le l} \subset \scr K$  with $p_i>2$ and $0\le f_i\in\tt L_{q_i}^{p_i}(T), 1\le i\le l,$ such that
 $$  \|\nn\si\|\le \sum_{i=1}^l f_i,\ \ \|a\|_\infty+\|a^{-1}\|_\infty + \|\nn \si\|_{\tt L_{q_1}^{p_1}(T,D)}<\infty. $$
\item[$(2)$]  $b= b^{(1)}+b^{(0)}$ with $\nn_\n b_t^{(1)}|_{\pp D}=0$,  $\|\nn b^{(1)}\|_\infty+\|1_{\pp D} \<b^{(1)}, \n\>\|_\infty<\infty$ and
$|b^{(0)}|\in \tt L_{q_0}^{p_0}(T,D)$ for some $(p_0,q_0)\in \scr K$ with $p_0>2$.
\item[$(3)$]  For any  $\phi\in C_b^2(\bar D)$ and $t\in (0,T]$, the PDE \eqref{NMM} has a unique solution $P_{\cdot,t}^{\si,b^{(1)}}\phi\in C_b^{1,2}([0,t]\times \bar D),$  such that    for some  constant $c>0$  we have
 \beq\label{AB1} \|\nn^i P_{s,t}^{\si,b^{(1)}} \phi\|_{\infty} \le  c (t-s)^{-\ff 1 2} \|\nn^{i-1} \phi\|_{\infty},\ \ 0\le s<t\le T, \ i=1,2,\phi\in C_b^2(\bar D),\end{equation}  where $\nn^{0}\phi:=\phi.$
 \end{enumerate}}

\paragraph{Remark 2.2.}  (1) Let $\rr\in C_b^2(\pp_{r_0}D)$ for some $r_0>0$. Since $\nn\rr|_{\pp D}=\n$,  $\|\nn b^{(1)}\|_\infty+\|1_{\pp D} \<b^{(1)}, \n\>\|_\infty<\infty$ implies
$\|1_{\pp_{r_0}D} \<b^{(1)},\nn\rr\>\|_\infty<\infty,$ which will be used in the proof of Lemma \ref{LNN} below.

(2) $(A_2^{\si,b})(3)$ holds   if  $D$ is  bounded with $\pp D\in C^{2+\aa}$ for some $\aa\in (0,1),$  and    there exists   $c>0$ such that  
\beq\label{**A}  \big\{|b_t^{(1)}(x)-b_s^{(1)}(y)|+\|a_t(x)-a_s(y)\|\big\}\le c (|t-s|^{\aa}+|x-y|^{\ff \aa 2}),\ \ s,t\in [0,T],x,y\in \bar D.\end{equation}
Indeed,    $\pp D\in C^{2+\aa}$ implies $\n\in C^{1+\aa}(\pp D)$,  so that  $\eqref{**A}$ implies  estimates  (3.4) and (3.6) in \cite[Theorem VI.3.1]{Carr}  with $\varrho=\infty$   for  the Neumann heat kernel $p_{s,t}^{\si,b^{(1)}}(x,y)$ of $P_{s,t}^{\si,b^{(1)}}$.  We note that according to its proof, the condition (3.3) therein is assumed for some $\aa\in (0,1)$ rather than all $\aa\in (0,1).$   In particular,   $\nn^2 p_{s,t}^{\si,b^{(1)}}  (\cdot, y)(x)$ and $\pp_s p_{s,t}^{\si,b^{(1)}}(x,y)$ are   continuous in $(s,x)\in [0,t]\times\bar D$, and   there exists a  constant $c>1$ such that
  \beg{align*} &|\nn^i p_{s,t}^{\si,b^{(1)}}(\cdot,y)(x)|\le c|t-s|^{-\ff{d+i}2} \e^{-\ff{|x-y|^2}{c(t-s)}},\ \ 0\le s<t\le T, x,y\in \bar D, i=0,1,2,\\
&|\pp_s p_{s,t}^{\si,b^{(1)}}(x,y)|=|L_s^{\si,b^{(1)}} p_{s,t}^{\si,b^{(1)}}(\cdot,y)(x)|\le c|t-s|^{-\ff{d+2}2} \e^{-\ff{|x-y|^2}{c(t-s)}},\ \ 0\le s<t\le T, x,y\in \bar D.\end{align*}
   These properties   imply  $\eqref{AB1}.$

\

The following are   main results of this section,  where Theorem \ref{T2.2}  improves the main result (Theorem 6.3) in \cite{YZ} for bounded  $C^3$ domain $D$, bounded drift $b$
  and  $\si= \I_d.$ Moreover, going back to the case without reflection (i.e. $D=\R^d$),  Theorem \ref{T2.3} covers the main result (Theorem 1.1) of \cite{LLW} where  $b^{(1)}=0$ is considered.

\beg{thm}[Weak well-posedness]\label{T2.1}  If  either $(A_1^{\si,b})$ or $(A_2^{\si,b})$ holds,    then $\eqref{E01}$ is  weakly well-posed. Moreover, for any $k\ge 1$ there exists a constant $c>0$ such that
  \beq\label{EPP} \E \Big[\sup_{t\in [0,T]}  |X_t^x|^k\Big] \le c(1+|x|^k),\ \ \E \e^{k l_T^x}\le c,\ \ x\in \bar D,\end{equation}
  where $(X_t^x,l_t^x)$ is the  $($weak$)$ solution   of $\eqref{E01} $  with  $X_0^x=x$.
  \end{thm}

\beg{thm}[Well-posedness]\label{T2.2} Assume that one of the following conditions holds:
 \beg{enumerate} \item[$(i)$]  $d=1$ and $(A_1^{\si,b})$ holds;
 \item[$(ii)$]   $(A_2^{\si,b})$ holds with $p_1>2$. \end{enumerate}
 Then  $\eqref{E01}$ is well-posed, and   for any $k\ge 1$, there exists a constant $c>0$ such that
\beq\label{INN} \E \Big[\sup_{t\in [0,T]}|X_t^x-X_t^y|^k\Big]\le c|x-y|^k,\ \   x,y\in \bar D.\end{equation}
Consequently, for any $p>1$   there exists a constant $c(p)>0$ such that
$P_tf(x):= \E[ f(X_t^x)]$ satisfies
\beq\label{GRD} |\nn P_t f|\le c(p)(P_t|\nn f|^p)^{\ff 1 p},\ \ f\in C_b^1(\bar D),\ \ t\in [0,T].  \end{equation}
\end{thm}

\beg{thm}[Functional inequalities] \label{T2.3}  Assume  that   $(A_2^{\si,b})$ holds with $p_1>2$. Then    there exist  a constant $C>0$ and a map $c: (1,\infty)\to (0,\infty)$ such that
\beq\label{GRD'} |\nn P_t f|\le \ff{c(p)}{\ss {t}} (P_t|f|^p)^{\ff 1 p},\ \ t\in [0,T], f\in \B_b(\bar D),\ p>1,\end{equation}
\beq\label{PCC} P_tf^2-(P_tf)^2\le tC  P_t|\nn f|^2,\ \ f\in C_b^1(\bar D),\ \ t\in [0,T],\end{equation}
\beq\label{LHA} P_t \log f(x)\le \log P_tf(y) + \ff{C|x-y|^2}{  t}, \ \ t\in [0,T], x,y\in \bar D, 0< f \in \B_b(\bar D).\end{equation}
  \end{thm}

 To prove these results,  we   first   establish Krylov's estimates under different conditions, then  prove   the weak and strong well-posedness by using Girsanov's transform and  Zvokin's transforms respectively.

 \subsection{Krylov's estimate and  It\^o's formula}

A crucial step in the study of  singular SDEs   is to establish Krylov's estimate \cite{Kry}.
 To this end, we first introduce the following lemma   taken from \cite[Theorem 2.1]{YZ0},
which extends \cite[Theorem 3.2]{XXZZ} where $b^{(1)}=0$ is considered. See \cite{XXZZ,Z2} and references within for earlier assertions.

\beg{lem}\label{LYZ} Assume $(A_0^{\si,b})$. For any $0\le t_0<t_1\le T$ and $f\in \tt L_{q}^p(t_0,t_1)$ for some $p,q>1,$  the   PDE
\beq\label{PPD5} (\pp_t +L_t)u_t^\ll=  \ll u_t^\ll+ f_t,\ \ t\in [t_0,t_1], u^\ll_{t_1}=0,\end{equation}
  has a unique solution in $\tt H_q^{2,p}(t_0,t_1)$. Moreover,    for any $\theta\in [0,2), p'\in [p,\infty]$ and $q'\in [q,\infty]$ with $\ff{d}p +\ff 2 q< 2-\theta +\ff d{p'}+\ff 2 {q'},$ there exist constants $\ll_0, c>0$ increasing in $\|b^{(0)}\|_{\tt L_{q_0}^{p_0}(T)}\,($i.e. they do not have to be changed when $b^{(0)}$ is replaced by $\tt b^{(0)}$ with $\|\tt b^{(0)}\|_{\tt L_{q_0}^{p_0}(T)}\le \|b^{(0)}\|_{\tt L_{q_0}^{p_0}(T)})$, such that for any $\ll\ge \ll_0$ and 
  $0\le t_0<t_1\le T, \ll\ge \ll_0$ and $f\in \tt L_{q}^p(t_0,t_1)$, the solution satisfies 
$$\ll^{\ff 1 2 (2-\theta+\ff d {p'}+\ff 1 {q'} -\ff d p-\ff 2 q)} \|u^\ll\|_{\tt H_{q'}^{\theta,p'}(t_0,t_1)}+ \|(\pp_t+\nn_{b^{(1)}})u^\ll\|_{\tt L_q^p(t_0,t_1)}+ \|u^\ll\|_{\tt H_q^{2,p}(t_0,t_1)}\le c \|f\|_{\tt L_q^p(t_0,t_1)}.$$
\end{lem}

By estimating the local time,     this result enables us to derive
the following   Krylov's estimate \eqref{KR} and Khasminskii's estimate \eqref{KAS}.

 \beg{lem} \label{L1}  Assume   $(A^{\si,b}_0)$. Let $(p,q)\in \scr K$.
  \beg{enumerate} \item[$(1)$] There exist a constant $i\ge 1$ depending only on $(p,q)$, and  a constant $c\ge 1$ increasing in $\|b^{(0)}\|_{\tt L_{q_0}^{p_0}(T)}$, such that  for    any   solution $X_t$ of $\eqref{E01}$, and any $0\le t_0\le t_1\le T,$
  the following estimates hold.
\beq\label{KR}  \E \bigg[\bigg(\int_{t_0 }^{t_1 } |f_s(X_s)|\d s\bigg)^m\bigg|\F_{t_0}\bigg]
 \le c^m m! \|f\|_{\tt L_{q}^{p}(t_0,t_1)}^m,\ \  f\in \tt L_q^p(t_0,t_1), m\ge 1,  \end{equation}
\beq\label{KAS}  \E  \big(\e^{ \int_{t_0}^{t_1} |f_t(X_t)|\d t}\big|\F_{t_0}\big)  \le \exp\big[c + c \|f\|_{\tt L_{q}^{p}(t_0,t_1)}^i\big],\ \      f\in \tt L_{q}^{p}(t_0,t_1),\end{equation}
\beq\label{KAS'} \sup_{t_0\in [0,T]}\E  \big(\e^{\ll (l_T-l_{t_0})}\big|\F_{t_0}\big)<\e^{c(1+\ll^2)},\ \ \ll>0.\end{equation}
\item[$(2)$]    For  any   $u\in C([0,T]\times \R^d)$ with continuous $\nn u$ and
 \beq\label{UUU} \|u\|_\infty+\|\nn u\|_\infty+ \|(\pp_t +\nn_{b^{(1)}}) u\|_{ \tt L_{q}^{p}(T)}+\|\nn^2 u\|_{\tt L_q^p(T)} <\infty,\end{equation}
 we have the following It\^o's formula for a solution $X_t$ to $\eqref{E01}$:
\beq\label{ITO}  \d   u_t (X_t)= (\pp_t+L_t)u_t (X_t) \d t + \<\nn u_t(X_t), \si_t(X_t)\d W_t\>
 + (\nn_\n   u_t)(X_t)\d l_t. \end{equation}
 \end{enumerate} \end{lem}

\beg{proof} (1) We first prove \eqref{KR} for $m=1$.
By first using $(|f|\land n)1_{B(0,n)}$ replacing $f$  then letting $n\to \infty$,    we may and do assume that $f$ is bounded   with compact support.   Next, by a standard approximation argument, we only need to prove for   $f\in C_0^\infty([t_0,t_1]\times \R^d)$.

Let $f\in C_0^\infty([t_0,t_1]\times \R^d)$.  By Lemma \ref{LYZ}, for any $(p',q')\in \scr K$, \eqref{PPD5} has a unique solution satisfying
\beq\label{PPD5'}  \beg{split}&\ll^{\vv}\big( \| u^\ll\|_\infty  +\|\nn u^\ll\|_\infty\big)+ \|(\pp_t+\nn_{b^{(1)}} )u^\ll\|_{\tt L_{q'}^{p'}(t_0,t_1)}   +\| u^\ll\|_{\tt H_{q'}^{2,p'}(t_0,t_1)}\\
&  \le c_1 \|   f\|_{\tt L_{q'}^{p'}(t_0,t_1)},\ \ \ll\ge \ll_0,\end{split}\end{equation}
where $\vv>0$ depends on $(p',q')$ and  $\ll_0,c>0$ are constants increasing in $\|b^{(0)}\|_{\tt L_{q_0}^{p_0}(T)}.$ 
To apply It\^o's formula, we make a standard mollifying approximation  of $u^\ll$, which is extended to $\R^{d+1}$ by letting $u_t^\ll:=u^\ll_{(t\lor t_0)\land t_1}$ for $t\in \R.$    Let  $0\le \varrho\in C_0^\infty(\R^{d+1})$ such that $\int_{\R^{d+1} } \varrho(z)\d z =1.$ For any $n\ge 1$, let
\beq\label{MLL} u^{\ll,n}_t(x) = n^{d+1} \int_{\R^{d+1}} u^\ll_{t-s}(x-y) \varrho(ns, ny)\d s\d y,\ \ t\in\R, x\in \R^d.\end{equation}
Then
$$\lim_{n\to\infty}  \big\{\|(\pp_t+\nn_{b^{(1)}}) (u^{\ll,n}-u^\ll)\|_{\tt L_{q'}^{p'}(t_0,t_1)} + \|u^{\ll,n}-u^\ll\|_{\tt H_{q'}^{2,p'}(t_0,t_1)} \big\}=0,\ \ (p',q')\in \scr K,$$
so that as shown in the proof of \cite[Lemma 5.4]{XZ}, 
\beq\label{P*D} f^{\{n\}}_t:= (\pp_t + L_t-\ll) u_t^{\ll,n} \end{equation}
 satisfies
\beq\label{FNN} \lim_{n\to\infty} \|f-f^{\{n\}}\|_{\tt L_{q'}^{p'}(t_0,t_1)}=0,\ \ (p',q')\in \scr K,\end{equation} and
\eqref{PPD5'} with $(p',q')=(p,q)$  implies
\beq\label{PPD6}  \| u^{\ll,n}\|_\infty  +\|\nn u^{\ll,n}\|_\infty  \le c \ll^{-\vv}  \|f\|_{\tt L_q^p(t_0,t_1)},\ \ n\ge 1,\ll>\ll_0.\end{equation}
 By Theorem 6.2.7(ii)-(iii) in \cite{BKRS},
the conditional distribution of $X_t$ under $P_{t_0} $   is absolutely continuous for $t>t_0$,  so that by the dominated  convergence theorem, \eqref{FNN} implies $\P$-a.s. 
\beq\label{*NN} \E\bigg(\int_{t_0}^{t_1\land \tau_k} f_s(X_s)\d s\bigg|\F_{t_0} \bigg)= \lim_{n\to\infty} \E\bigg(\int_{t_0}^{t_1\land \tau_k} f_s^{\{n\}}(X_s)\d s\bigg|\F_{t_0}\bigg).\end{equation} 
Let
 $$\tau_k:=\inf\bigg\{t\in [t_0,T]: l_t-l_{t_0}+\int_{t_0}^t |b_s(X_s)|\d s\ge k\bigg\},\ \ k\ge 1.$$
  Applying It\^o's formula to $u^{\ll,n}$, we deduce from  \eqref{P*D} and   \eqref{PPD6} that 
\beq\label{TTP} \beg{split} &2c\ll^{-\vv}\|f\|_{\tt L_q^p(t_0,t_1)}\ge \E\big\{ u^{\ll,n}_{t_1\land\tau_k}(X_{t_1\land\tau_k})- u^{\ll,n}_{t_0}(X_{t_0}) \big|\F_{t_0}\big\}\\
&=\E\bigg(\int_{t_0 }^{t_1\land\tau_k} (\pp_s+L_s )u^{\ll,n}_s(X_s)\d s + \int_{t_0 }^{t_1\land\tau_k} \{\nn_{\n(X_s)} u^{\ll,n}_s\}(X_s)\d l_s \bigg|\F_{t_0}\bigg)\\
&\ge \E\bigg(\int_{t_0 }^{t_1\land\tau_k} f^{\{n\}}_s (X_s)\d s \bigg|\F_{t_0}\bigg)-c  \|f\|_{\tt L_q^p(t_0,t_1)} \big\{\ll+ \ll^{-\vv}\E(l_{t_1\land\tau_k} -l_{t_0 }|\F_{t_0})\big\}. \end{split}\end{equation}
  Therefore,
\beq\label{KR1} \beg{split} &\E\bigg(\int_{t_0}^{t_1\land\tau_k} f^{\{n\}}_s(X_s)\d s\bigg|\F_{t_0}\bigg)\\
&\le c\|f\|_{\tt L_q^p(t_0,t_1)}\big\{2+\ll + \ll^{-\vv} \E(l_{t_1\land\tau_k} -l_{t_0 }|\F_{t_0})\big\},\ \ n,k\ge 1,\ll>0.\end{split}\end{equation}
Combining this with \eqref{*NN}, we obtain   
\beq\label{KR1'} \beg{split} &\E\bigg(\int_{t_0 }^{t_1\land \tau_k }   f_s(X_s)\d s\bigg|\F_{t_0}\bigg)= \lim_{n\to\infty} \E\bigg(\int_{t_0 }^{t_1\land \tau_k }   f_s^{\{n\}}(X_s)\d s\bigg|\F_{t_0}\bigg)\\
&\le c\|f\|_{\tt L_q^p(t_0,t_1)}\big\{2+\ll + \ll^{-\vv} \E(l_{t_1 } -l_{t_0 }|\F_{t_0})\big\},\ \ \ll>0, k\ge 1.\end{split}\end{equation}

On the other hand, by \eqref{GB1} and the boundedness of $\si$, we find a constant $c_1>0$ such that
\beq\label{GHH} \d \tt\rr(X_t)\ge -c_1\d t- c_1 |b_t^{(0)}(X_t)|\d t +\d l_t +\<\nn \tt\rr(X_t), \si_t(X_t)\d W_t\>.\end{equation}
So, \eqref{KR1'} with $(p,q)=(p_0,q_0)$  implies 
\beg{align*}& \E(l_{t_1 \land\tau_k}-l_{t_0}|\F_{t_0})\le c_1(t-t_0) + c_1 \E\bigg(\int_{t_0}^{t_1\land\tau_k}|b_s^{(0)}(X_s)|\d s\bigg|\F_{t_0}\bigg)+\|\tt\rr\|_\infty\\
&\le c_2(1+\ll) + c_2\ll^{-\vv} \E(l_{t_1\land\tau_k}-l_{t_0}|\F_{t_0}),\ \ t\in [t_0,T],\ \ \ll>0, k\ge 1\end{align*}   for some constant $c_2>0$  increasing in $\|b^{(0)}\|_{\tt L_q^p(T)}$.
 Taking $\ll>0$ large enough such that $c_2\ll^{-\vv}\le \ff 1 2$,
we arrive at
$$\E(l_{t_1\land\tau_k}-l_{t_0}|\F_{t_0})\le c_3,\ \ k\ge 1$$
for some constant $c_3>0$  increasing in $\|b^{(0)}\|_{\tt L_q^p(T)}$. Letting $k\to\infty$  gives
\beq\label{LCT} \E(l_{t_1}- l_{t_0}|\F_{t_0})\le c_3,\ \ t_0\le t_1\le T.\end{equation}
This and \eqref{KR1'} with $k\to\infty$  imply \eqref{KR} for $m=1$, which further yields the inequality for any $m\ge 1$ as shown in the proof of \cite[Lemma 3.5]{XZ}. 
Moreover, taking $q'\in (2,q)$ such that $(p,q')\in \scr K$,  \eqref{KR} for $m=1$ with $(p,q')$ replacing $(p,q)$ yields
 $$\E\bigg(\int_{t_0}^{t_1}f_s(X_s)\d s\bigg|\F_{t_0}\bigg)\le c \|f\|_{\tt L_{q'}^p(t_0,t_1)}\le c(t_1-t_0)^{\ff{q-q'}{qq'}}\|f\|_{\tt L_q^p(t_0,t_1)}.$$
 This and \cite[Lemma 3.5]{XZ} with $\tt L_{q'}^p$ replacing $L_q^p$ imply \eqref{KAS} for $i=\ff q {q-q'}.$
 Finally, combining \eqref{KAS} with \eqref{GHH}, $b^{(0)}\in \tt L_{q_0}^{p_0}(T)$ and $\|\si^*\nn \tt\rr\|_\infty<\infty$, we derive \eqref{KAS'}.

 (2)   We first extend $u$ to $\R^{d+1}$ by letting
 $u_t=  u_{t^+\land T} $ for $t\in \R$, and consider its mollifying approximation $u^{\{n\}}$ defined above. Then $\|\si\|_\infty<\infty$ and \eqref{UUU} imply
 \beq\label{XN0} \lim_{n\to\infty}\big\{\|u-u^{\{n\}}\|_{\infty}+ \|\nn (u-u^{\{n\}})\|_\infty+ \|(\pp_t+ L_t) (u-u^{\{n\}})\|_{\tt L_q^{p}(T)}\big\}=0.\end{equation}
 Combining this with $\|\si\|_\infty<\infty$ and \eqref{KR}, we obtain
\beq\label{XN2} \beg{split} & \lim_{n\to\infty} \sup_{t\in [0,T]} |u_t^{\{n\}}(X_t)- u_t(X_t)|= 0,\ \P\text{-a.s.}\\
& \lim_{n\to\infty} \int_0^t \nn_{\n} u_s^{\{n\}}(X_s)\d l_s =  \int_0^t \nn_{\n} u_s(X_s)\d l_s,\ \P\text{-a.s.}\\
 &\lim_{n\to\infty} \E    \int_0^{T}\big|(\pp_s+ L_s)(u_s^{\{n\}}-u_s)\big| (X_s) \d s=0,\\
 &\lim_{n\to\infty}\E \sup_{t\in [0,T]} \bigg|\int_0^t\<\nn(u_s^{\{n\}}-u_s)(X_s), \si_s(X_s)\d W_s\>\bigg|=0.\end{split}\end{equation}
Therefore, we prove \eqref{ITO} by letting $n\to\infty$ in the following It\^o's formula:
 \beg{align*} u^{\{n\}}_t(X_t)= &\,u_0^{\{n\}}(X_0)+\int_0^t (\pp_s+L_s)(u_s^{\{n\}})(X_s)\d s\\
  &+ \int_0^t \<\nn u_s^{\{n\}}(X_s), \si_s(X_s)\d W_s\>+\int_0^t  (\nn_\n  u_s^{\{n\}})(X_s)\d l_s,\ \ t\in [0,T].\end{align*}
  \end{proof}


To improve Lemma \ref{L1} for $(p,q)\in \scr K$ with $\ff d p+\ff 2 q<2$, we first  extend  Lemma \ref{LYZ}  to the Neumann boundary case.
For any $k\in\mathbb N$,  let $C_b^{0,k}([t_0,t_1]\times\bar D; \R^d)$ be the space of     $f\in C_b([t_0,t_1]\times \bar D; \R^d)$ with bounded and continuous derivatives in $x\in \bar D$ up to order $k$.  Let $ C_b^{1,2}([t_0,t_1]\times\bar D; \R^d)$ denote the space of  $f\in C_b^{0,2}([t_0,t_1]\times\bar D; \R^d)$ with bounded and continuous $\pp_t f$.  

\beg{lem}\label{LNN} Assume $(A_2^{\si,b})$  but without  the condition  on $\|\nn\si \|$.  Then $(A_0^{\si,b})$ and the following assertions hold.
 \beg{enumerate} \item[$(1)$] For any $\ll\ge 0$,  $0\le t_0<t_1\le T$  and  $ \tt b, f\in C_b^{0,2}([t_0,t_1]\times\bar D; \R^d)$,
the PDE
\beq\label{PDEE} (\pp_t+ L_t^{\si,b^{(1)}}+\nn_{\tt b_t}-\ll)\tt u_t^{\ll} = f_t,\ \ \tt u_{t_1}^\ll=\nn_\n \tt u_t^\ll|_{\pp D} =0, t\in [t_0,t_1]\end{equation}
has a unique solution    $\tt u^\ll\in C_b^{1,2}([t_0,t_1]\times\bar D; \R^d).$
\item[$(2)$] For any $ (p,q),(p',q')\in \scr K$ and $\tt b\in  C_b^{0,2}([0,T]\times\bar D; \R^d)$, there exist a constant $\vv>0$ depending only on $(p,q)$ and $(p',q')$, and  constants $\ll_0,c >0$ increasing in $\|\tt b\|_{\tt L_{q'}^{p'}(T,D)},$  such that  for any $0\le t_0<t_1\le T$ and $f\in  C_b^{0,2}([t_0,t_1]\times\bar D; \R^d),$
\beq\label{KG}  \ll^{\vv}   (\|\tt u^\ll\|_{\infty}  +   \|\nn \tt u^\ll\|_{\tt L_{q}^{p}(t_0,t_1,D)})
 \le c  \|f\|_{\tt L_{q/2}^{p/2}(t_0,t_1,D)},\ \  \ll\ge \ll_0 \ ({\rm when}\ p>2),\end{equation}
\beq\label{KG'}  \ll^{\vv}   \|\nn \tt u^\ll\|_{\infty}
 \le c  \|f\|_{\tt L_{q}^{p}(t_0,t_1,D)}, \ \ \ \ll\ge \ll_0, \end{equation}
 and  there exists decomposition $\tt u^\ll= \tt u^{\ll,1}+\tt u^{\ll,2}$ such that
\beq\label{J2} \beg{split}  &\|\nn^2\tt u^{\ll,1}\|_{\tt L_q^p(t_0,t_1,D)}
 + \|(\pp_t+\nn_{b^{(1)}})\tt u^{\ll,1}\|_{\tt L_q^p(t_0,t_1,D)} +  \|\nn^2\tt u^{\ll,2}\|_{\tt L_{q'}^{p'}(t_0,t_1,D)}
\\
&+ \|(\pp_t+\nn_{b^{(1)}})\tt u^{\ll,2}\|_{\tt L_{q'}^{p'}(t_0,t_1,D)}  \le c  \|f\|_{\tt L_q^p(t_0,t_1,D)},\ \ \ll\ge\ll_0.\end{split}\end{equation}\end{enumerate}
\end{lem}

\beg{proof}  (1) Let $\V:=   C_b^{0,2}([t_0,t_1]\times \bar D;\R^d)$,  which   is a Banach space under the norm
$$\|u\|_{\V,N}:=\sup_{t\in [t_0,t_1]} \e^{-N(t_1-t)} \big\{ \|u_t\|_\infty+ \|\nn u_t\|_\infty+ \|\nn^2 u_t\|_\infty\big\},\ \ u\in\V$$
for $N>0$. To solve \eqref{PDEE}, for any $\ll\ge 0$ and $u\in \V$, let
$$\Phi^\ll_s(u):= \int_s^{t_1} \e^{-\ll(t-s)} P_{s,t}^{\si,b^{(1)}}\{ \nn_{\tt b_t}u_t-f_t\}\d t,\ \ s\in [t_0,t_1].$$
Then  $(A_2^{\si,b})$ implies   $\Phi^\ll(u)\in  C_b^{1,2}([t_0,t_1]\times \bar D)$ with
\beq\label{NMM'} (\pp_s + L_s^{\si,b^{(1)}} -\ll)\Phi_s^\ll(u) = f_s-\nn_{\tt b_s} u_s,\ \ s\in [t_0,t_1], \nn_{\n} \Phi_t^\ll(u)|_{\pp D}=0, \Phi_{t_1}^\ll(u)=0.\end{equation}
 So, it suffices to prove that $\Phi^\ll$ has a unique fixed point $\tt u^\ll\in \V:$
\beq\label{KG0} \tt u_s^\ll=   \int_s^{t_1} \e^{-\ll(t-s)}P_{s,t}^{\si,b^{(1)}}\big\{\nn_{\tt b_t}\tt u_t^\ll-f_t\big\}\d t,\ \ s\in [t_0,t_1],\end{equation}
which,   according to   \eqref{NMM'}, is  the unique solution of  \eqref{PDEE} in $C_b^{1,2}([t_0,t_1]\times\bar D;\R^d)$.

For any $u,\bar u\in \V$, by $\|\tt b\|_\infty<\infty$, we   find a constant $c_1>0$ such that
$$\|\Phi_s^\ll(u)-\Phi_s^\ll(\bar u)\|_\infty \le \int_s^{t_1} \|\tt b_t\|_\infty \|\nn(u_t-\bar u_t)\|_\infty \d t \le  c_1\int_s^{t_1} \|\nn(u_t-\bar u_t)\|_\infty \d t.$$
Similarly,   \eqref{AB1} with $i=1$ implies
\beg{align*}&\|\nn\{\Phi^\ll(u)_s-\Phi^\ll(\bar u)_s\}\|_\infty \le c  \int_s^{t_1} (t-s)^{-\ff 1 2} \|\tt b_t\|_\infty  \|\nn(u_t-\bar u_t)\|_\infty \d t \\
&\le  c_1\int_s^t(t-s)^{-\ff 1 2} \|\nn(u_t-\bar u_t)\|_\infty \d t,\end{align*}
while \eqref{AB1} with $i=2$ and $\|\tt b\|_\infty+\|\nn \tt b_t\|_\infty<\infty$ yield
\beg{align*} &\|\nn^2\{\Phi_s^\ll(u)-\Phi_s^\ll(\bar u)\}\|_\infty \le c  \int_s^{t_1} (t-s)^{-\ff 1 2} \big\|\nn\{\nn_{\tt b_t}(u_t-\bar u_t)\}\big\|_\infty \d t\\
&\le c_1 \int_s^{t_1} (t-s)^{-\ff 1 2} \big\{\|\nn (u_t-\bar u_t)\|_\infty+ \|\nn^2(u_t-\bar u_t)\|_\infty \big\}   \d t.\end{align*}
Combining these with \eqref{NMM'} and the boundedness of $a$ and $\tt b\in C_b^{0,1}([t_0,t_1]\times\bar D;\R^d)$, we find a constant $c_2>0$ such that
\beg{align*} &\|\Phi^\ll(u)-\Phi^\ll(\bar u)\|_{\V,N} \\
&\le c_2 \sup_{s\in [t_0,t_1]}  \int_s^{t_1} \e^{-N(t_1-s)}   (t-s)^{-\ff 1 2}\Big\{\|u_t-\bar u_t\|_\infty \\
&\qquad\qquad +\|\nn(u_t-\bar u_t)\|_\infty+ \|\nn^2 (u_t-\bar u_t)\|_\infty\Big\} \d t\\
&\le c_2 \|u-\bar u\|_{\V,N} \sup_{s\in [t_0,t_1]} \int_s^{t_1} \e^{-N(t-s)}(t-s)^{-\ff 1 2}\d t.\end{align*}
So, $\Phi^\ll$ is contractive under the norm $\|\cdot\|_{\V,N}$ for large enough $N>0$, and hence has a unique fixed point $\tt u^\ll$ in $\V$.

(2) To prove \eqref{KG} and \eqref{J2}, we extend the PDE \eqref{PDEE} to a global one such that   estimates in Lemma \ref{LYZ} apply.
By $(A_2^{\si,b})$,  there exists  $r_0>0$     such that
$$\varphi: \pp_{-r_0 D}\to \pp_{r_0}D; \ \  \theta- r\n(\theta)\mapsto \theta + r \n(\theta),\ \ r\in [0,r_0], \theta\in \pp D$$
is a $C_b^{1,L}$-diffeomorphism (i.e. it is a homeomorphism with $\nn \varphi$ bounded and Lipschitz continuous) and $\rr_D:={\rm dist}(\cdot, D)\in C^2_b(D_{r_0}),$ recall that $D_{r_0}  =\{\rr_D\le r_0\}.$    For any vector field $v$ on $\pp_{r_0}D$,
  $v^\star:=(\varphi^{-1})^*v$ is the vector  field on $\pp_{-r_0}^0D:=\pp_{-r_0}D\setminus \pp D$ given by
$$\<v^\star, \nn g\>(x):= \<v, \nn (g\circ \varphi^{-1})\>(\varphi(x)),\ \ x\in \pp_{-r_0}^0D,\ g\in C^1(\pp_{-r_0}^0D).$$
We then extend $b_t^{(1)}$ and $\tt b_t$ to $\R^d$ by taking
\beq\label{BAB}   b_t^{(1)}:= 1_{\bar D} b_t^{(1)}+ h(\rr_D/2)1_{\pp_{-r_0}^0D} (b_t^{(1)})^\star, \ \   \tt b_t:= 1_{\bar D} \tt b_t+ 1_{\pp_{-r_0}^0D} (\tt b_t)^\star,\end{equation}
where   $h\in C^\infty(\R)$ such that  $0\le h\le 1, h|_{(-\infty,r_0/4]}=1$ and $h|_{[r_0/2,\infty)}=0.$ Since $(A_2^{\si,b})$ implies
  $\|1_{\bar D}\nn b^{(1)}\|_\infty<\infty$   and $\nn_\n b^{(1)}|_{\pp D}=0$, we have $\|\nn b^{(1)}\|_\infty<\infty.$
Let
\beq\label{SLL} \tt\varphi(x):= x 1_{\bar D}(x) +  \varphi(x) 1_{\pp_{-r_0}^0 D}(x),\ \  x\in  D_{r_0}.   \end{equation}
We  extend $\tt u^\ll$ to $[t_0,t_1]\times \R^d$ by  setting
\beq\label{UUK} u_t^\ll= h(\rr_D) (\tt u_t^\ll\circ\tt\varphi),\ \ t\in [t_0,t_1].\end{equation}
  We claim  that
\beq\label{NB1}u_t^\ll\in    C^{1,L}_{b}(\R^d),\ \ t\in [t_0,t_1],\end{equation}
  where $C^{1,L}_{b}(D_{r_0})$ is the class of $C_b^1$-functions $f$ on $D_{r_0}$ with Lipschitz continuous $\nn f.$
Indeed,   since  $\varphi$ is a  $C_b^{1,L}$-diffeomorphism from $\pp_{-r_0}D$ to $\pp_{r_0}D$, $\tt\varphi\in C_b^{1,L}(D_{r_0}\setminus \pp D)$ with bounded and continuous first and second order derivatives,  which together with $\tt u_t^\ll\in C_b^2(\bar D)$
yields $u_t^\ll \in C_b^{1,L}(\R^d\setminus \pp D).$ So, we only need to verify that  $\tt u_t^\ll\circ\tt\varphi \in C^{1,L}_{b}(D_{r_0}).$
 To this end,  for any $x\in \pp_{-r_0} D$ and $v\in \R^d$, let
$$\pi_x v:= v- \<v, \n(\theta(x))\>\n(\theta(x))$$ be the projection of $v\in T_x\R^d$ to the tangent space of $\pp D$, recall that  $\theta(x)$ is the projection of $x$ to $\pp D$, i.e. $x=\theta(x)-\rr_D(x)\n(\theta(x))$ for $\rr_D(x):={\rm dist}(x,D).$
We have
\beq\label{NB20} \beg{split} &\nn_v \tt\varphi (x) =  \nn_{\<v, \n(\theta(x))\> \n(\theta(x))} \tt\varphi (x) +\nn_{\pi_x v} \tt\varphi (x) \\
&= 1_{\pp D}(x) |\<v,\n(\theta(x))\>|\n(\theta(x)) +\{1_D- 1_{\pp_{-r_0}^0D}\}(x) \<v, \n(\theta(x))\>\n(\theta(x)) \\
&\qquad +\pi_x v
+\rr_D(x) (\nn_{\pi_x v}\n)(\theta(x)).
 \end{split}\end{equation}
  Since $\tt u^\ll_t\in C_b^2(\bar D)$ with $\nn_{\n } \tt u_t^\ll|_{\pp D}=0$, \eqref{NB20} yields
\beq\label{*QR} \beg{split}  &\nn_v(\tt u_t^\ll \circ\tt\varphi)(x) =  (\nn_v \tt u_t^\ll)\circ \tt\varphi(x) \\
&- 21_{\pp_{-r_0}^0D}(x)\<v, \n(\theta(x))\>\cdot\<\n(\theta(x)), (\nn \tt u_t^\ll)\circ \tt\varphi(x)\>\\
& +\rr_D(x) \big(\nn_{(\nn_{\pi_x v}\n)(\theta(x))}\tt u_t^\ll)\circ\tt\varphi(x),
\ \ x\in  D_{r_0}.\end{split}\end{equation}
Combining this with $\nn \tt u_t^\ll\in C_b^1(\bar D),  \nn_\n \tt u_t^\ll|_{\pp D}=0$ and $\n,\nn \n$ are Lipschitz continuous  on $\pp_{-r_0}D)$ due to $\pp D\in C_b^{2,L}$,
we conclude that $\nn (\tt u_t^\ll\circ \tt\varphi)$ is Lipschitz continuous on $D_{r_0}.$

Next, we construct the PDE satisfied by $u^\ll$. By \eqref{NB20}, we see that $(\nn\tt\varphi)(\nn\tt\varphi)^*=Q$ holds on $D_{r_0}\setminus \pp D$, where  $Q$ is a $d\times d$ symmetric matrix valued function given by
\beg{align*} &\<Q(x)v_1,v_2\>:=   \<v_1,v_2\> + \rr_D(x)^2 \big\<(\nn_{\pi_x v_1}\n)(\theta(x)),\ (\nn_{\pi_x v_2}\n)(\theta(x))\big\>\\
&+ \rr_D(x) \Big\{\big\<v_1-2 1_{\pp_{-r_0}D}(x) \<v_1, \n(\theta(x))\>\n(\theta(x)),\ (\nn_{\pi_x v_2}\n)(\theta(x))\big\>\\
&\qquad+\<v_2-2 1_{\pp_{-r_0}D}(x) \<v_2, \n(\theta(x))\>\n(\theta(x)),\ (\nn_{\pi_x v_1}\n)(\theta(x))\big\>\Big\}\end{align*}
for $x\in D_{r_0}, v_1,v_2\in \R^d.$
Then by taking $ r_0   >0$ small enough,  on $D_{r_0}$ the matrix-valued functional  $Q$ is bounded, invertible, Lipchitz continuous, and symmetric with
$Q^{-1}(x) \ge \ff 1 2 {\bf I}_d$ for $ x \in D_{r_0}.$ 
We  extend  $a_t:=\ff 1 2 \si_t\si_t^*$ from $\bar D$  to $\R^d$ by letting
\beq\label{BAA}  a_t:=  h (\rr_D/2) ( a_t\circ \tt\varphi)Q^{-1}+(1-h (\rr_D/2) ) {\mathbf I}_d.\end{equation}
Since \eqref{*UN} holds for $x,y\in \bar D$, with this extension of $a$  it holds for all $x,y\in \R^d$.  Combining this with \eqref{BAB}, Remark 2.1(a) for the existence of $\tt\rr$,  and noting that  $b_t=b_t^{(1)}+ 1_{\bar D} b_t^{(0)}$ extends $b$ from $\bar D$ to $\R^d$,  we see that $(A_0^{\si,b})$ holds.

Since   $h(\rr_D/2),h(\rr_D)\in C_b^2(\R^d)$ with $h(\rr_D/2)=1$ on $\{h(\rr_D)\ne 0\}$,
and      since   $(\nn \tt\varphi)^2= Q$   on   $D_{r_0}\setminus\pp D$,
by
 \eqref{PDEE}, \eqref{BAB}, \eqref{BAA} and \eqref{NB1}, we see that $u_t^\ll$  in \eqref{UUK} solves the PDE
\beq\label{ULL} \beg{split}  (\pp_t +{\rm tr}\{a_t\nn^2\}+\nn_{b_t^{(1)}+\tt b_t}  )u_t^\ll= \ll  u_t^\ll  +   f_t^{(1)}+f_t^{(2)},\ \
  t\in [t_0,t_1],   u_{t_1}^\ll=0,\end{split}\end{equation}
where outside the null set $\pp D$,
  \beg{align*}   &f_t^{(1)}:=   (h\circ \rr_D)f_t\circ\tt\varphi +   2\big\<a_t\nn(h\circ \rr_D), \nn \{\tt u_t^\ll\circ \tt\varphi\} \big\> ,\\
  & f_t^{(2)}:=     (\tt u_t^\ll\circ\tt\varphi)(L_t^{\si,b^{(1)}}  +\nn_{ \tt b_t}) (h\circ\rr_D). \end{align*}
By \eqref{NB20}, $h\in C^\infty([0,\infty))$ with support ${\rm supp} h\subset [0, r_0/2]$,  $\|a\|_\infty+\|1_{\pp_{r_0} D} \nn_{b^{(1)}}\rr\|_\infty<\infty $   according to $(A_2^{\si,b})$ and Remark 2.2(1),
   we find a constant $c>0$ such that
\beg{align*}&|f_t^{(1)}|\le 1_{\{\rr_D\le \ff{r_0}2\}} (|f_t| + |\nn\tt u_t^\ll|)\circ\tt\varphi,\\
   &|f_t^{(2)}|\le c 1_{\{\rr_D\le \ff{r_0}2\}}  \big\{(1+|\tt b_t|)|\tt u_t^\ll|    \big\}\circ\tt\varphi.\end{align*}
Since $|f|+ |\tt b|+|\tt u^\ll|$ is bounded on $[0,T]\times \bar D$, so is  $|f^{(1)}|+|f^{(2)}|$ on $[0,T]\times\R^d$.  Hence, by Lemma \ref{LYZ},    the PDE \eqref{ULL} has a unique solution in $\tt H_q^{2,p}(t_0,t_1)$,
  for each $i=1,2$ and $\ll\ge 0$, the PDE
\beq\label{ULLi}   (\pp_t + {\rm tr}\{a_t\nn^2\}+\nn_{b_t^{(1)}+\tt b_t})u_t^{\ll,i}= \ll  u_t^{\ll,i}  +   f_t^{(i)},
\ \ t\in [t_0,t_1],   u_{t_1}^{\ll,i}=0\end{equation}
has a unique solution in $\tt H_q^{2,p}(t_0,t_1)$ as well,   and
  there exist  constants $c_1,c_2>0$ increasing in $\|\tt b\|_{\tt L_{q'}^{p'}(T,D)} $  such that
\beq\label{AO1}  \beg{split}   &\ll^{1-\ff d p-\ff 2 q}  \|u^{\ll,1}\|_\infty+ \ll^{\ff 1 2 (1-\ff d p-\ff 2 q)} \|\nn u^{\ll,1}\|_{\tt L_{q}^{p}(t_0,t_1)}\\
 &\le c_1\|f^{(1)}\|_{\tt L_q^p(t_0,t_1)}
   \le c_2  \big(\|f\|_{\tt L_{q/2}^{p/2}(t_0,t_1,D)} + \|\tt u_t^\ll\|_{\tt L_{q}^p(t_0,t_1,D)}\big),\ \
   p >2,\end{split}  \end{equation}
\beq\label{AO2} \beg{split}    &\ll^{\ff 1 2 (1-\ff d p-\ff 2 q)}
  \|\nn u^{\ll,1}\|_{\infty}   +\|\nn^2 u^{\ll,1}\|_{\tt L_q^p(t_0,t_1)}
 + \|(\pp_t+\nn_{b^{(1)}}) u^{\ll,1}\|_{\tt L_q^p(t_0,t_1)}\\
&\le  c_1\|f^{(1)}\|_{\tt L_q^p(t_0,t_1)}
   \le c_2 ( \|f\|_{\tt L_{q}^{p}(t_0,t_1,D)}+  \|\tt u^\ll\|_{\tt L_{q}^{p}(t_0,t_1,D)}),\end{split} \end{equation}
 and
\beq\label{AO3} \beg{split} &\ll^{\ff 1 2 (1-\ff d {p'}-\ff 2 {q'})}
  (\|u^{\ll,2}\|_\infty+\|\nn u^{\ll,2}\|_{\infty})   +\|\nn^2 u^{\ll,2} \|_{\tt L_{q'}^{p'}(t_0,t_1)}\\
& + \|(\pp_t+\nn_{b^{(1)}}) u^{\ll,2}\|_{\tt L_{q'}^{p'}(t_0,t_1)}
  \le c_1 \|  f^{(2)}\|_{\tt L_{q'}^{p'}(t_0,t_1)}
 \le c_2   (1+\|\tt b\|_{\tt L_{q'}^{p'}(t_0,t_1,D)}) \|\tt u^\ll\|_\infty.\end{split} \end{equation}
By taking large enough $\ll_0>0$ increasing in $\|\tt b\|_{\tt L_{q'}^{p'}(T,D)}$, we derive from \eqref{AO1} and \eqref{AO3} that
\beg{align*} &\|u^{\ll, 1}\|_\infty +\|\nn u^{\ll, 1}\|_{\tt L_q^p(t_0,t_1)}\le \ff 1 2 \big(\|f\|_{\tt L_{q/2}^{p/2}(t_0,t_1,D)} + \|\tt u_t^\ll\|_{\tt L_{q}^p(t_0,t_1,D)}\big),\\
&\|u^{\ll, 2}\|_\infty +\|\nn u^{\ll, 2}\|_{\infty} \le \ff 1 2 \|\tt u^\ll\|_\infty,\ \ \ll\ge\ll_0.\end{align*}
Noting that the uniqueness of   \eqref{ULL} and \eqref{ULLi} implies $u_t^\ll= u_t^{\ll,1}+u_t^{\ll,2}$, this and the definition of $ u_t^\ll$ yield
\beg{align*} &\|\tt u^\ll \|_\infty + \|\nn \tt u^{\ll}\|_{\tt L_q^p(t_0,t_1,D)} \le \sum_{i=1}^2( \|u_t^{\ll, i}\|_\infty +\|\nn u^{\ll, i}\|_{\tt L_q^p(t_0,t_1)})\\
&\le \ff 1 2 \big\{\|\tt u^\ll\|_\infty+ \|f\|_{\tt L_{q/2}^{p/2}(t_0,t_1,D)} + \|\tt u_t^\ll\|_{\tt L_{q}^p(t_0,t_1,D)}\big\},\end{align*}
so that
$$\|\tt u^\ll \|_\infty + \|\nn \tt u^{\ll}\|_{\tt L_q^p(t_0,t_1,D)} \le \|f\|_{\tt L_{q/2}^{p/2}(t_0,t_1,D)},\ \ \ll\ge \ll_0.$$
This together with  \eqref{AO1}-\eqref{AO3} implies \eqref{KG}, \eqref{KG'} and \eqref{J2} for some $c,\vv>0$.
 \end{proof}

\beg{lem} \label{L1'}  Assume   $(A_2^{\si,b})$   but without  the condition  on $\|\nn\si \|$.   For any $(p,q)\in \scr K$ with $p>2$,
there exist  a constant $i\ge 1$ depending only on $(p,q)$, and a constant $c\ge 1$ increasing in $\|b^{(0)}\|_{\tt L_{q_0}^{p_0}(T,D)}$,  such that  for    any solution $(X_t)_{t\in [0,T]}$ of  $\eqref{E01}$,   and any $0\le t_0\le t_1\le T,$
\beq\label{KR-2} \E\bigg(\int_{t_0 }^{t_1 } |f_s(X_s)|\d s\bigg|\F_{t_0}\bigg)^m
 \le c^m m!  \|f\|_{\tt L_{q/2}^{p/2}(t_0,t_1)}^m,\ \  f\in \tt L_{q/2}^{p/2}(t_0,t_1), m\ge 1,\end{equation}
 \beq\label{KAS-2} \E \big(\e^{ \int_{t_0}^{T} |f_t(X_t)|\d t}\big|\F_{t_0}\big)  \le \exp\Big[c +c  \|f\|_{\tt L_{q/2}^{p/2}(t_0,T)}^i\Big],\ \
  f\in \tt L_{q/2}^{p/2}(t_0,T), t_0\in [0,T]. \end{equation}
  \end{lem}

 \beg{proof} As explained in step (1) of the proof of Lemma \ref{L1},   it suffices to prove \eqref{KR-2} for $m=1$ and $f\in C_0^\infty([t_0,t_1]\times \R^d).$

 Let   $(b^{0,n})_{n\ge 1}$ be the mollifying approximations  of   $b^{(0)}=1_{\bar D}b^{(0)}.$
  We have
\beq\label{APP-1}\|b^{0,n}\|_{\tt L_{q_0}^{p_0}(T)}\le  \|b^{(0)}\|_{\tt L_{q_0}^{p_0}(T)},\ \  \lim_{n\to\infty}  \|b^{0,n}- b^{(0)}\|_{\tt L_{q_0}^{p_0}(T)} =0.\end{equation}
By Lemma \ref{LNN}  for $(f,0,\cdots,0)$ replacing $f$, there exist constants $c, \ll_0>0$ such that for any $\ll\ge \ll_0$, the following PDE on $\bar D$
\beq\label{N2PDE}  (\pp_t+L_t^{\si,b^{(1)}} +\nn_{b_t^{0,n}}-\ll) u_t^{\ll, n} =f_t,\ \ t\in [t_0,t_1),
 \ \nn_\n u_t^{\ll, n}|_{\pp D}=0, u_{t_1}^{\ll, n}=0\end{equation}
 has a unique solution in $C^{1,2}([t_0,t_1]\times\bar D)$,  and   for some constant $c_1>0$ we have
\beq\label{ESTN} \|u^{\ll, n}\|_\infty\le c_1  \|f\|_{\tt L_{q/2}^{p/2}(t_0,t_1,D)},\ \
 \|\nn u^{\ll,n}\|_\infty \le c_1 \|f\|_{\infty}, \ \ \ll\ge \ll_0,n\ge 1.  \end{equation}
Moreover, since $(A_2^{\si,b})$ implies $(A_0^{\si,b})$ due to Lemma \ref{LNN}, by \eqref{KR} for $f=|b^{(0)}-b^{0,n}|,$ we find a constant $c_2>0$ such that
\beq\label{EST-N'}  \E\bigg(\int_{t_0}^{t_1} |b^{(0)}-b^{0,n}|(X_s)\d s\bigg|\F_{t_0}\bigg) \le c_2 \|b^{(0)}-b^{0,n}\|_{\tt L_{q_0}^{p_0}(t_0,t_1)},\ \ n\ge 1.  \end{equation}
By \eqref{N2PDE} and $u^{\ll, n}\in C_b^{1,2}([t_0,t_1]\times \bar D)$, we have the following It\^o's formula
\beg{align*} \d u_t^{\ll, n}(X_t) &= (\pp_t +   L_t )u_t^{\ll, n}(X_t)\d t+\d M_t\\
 &= \{f_t+\nn_{b_t^{(0)}-b_t^{0,n}} u_t^{\ll ,n}\}(X_t)\d t+\d M_t\end{align*}
for some martingale $M_t$. Combining this with
\eqref{ESTN}   and \eqref{EST-N'}, we obtain
$$ \E\bigg(\int_{t_0}^{t_1}f_t(X_t)\d t\bigg|\F_{t_0}\bigg)\le c_1   \|f\|_{\tt L_{q/2}^{p/2}(t_0,t_1)} +
c_1c_2 \|f\|_{\infty}\|b_t^{(0)}-b_t^{0,n}\|_{\tt L_{q_0}^{p_0}(t_0,t_1)}.$$
Therefore, by \eqref{APP-1}, we may let     $n\to\infty$ to derive \eqref{KR-2} for $m=1$.
\end{proof}

 \subsection{Weak well-posedness: proof of Theorem \ref{T2.1} }

 We first introduce some known results for the    reflecting SDE   with random coefficients:
\beq\label{E011} \d X_t= J_t(X_t) \d t+ S_t(X_t)\d W_t + \n(X_t)\d l_t,\ \ t\in [0,T],\end{equation}
where  $(W_t)_{t\in [0,T]}$ is an $m$-dimensional Brownian motion on a complete filtration probability space $(\OO,\{\F_t\}_{t\in [0,T]},\P)$,
$$J: [0,T]\times\OO\times\R^d\to \R^d,\ \  S: [0,T]\times\OO\times\R^d\to \R^d\otimes\R^m$$
are progressively measurable, and $l_t$ is the local time of $X_t$ on $\pp D$.
Let $\LL$ be the set of increasing functions $h: (0,1]\to (0,\infty)$ such that $\int_0^{1} \ff{\d s} {h(s)}=\infty$, and let $\GG$ be the class of increasing functions $\gg: [0,\infty)\to [1,\infty)$ such that
$\int_0^\infty \ff{\d s}{\gg(s)}=\infty.$
When $D$ is convex the following result goes back to \cite{Tanaka}, and in general it is mainly summarized from \cite[Theorem 1, Corollary 1 and Theorem 2]{Hino}, where the condition in the first assertion is more general than that stated in \cite[Theorem 1.1]{Hino}:
$$\|S_t(x)-S_t(y)\|^2_{HS}+2 \<x-y, J_t(x)-J_t(y)\>\le g_t h(|x-y|^2),\ \ t\in [0,T],  x,y\in\bar D,$$
since in the proof of this assertion, one only uses the upper bound of
$$\|S_t(X_t)-S_t(Y_t)\|^2_{HS}+2 \<X_t-Y_t, J_t(X_t)-J_t(Y_t)\>,$$  so that the present condition is enough for the pathwise uniqueness.
In Theorem \ref{T01}(3), the term ${\rm tr}\{S_tS_t^*\nn^2 V_t\}$ was formulated in \cite[Theorem 1.1]{Hino} as $\|S_t(x)\|^2 \DD V_t(x)$,
which should be changed into the present one   according to It\^o's formula of $V_t(X_t)$.  Moreover, when $S$ and $J$ are bounded and deterministic, the weak existence is given in \cite[Theorem 2.1]{RS}.

\beg{thm}[\cite{Hino,RS,Tanaka}] \label{T01}  Assume {\bf (D)}.
\beg{enumerate} \item[$(1)$] For any two solutions $X_t$ and $Y_t$ of \eqref{E011} with $X_0=Y_0\in \bar D$, if   there exist $h\in \LL$ and a positive $L^1([0,T])$-valued random variable $g$ such that $\P$-a.s.
$$\|S_t(X_t)-S_t(Y_t)\|^2_{HS}+2 \<X_t-Y_t, J_t(X_t)-J_t(Y_t)\>\le g_t h(|X_t-Y_t|^2),\ \ t\in [0,T], $$
then $X_t=Y_t$   up to life time.
\item[$(2)$] If $\P$-a.s. $S$ and $J$ are  continuous and  locally bounded on $ [0,\infty)\times\bar D$,   then  for any initial value in $\bar D$, $\eqref{E011}$ has a weak solution up to life time. If  
$S$ and $J$ are     bounded and deterministic $S$ and $J$ on $[0,T]\times\bar D$, $\eqref{E011}$  has a global weak solution.
\item[$(3)$] If  either $D$ is bounded, or there exist  $1\le V\in C^{1,2}([0,T]\times \bar D)$ with
$$\lim_{x\in \bar D, |x|\to \infty} \inf_{t\in [0,T]}V_t(x)=\infty, \ \   \nn_{\n} V_t|_{\pp D}  \le 0,$$ and a positive $L^1([0,T])$-valued random variable $g$ such that $\P$-a.s.
\beg{align*} &{\rm tr}\{S_tS_t^*\nn^2 V_t\} + 2 \<\nn V(x), J_t(x)\>+ 2\pp_t V_t(x)\\
&\le g_t \gg(V(x)),\ \ t\in [0,T], x\in \bar D\end{align*}
holds for some $\gg\in \GG$, then any solution to  $\eqref{E011}$ is non-explosion. \end{enumerate}
\end{thm}

Next, we apply Theorem \ref{T01} to \eqref{E01} with coefficients satisfying the following assumption, where
$(1_b)$ is known as monotone or semi-Lipschitz condition, which comparing with $(1_a)$ allows $\si$ to be unbounded.

\emph{\beg{enumerate} \item[{\bf (H1)}]     $b$ and $\si$  satisfying the following conditions.
\item[$(1)$] One of the following conditions holds:
\item[$(1_a)$]  $(A_0^{\si,b})$ holds with $\|\nn \si\|^2\in \tt L_q^p(T)$ for some $(p,q)\in \scr K$, or $(A_2^{\si,b})$ holds. Moreover,   there exists a   constant   $K>0$  such that
\beq\label{H1'0}        \<x-y,b_t(x)-b_t(y)\>
   \le K  |x-y|^2,\ \ t\in [0,T], x,y\in \bar D. \end{equation}
\item[$(1_b)$]   There exists an increasing function $h: [0,\infty)\to [0,\infty)$ with $\int_0^{1}\ff{\d r}{r+h(r)}=\infty$,   such that
\beq\label{H1'00}       2 \<x-y,b_t(x)-b_t(y)\>^+  +\|\si_t(x)-\si_t(y)\|_{HS}^2\le h(|x-y|^2),\ \ t\in [0,T], x,y\in \bar D. \end{equation}
\item[$(2)$]     $\|\si \|\le c(1+|\cdot|^2)$ holds for some constant $c>0$,   there exist     $x_0\in D$ and $\tt\pp D\subset \pp D$ such that
\beq\label{CVX'} \<x-x_0,\n(x)\>\le 0,\ \ x\in \pp D\setminus\tt\pp D, \ \n(x)\in \scr N_x;\end{equation}
and when $\tt\pp D\ne \emptyset$ there exists a function $\tt\rr\in C_b^2(\bar D)$ such that
\beq\label{GB1'} \<\nn \tt\rr,\n\>|_{\pp D}\ge 1_{\tt\pp D},\ \ \sup_{ [0,T]\times\bar D} \big\{ \|\si^*\nn \tt\rr\|+
\|{\rm tr}\{\si\si^*\nn^2\tt\rr\}\|+ \<b,\nn \tt\rr\>^-\big\}\le K.\end{equation}
    \end{enumerate}}
According to \eqref{CVX} and Remark 2.1(a),  {\bf (H1)}(2) holds with $\tt\rr=0$  if either $D$ is convex,  and it holds with $\tt\rr=\rr$ in $\rr_{r_0/2}D$ for some $r_0>0$ when   $\pp D\in C_b^2$ and $\|\si\|+  \<b,\nn\rr\>^-$  is bounded on $[0,T]\times \pp_{r_0}D$.

\beg{lem}\label{L0} Assume  {\bf (D)} and {\bf (H1)}$(1)$. Then the reflecting SDE
   $\eqref{E01}$ is well-posed up to life time.
 If     {\bf (H1)}$(2)$ holds, then the solution is non-explosive, and for any  $k>0$    there exists a constant  $c >0$ such that
\beq\label{EST}  \E\Big[ \sup_{t\in [0,T]} |X_t^x|^k\Big]\le c(1+|x|^k),\ \ x\in \bar D, t\in [0,T],\end{equation}
\beq\label{LC}\sup_{x\in\bar D} \E \big(\e^{k (\tt l_{t_1}^x-\tt l_{t_0}^x)}|\F_{t_0}\big)\le c,\  \ 0\le t_0\le t_1\le T,\end{equation}
where $(X_t^x,l_t^x)$ is  the solution with $X_0^x=x$, and  $\tt l_t^x:=\int_0^t 1_{\tt\pp(D)}(X_s^x)\d l_s^x.$
   \end{lem}

To prove this result,  we need the following lemma on the maximal functional for  nonnegative functions $f$ on $\bar D$:
$$\scr M_D f(x):= \sup_{r\in (0,1)} \ff 1 {|B(0,r)|}\int_{B(0,r)} (1_{D}f)(x+y)\d y,\ \ x\in\bar D.$$

\beg{lem}\label{NN} Let $\pp D\in C_b^2$.
\beg{enumerate}
\item[$(1)$] For any real function $f$ on $\bar D$ with $|\nn f|\in L^1_{loc}(\bar D)$,
$$|f(x)-f(y)|\le c|x-y|\big(\scr M_D |\nn f|(x)+ \scr M_D |\nn f|(y)+\|f\|_\infty \big),\ \ {\rm a.e.}\ x,y\in \bar D.$$
\item[$(2)$] There exists a constant $c>0$ such that for any  nonnegative measurable function $f$ on $[0,T]\times\bar D$,
$$\|\scr M_D f\|_{\tt L_q^p(T,\bar D)}  \le c\|f\|_{\tt L_q^p(T,\bar D)},\ \ p,q\ge 1.$$
\end{enumerate}
\end{lem}

\beg{proof} We only prove (1), since (2) follows from   \cite[Lemma 2.1(ii)]{XXZZ} with $1_{\bar D}f$ replacing $f$.
Let $\tt\varphi$ be in \eqref{SLL}.  Take $0\le h\in C_b^\infty(\R)$ with $h(r)=1$ for $r\le r_0/4$ and $h(r)=0$ for $r\ge r_0/2.$
We then extend a function $f$ on $\bar D$ to $\tt f$ on $\R^d$ by letting
$$\tt f(x):= \{h\circ \rr_D\} f\circ\tt\varphi,$$ where $\rr_D$ is the distance function to $D$.
Then there exists a constant $c>0$ such that
$$|\nn \tt f|\le 1_{\bar D} |\nn f| + c1_{\pp_{-r_0/2}D}(|f\circ\tt\varphi|+ |\nn f|\circ\tt\varphi).$$
By \cite[Lemma 5.4]{Z2} and the  integral transform $x\mapsto \tt\varphi(x)$ with $\|(\nn\tt\varphi)^{-1}\|$ bounded on $\pp_{-r_0}D$,
we find constants $c_1,c_2>0$ such that for any $x,y\in \bar D$,
\beg{align*} &|f(x)-f(y)|=|\tt f(x)-\tt f(y)|\\
&\le c_1 |x-y|\big(\scr M |\nn \tt f|(x)+   \scr M |\nn f|(y)+\|f\|_\infty\big\}\\
&\le c_2|x-y| \big\{\scr M_D |\nn f|(x)+\scr M_D |\nn f|(y)+ \|f\|_\infty\big\},\end{align*}
where $\scr M:= \scr M_D$ for $D=\R^d.$
\end{proof}

 \beg{proof}[Proof of Lemma \ref{L0}]  (1) We first prove the existence and uniqueness up to life time.  Since $\si$ and $b$ are locally bounded, by a truncation argument we may and do assume that $\si$ and $b$ are bounded. Indeed, let for any $n\ge 1$ we take
  $$\si^{\{n\}}_t(x):= \si_t\big(\{1\land (n/|x|)\} x\big),\ \ b_t^{\{n\}}(x):= h(|x|/n) b_t (x),\ \ t\ge 0, x\in\bar D,$$
  where $h\in C_0^\infty([0,\infty)$ with $0\le h\le 1$ and $h|_{[0,1]}=1.$
  Then $\si^{\{n\}}$ and $b^{\{n\}}$ are bounded on $[0,T]\times\bar D$ and for some constant $K_n>0$,
  \beg{align*} &\<b_t^{\{n\}}(x)- b_t^{\{n\}}(y), x-y\>^+\\
  &\le h(|x|/n)\<b_t(x)- b_t(y), x-y\>^++  \big|h(|x|/n)- h(|y|/n)\big|\<b_t(y),x-y\>^+\\
  &\le \<b_t(x)- b_t(y), x-y\>^+ + K_n |x-y|^2,\ \ t\in [0,T], x,y\in \bar D, |y|\le |x|.\end{align*}
  So, by the symmetry of $\<b_t^{\{n\}}(x)- b_t^{\{n\}}(y), x-y\>^+$ in $(x,y)$,  under  $(1_a)$, $\si$ and $b^{\{n\}}$ are bounded on $[0,T]\times\bar D$ and satisfy \eqref{H1'0}   with $K+K_n$ replacing $K$; while
   $(1_b)$ and $$|\{1\land (n/|x|)\} x-\{1\land (n/|y|)\} y|\le |x-y|$$ imply that  $\si^{\{n\}}$ and $ b^{\{n\}}$  are bounded and satisfy \eqref{H1'00} for $2h(r)+K_nr$ replacing $h(r)$. Therefore,
  if the well-posedness is proved under {\bf (H1)} for bounded $b$ and $\si$, the SDE is well-posed up to the hitting time of $\pp B(0,n)$ for any $n\ge 1$, i.e. it is well-posed up to life time.

  When $\si$ and $b$ are bounded, the weak existence is implied by Theorem \ref{T01}(2).
 By the Yamada-Watanabe principle, it suffices to verify the pathwise uniqueness.  Let
 $X_t$ and $Y_t$ be two solutions starting from $x\in \bar D$. By Lemma \ref{NN}(1) and {\bf (H1)}(1),
$$\|\si_t(X_t)-\si_t(Y_t)\|^2_{HS}+2 \<X_t-Y_t, b_t(X_t)-b_t(Y_t)\>\le \beg{cases} g_t  |X_t-Y_t|^2,&\text{under\ } (1_a),\\
h(|X_t-Y_t|^2), \  &\text{under\ } (1_b),\end{cases} $$ where for some constant $c >0$
$$g_t:=   c \big\{ 1+ \scr M_D \|\nn \si_t\|^2(X_t)+ \scr M_D \|\nn \si_t\|^2(Y_t)\big\}. $$
So, by Theorem \ref{T01}(1),   it suffices to prove   $\int_0^T g_t\d t<\infty$ under $(1_a)$.
By Lemma \ref{NN}, this follows from   \eqref{KR} under condition $(A_0^{\si,b})$   with $\|\nn \si\|^2\in \tt L_{q}^{p}(T)$ for   some $(p,q)\in \scr K$,   or  \eqref{KR-2} under condition $(A_2^{\si,b})$.

(2)  To prove the non-explosion, we simply denote $(X_t,l_t)=(X_t^x,l_t^x)$ and let
 $$\tau_n:=\inf\{t\ge 0: |X_t|\ge n\},\ \ n\ge 1.$$
 By {\bf (H1)}(2),  we find a constant $c_1>0$ such that
\beq\label{LJ}\d \tt\rr(X_t)\ge -K \d t+\d M_t +\d \tt l_t,\ \ t\in [0,T]\end{equation}
holds for $\d M_t:= \<\si_t(X_t)^*\nn \tt\rr(X_t), \d W_t\>$ satisfying $\d\<M\>_t\le K^2\d t$.
  This implies  \eqref{LC}.
Next, by {\bf (H1)}, we find a constant  $c_1 >0$ such that
\beg{align*} &2\<b_t(x), x-x_0\>+ \|\si_t(x)\|_{HS}^2\\
&=2\<b_t(x)-b_t(x_0), x-x_0\> + \|\si_t(x)-\si_t(x_0)\|_{HS}^2 \\
&\qquad + 2 \<b_t(x_0), x-x_0\> +\|\si_t(x_0)\|_{HS}^2 + 2\<\si_t(x_0),\si_t(x)\>_{HS}\\
&\le  c_1(1+|x-x_0|^2),\  \ x\in\bar D.\end{align*}
Then by {\bf (H1)}(2) and It\^o's formula, for any $k\ge 2$ we find a constant $c_2>0$ such that
$$\d |X_t-x_0|^k\le c_2 (1+|X_t-x_0|^k)\d t + \d \tt M_t+ k |X_t-x_0|^{k-1}\d \tt l_t,$$
where $\tt M_t$ is a local  martingale with $\d\<\tt M\>_t\le c_2 (1+|X_t-x_0|^k)^2\d t.$
By BDG's inequality and \eqref{LC}, we find   constants $c_3, c_4>0$ such that
$$\eta_t^{\{n\}}:= \sup_{s\in [0,t\land\tau_n]}(1+|X_s-x_0|^k),\ \ n\ge 1, t\in [0,T]$$ satisfies
\beg{align*}  \E \eta_{t}^{\{n\}}   &\le 1+ |x-x_0|^k+   c_3\E\int_0^{t}\eta_{t}^{\{n\}} \d s + 2c_3 \E^x \bigg(  \int_0^{t}  |\eta_{t}^{\{n\}}| ^2 \d s\bigg)^{\ff 1 2} +k\E \Big[|\eta_{t}^{\{n\}}|^{\ff{k-1}k} \tt l_t\Big] \\
&\le  \ff 1 2\E \eta_{t}^{\{n\}}   + c_4(1+|x|^k) + c_4\int_0^t \E\eta_{s}^{\{n\}}\d s,\ \ t\in [0,T].\end{align*}
By Gronwall's lemma, we obtain
$$ \E [\eta_{t}^{\{n\}}]\le 2c_4(1+|x|^k)\e^{2c_4 t},\ \ t \in [0,T], x\in\bar D, n\ge 1, $$
which   implies the non-explosive of $X_t$ and  \eqref{EST} for some constant $c>0$.
 \end{proof}

  \beg{proof}[Proof of Theorem $\ref{T2.1}$] Let $X_0=x\in \bar D$.  We consider the following two cases respectively.

  (a) Let $(A^{\si,b}_1)$ hold.
  Then  {\bf (H1)}   holds for $b^{(1)}$ replacing $b$.
 By Lemma \ref{L0},  the reflecting  SDE
 \beq\label{QPP} \d X_t = b^{(1)}_t(X_t)\d t+ \si_t(X_t)\d W_t+ \n(X_t)\d l_t\end{equation}
  is well-posed with \eqref{EST} holding for all $k\ge 1$ and some constant $c>0$ depending on $k$.   By Lemmas \ref{L1}-\ref{L1'}, \eqref{LC} and  $(A_0^{\si,b})$ with $|b^{(0)}|^2\in \tt L_{q}^{p}(T)$,     we see that  \eqref{KAS} holds for $f:= |b^{(0)}|^2$, so that for some map $c: [1,\infty)\to (0,\infty)$   independent of the initial value $x$,
  \beq\label{YFF} \sup_{x\in \bar D} \E^x |R_T|^k\le c(k),\ \ k\ge 1\end{equation} holds for
    $$R_t:= \e^{\int_0^t\<\{\si_s^* (\si_s\si_s^*)^{-1} b_s^{(0)}\}(X_s), \d W_s\>- \ff 1 2 \int_0^t |\si_s^* (\si_s\si_s^*)^{-1} b_s^{(0)}|^2(X_s)\d s},\ \ t\in [0,T].$$
   By Girsanov's theorem,
  $$\tt W_t:= W_t- \int_0^t \{\si_s^* (\si_s\si_s^*)^{-1} b_s^{(0)}\}(X_s) \d s,\ \ t\in [0,T]$$
  is an $m$-dimensional Brownian motion under the probability measure $\Q:= R_T\P$.  Rewriting \eqref{QPP} as
  $$\d X_t = b_t(X_t)\d t+ \si_t(X_t)\d \tt W_t+ \n(X_t)\d l_t,$$ we see that $(X_t, l_t,\tt W_t)_{t\in [0,T]}$ under probability $\Q$ is a weak solution of \eqref{E01}.
  Moreover, letting $\E_\Q$ be the expectation under $\Q$, by  \eqref{EST} and \eqref{YFF}, for any $k\ge 1$ we find a constant $\tt c(k)>0$ independent of $x$ such that
  \beg{align*} &\E_\Q\Big[\sup_{t\in [0,T]} |X_t|^k\Big]= \E\Big[R_T\sup_{t\in [0,T]} |X_t|^k\Big]\\
  &\le \big(\E\big[R_T^2\big]\big)^{\ff 1 2}
  \Big(\E\sup_{t\in [0,T]} |X_t|^{2k}\Big]\Big)^{\ff 1 2}\le \tt c(k) (1+|x|^k),\ \ x\in\bar D\end{align*}
  for some constant $c>0$. Similarly, \eqref{LC} and \eqref{YFF} imply
  $ \E_\Q \e^{k l_T}\le C(k)$ for $k\ge 1$ and  constant $C(k)>0$ independent of $x$.
  So,   \eqref{EPP} holds for this weak solution.

To prove the weak uniqueness,
  let $(\bar X_t, \bar l_t,\bar W_t)_{t\in [0,T]}$ under probability $\bar \P$ be another weak solution of \eqref{E01} with $\bar X_0=x$, i.e.
  \beq\label{QPP'} \d \bar X_t = b_t(\bar X_t)\d t+ \si_t(\bar X_t)\d \bar W_t+ \n(\bar X_t)\d \bar l_t,\ \ t\in [0,T],  \bar X_0=x.\end{equation}
   It suffices to show
   \beq\label{UNI}  \L_{(\bar X_t,\bar l_t)_{t\in [0,T]}|\bar\P}=\L_{(X_t,l_t)_{t\in [0,T]}|\Q}.\end{equation}
  By Lemma \ref{L1} the estimate \eqref{KAS} holds for
  $\bar X_t$ and $f=|b^{(0)}|^2$, so that
\beq\label{GS} \E_{\bar\P}\e^{\ll \int_0^{T} |b_t^{(0)}(\bar X_t)|^2\d t} <\infty,\ \ \ll >0.\end{equation}
  By Girsanov's theorem, this and $(A^{\si,b}_0)$ imply that
  $$G_t(\bar X,\bar W):= \bar W_t+\int_0^{t} \{\si_s^* (\si_s\si_s^*)^{-1} b_s^{(0)}\}(\bar X_s)\d s,\ \ t\in [0,T]$$
  is an $m$-dimensional Brownian motion under the probability $\bar \Q:= R(\bar X,\bar W)\bar\P,$ where
  $$R(\bar X,\bar W):= \e^{-\int_0^{T} \<\{\si_s^* (\si_s\si_s^*)^{-1} b_s^{(0)}\}(\bar X_s),  \d \bar W_s\> -\ff 1 2 \int_0^{T} |\{\si_s^* (\si_s\si_s^*)^{-1} b_s^{(0)}\}(\bar X_s)|^2\d s}.$$
  Reformulating \eqref{QPP'} as
  $$\d\bar X_t=b_t^{(1)}(\bar X_t)\d t +\si_t(\bar X_t) \d G_t(\bar X,\bar W) +\n(\bar X_t)\d \bar l_t,\ \ t\in [0,T],$$
  and applying the well-posedness of \eqref{QPP} which implies the joint weak uniqueness, we conclude that
  $$\L_{(\bar X_t, \bar l_t, G_t(\bar X,\bar W))_{t\in [0, T]}|\bar \Q}= \L_{(X_t,l_t,W_t)_{t\in [0, T]}|\P}.$$
  Noting that
  $$ R(\bar X, \bar W)^{-1}=  \e^{-\int_0^T |\{\si_s^* (\si_s\si_s^*)^{-1} b_s^{(0)}\}(\bar X_s)|^2\d s} R(\bar X, G(\bar X,\bar W))^{-1},$$
  this implies that for any bounded continuous function $F$ on $C([0,T];\R^d\times [0,\infty))$,
 \beg{align*} &\E_{\bar\P}[F(\bar X,\bar l)]= \E_{\bar \Q} [R(\bar X, \bar W)^{-1} F(\bar X,\bar l)]\\
 &=\E_{\bar \Q} [R(\bar X, G(\bar X,\bar W))^{-1} \e^{-\int_0^T |\{\si_s^* (\si_s\si_s^*)^{-1} b_s^{(0)}\}(\bar X_s)|^2\d s}F(\bar X,\bar l)]\\
   & =  \E_{\P} [R(X, W)^{-1}\e^{-\int_0^T |\{\si_s^* (\si_s\si_s^*)^{-1} b_s^{(0)}\}(X_s)|^2\d s} F(X,l)]\\
   &=\E_\P[R_T F(X,l)]= \E_\Q[F(X,l)].\end{align*}
  Therefore,   \eqref{UNI} holds.

  (b) Let $(A_2^{\si,b})$  hold.   By Lemma \ref{L1'}, \eqref{YFF} and \eqref{GS} hold,  so that the desired assertions follow from Girsanov's transforms as shown in step (a).
   \end{proof}

   \subsection{Well-posedness: proof of Theorem \ref{T2.2}}

   The weak existence  is implied by Theorem \ref{T2.1}. By the Yamada-Watanabe principle,
    it suffices  to prove  estimate \eqref{INN} which in particular implies the pathwise uniqueness as well as estimate \eqref{GRD}:
    \beg{align*} &|\nn P_t f|(x):= \limsup_{\bar D\ni y\to x} \ff{|P_tf(x)-P_tf(y)|}{ |x-y|}\le \limsup_{\bar D\ni y\to x} \E\bigg[\ff{| f(X_t^x)- f(X_t^y)|}{|x-y|}\bigg]\\
& \le \limsup_{\bar D\ni y\to x} \Big(\E\ff{|f(X_t^x)-f(X_t^y)|^p}{|X_t^x-X_t^y|^p}\Big)^{\ff 1 p} \Big(\ff{\E[|X_t^x-X_t^y|^{\ff p{p-1}}]}{ |x-y|^{\ff p{p-1}}}\Big)^{\ff {p-1}p}\\
&\le c(p) \big(P_t|\nn f|^p\big)^{\ff 1 p}(x),\ \ x\in \bar D, t\in [0,T],f\in C_b^1(\bar D).\end{align*}
Let $(X_t^{(i)}, l_t^{(i)})$ be two solutions of \eqref{E01} with $X_0^{(i)}=x^{(i)}\in\bar D, i=1,2.$  Below we prove \eqref{INN} in  situations $(i)$ and $(ii)$ respectively.

\beg{proof}[Proof of Theorem \ref{T2.2} under  $(i)$]       In this case,   $D$ is an interval or a half-line.
 For any $\ll>0$, let $u_t^\ll$ be the unique solution to \eqref{PPD5} with $t_0=0, t_1=T$ and  $f=-  b^{(0)}$, that is,
\beq\label{PPD5N} (\pp_t+L_t) u_t^\ll = \ll u_t^\ll-  b^{(0)}_t,\ \ t\in [0,T], u_T^\ll=0.\end{equation}
    By \eqref{PPD5'} with $f=- b^{(0)}\in \tt L_{2q_0}^{2p_0}(T)$, we take large enough $\ll>0$ such that
 \beq\label{PPD5''} \|u^\ll\|_\infty + \|\nn u^\ll\|_\infty \le \ff 1 2,\ \ \|u^\ll\|_{\tt H_{2q_0}^{2,2p_0}(T)} <\infty.\end{equation}
  Then
 $\Theta_t^\ll(x):= x+ u_t^{\ll}(x)$ is a diffeomorphism and there exists a constant $C>0$ such that
\beq\label{HHH} \ff 1 2 |x-y|\le |\Theta_t^\ll(x)-  \Theta_t^\ll(y)| \le 2 |x-y|,\ \ x,y\in \R, t\in [0,T].  \end{equation}
 Let $(X_t^{(i)},l_t^{(i)})$ solve \eqref{E01} for $X_0^{(i)}=x^{(i)}\in \bar D, i=1,2$, and let
$$Y_t^{(i)}:=   \Theta_t^\ll(X_t^{(i)})= X_t^{(i)}+u_t^{\ll}(X_t^{(i)}), \ \ i=1,2.$$
 By   It\^o's formula in Lemma \ref{L1}(2),
\beq\label{YYE}  \d Y_t^{(i)}= B_t (Y_t^{(i)})\d t+ \Sigma_t (Y_t^{(i)}) \d W_t+ \{1+\nn u_t^{\ll}(X_t^{(i)})\}\n(X_t^{(i)}) \d l_t^{(i)}, \ \ i=1,2 \end{equation}
holds for
 \beq\label{B**}   B_t(x):=\{b_t^{(1)}+ \ll u_t^\ll \} \big(\{\Theta_t^\ll\}^{-1}(x)\big),\ \ \ \Sigma_t(x):= \big\{( 1+ \nn u_t^\ll)\si_t\big\} \big(\{\Theta_t^\ll\}^{-1}(x)\big). \end{equation}
By  $(A^{\si,b}_1)$, \eqref{PPD5''}, \eqref{B**}  and $\|\nn b^{(1)}\|_\infty<1$ due to $(A^{\si,b}_0)$,   we find $0\le F_i\in \tt L_{q_i}^{p_i}(T), 0\le i\le l$, such that  
 \beq\label{BDD} \|\nn B\|_\infty<\infty,\ \      \|\nn \Sigma \|^2 \le \sum_{i=0}^l F_i.\end{equation}
Since  $ d=1$,      for any $x\in \pp D$ and $y\in D$ we have
$y-x =|y-x| \n(x),$  so that  \eqref{PPD5''} implies
\beq\label{AACD}  \big\<\Theta_t^{\ll}(y)- \Theta_t^{\ll}(x), \big\{1 + \nn u_t^{\ll}(x)\big\}\n (x)\big\> \ge |y-x|(1- \|\nn u^{\ll}\|_\infty)^2  \ge 0.\end{equation}
Combining this with  \eqref{YYE} and It\^o's formula,  up to a local martingale  we have
\beg{align*}  \d|Y_t^{(1)}- Y_t^{(2)}|^{2k}
 \le   2k  |Y_t^{(1)}-Y_t^{(2)}|^{2k}  \bigg\{ \ff{|B_t(Y_t^{(1)})- B_t(Y_t^{(2)})| }{ |Y_t^{(1)}-Y_t^{(2)}|}+    \ff{k\|\Sigma_t(Y_t^{(1)})-\Sigma_t(Y_t^{(2)})\|_{HS}^2}{|Y_t^{(1)}-Y_t^{(2)}|^2}\bigg\}\d t.
  \end{align*}
So, by Lemma \ref{NN}, we find a constant $c_1>0$ and a local martingale $M_t$ such hat
$$|Y_t^{(1)}- Y_t^{(2)}|^{2k}\le |Y_0^{(1)}- Y_0^2|^{2k}+ c_1 \int_0^t |Y_s^{(1)}- Y_s^{(2)}|^{2k}\d \scr L_s+\d M_t,$$
where
\beq\label{LLT} \scr L_t:= \int_0^t\Big\{1+  \scr M_D \|\nn \Sigma_s\|^2 (Y_s^{(1)})+ \scr M_D \|\nn \Sigma_s\|^2 (Y_s^{(2)}) \Big\}\d s.\end{equation}
Combining this with  \eqref{BDD}, \eqref{KAS}, Lemma \ref{NN}  and   the stochastic Gronwall lemma (see \cite{[13]} or \cite{XZ}),
for any $k>1$ and $p\in (\ff 1 2,1)$, we find   constants $c_2,c_3 >0$ such that
\beg{align*} &\Big(\E\Big[\sup_{s\in [0,t]} \Theta_s^\ll(X_s^{(1)})- \Theta_s^\ll(X_s^{(2)})|^{k}\Big]\Big)^{2} =\Big(\E\sup_{s\in [0,t]} |Y_s^{(1)}- Y_s^{(2)}|^{k}\Big)^{2}\\
&\le c_2 |Y_0^{(1)}-Y_0^{(2)}|^{2k}\big(\E \e^{\ff {c_1p}{p-1} \scr L_t}\big)^{\ff {p-1} p}\le   c_3|\Theta_0^\ll(x^{(1)})-\Theta_0^\ll(x^{(2)})|^{2k}.\end{align*}   This together with \eqref{HHH} implies \eqref{INN} for some constant $c>0$. \end{proof}

 To prove \eqref{INN} under  $(A_2^{\si,b})$,
we need the following lemma due to  \cite[Lemma 5.2]{YZ}, which is contained in the proof of \cite[Lemma 4.4]{DI}. Let $\nn_1$ and $\nn_2$ be the gradient operators in the first and second variables on $\R^d\times\R^d$.

\beg{lem}\label{L3} There exists a function $g\in C^1(\R^d\times \R^d)\cap C^2((\R^d\setminus\{0\})\times \R^d)$ having the following properties
for some constants $k_2>1$ and $k_1\in (0,1):$
\beg{enumerate}
\item[$(1)$] $k_1|x|^2\le g(x,y)\le k_2 |x|^2,\ \ x,y\in \R^d;$
\item[$(2)$] $\<\nn_1 g(x,y), y\>\le 0,\ \ |y|=1, \<x,y\>\le k_1|x|;$
\item[$(3)$] $|\nn_1^i\nn_2^j g(x,y)|\le k_2 |x|^{2-i},\ \ i,j\in \{0,1,2\}, i+j\le 2, x,y\in\R^d.$
\end{enumerate}
\end{lem}

\

\beg{proof}[Proof of Theorem \ref{T2.2} under   $(ii)$]
Let $b^{0,n}$ be the mollifying approximation of $b^{(0)}=1_{\bar D}b^{(0)}$.   By Lemma \ref{LNN}, there exists $\ll_0>0$ such that for any $\ll\ge \ll_0$ and $n\ge 1$, the PDE
\beq\label{POE2} (\pp_t+L_t+\nn_{b^{0,n}_t-b_t^{(0)}} -\ll)u_t^{\ll,n}=- b^{0,n}_t,\ \ u_T^{\ll,n}=\nn_{\n}u_t^{\ll,n}|_{\pp D}=0,\end{equation}
has a unique solution in $C_b^{1,2}([0,T]\times\bar D)$, and there exist constants $\vv,c>0$ such that
\beq\label{POE3} \beg{split} &\ll^\vv\big(\|u^{\ll,n}\|_{\infty} +\|\nn u^{\ll,n}\|_\infty\big) + \|(\pp_t +\nn_{b^{(1)}}) u^{\ll,n}\|_{\tt L_{q_0}^{p_0}(T,D)} +\| \nn^2 u^{\ll,n}\|_{\tt L_{q_0}^{p_0}(T,D)} \\
&\le c \|b^{(0)}\|_{\tt L_{q_0}^{p_0}(T,D)},\ \ \ll\ge \ll_0, n\ge 1.\end{split} \end{equation} Then for large enough   $\ll_0>0$,    $\Theta_t^{\ll,n}:=id+u_t^{\ll,n}$ satisfies
\beq\label{POE3'} \ff 1 2 |x-y|^2\le |\Theta_t^{\ll,n}(x)-\Theta_t^{\ll,n}(y)|^2\le 2 |x-y|^2,\ \  \ll\ge \ll_0, x,y\in \bar D.\end{equation}
Since $\pp D\in C_b^{2,L}$, there exists a constant    $r_0>0$   such that $\rr\in C_b^2(\pp_{r_0}D)$ with $\nn^2\rr$ Lipschitz continuous on $\pp_{r_0}D$. Take  $h\in C^\infty([0,\infty);[0,\infty))$ such that  $h'\ge 0$, $h(r)=r$ for $r\le r_0/2$ and $h(r)=r_0$ for $r\ge r_0$.

Let  $(X_t^{(i)},l_t^{(i)})$ solve  \eqref{E01} starting at $x^{(i)}\in \bar D$ for $i=1,2$.      Alternatively to $|X_t^{(1)}-X_t^{(2)}|^2$, we  consider the process
$$H_t:=  g \big(\Theta_t^{\ll, n} (X_t^{(1)})- \Theta_t^{\ll, n} (X_t^{(2)}),\nn (h\circ\rr)(X_t^{(1)})\big),\ \ t\in [0,T],$$
where  $g$ is in   Lemma \ref{L3}.
By Lemma \ref{L3}(1)  and \eqref{POE3'}, we have
\beq\label{POE5} \ff {k_1} 2   |X_t^{(1)}-X_t^{(2)}|^2 \le H_t\le 2 k_2    |X_t^{(1)}-X_t^{(2)}|^2,\ \ t\in [0,T].\end{equation}
Simply denote
$$\xi_t:= \Theta_t^{\ll, n} (X_t^{(1)})- \Theta_t^{\ll, n} (X_t^{(2)}),\ \  \eta_t:= \nn (h\circ\rr)(X_t^{(1)}).$$
By It\^o's formula,   \eqref{POE2} and $\nn_\n \Theta_t^{\ll, n}|_{\pp D}=\n$ due to $\nn_\n u_t^{\ll, n}|_{\pp D}=0$, we have
\beq\label{POEE}  \beg{split} &\d  \xi_t=  \big\{\ll   u_t^{\ll, n}(X_t^{(1)})-\ll    u_t^{\ll, n}(X_t^{(2)})+ (b_t^{(0)}-b_t^{0,n})(X_t^{(1)})-(b_t^{(0)}-b_t^{0,n})(X_t^{(2)})\big\}\d t\\
 &\qquad +\big\{[(\nn \Theta_t^{\ll, n})\si_t](X_t^{(1)})- [(\nn \Theta_t^{\ll, n})\si_t](X_t^{(2)})\big\}\d W_t+ \n(X_t^{(1)})\d l_t^{(1)}-\n(X_t^{(2)})\d l_t^{(2)},\\
 &\d \eta_t= L_t \nn   (h\circ \rr)(X_t^{(1)})\d t  + \big\{[\nn^2(h\circ\rr)]\si_t\big\}(X_t^{(1)}) \d W_t+ \{\nn_\n \nn(h\circ\rr)\}(X_t^{(1)})\d l_t^{(1)}.\end{split}\end{equation}
Hence, It\^o's formula for $H_t$ reads
 \beq\label{POE6} \d H_t = A_t\d t + B_t^{(1)} \d l_t^{(1)}- B_t^{(2)} \d l_t^{(2)} +\d M_t,\end{equation}
where
 \beq\label{POA}   \beg{split}A_t:= &\,\big\<\nn_1  g(\xi_t,\eta_t),\ \ll u_t^{\ll,n}(X_t^{(1)})-\ll u_t^{\ll,n}(X_t^{(2)})\big\>\\
 &+\big\<\nn_1  g(\xi_t,\eta_t),\ \nn_{b_t^{(0)}-b_t^{0,n}} \Theta_t^{\ll,n} (X_t^{(1)}) - \nn_{b_t^{(0)}-b_t^{0,n}} \Theta_t^{\ll,n} (X_t^{(2)})\big\> \\
&+ \big\<\nn_2  g(\xi_t,\eta_t),\ L_t \nn(h\circ\rr)(X_t^{(1)})\Big\> +\Big\<\nn_1^2    g(\xi_t,\eta_t), N_tN_t^* \big\>_{HS} \\
& +\big\<\nn_1\nn_2   g(\xi_t,\eta_t),\ N_t \si_t(X_t^{(1)})^*\nn^2(h\circ\rr)(X_t^{(1)})\big\>_{HS} \\
& + \big\<\nn_2^2    g(\xi_t,\eta_t),\  \big\{[\nn^2(h\circ\rr)]\si_t\si_t^*  \nn^2(h\circ\rr)\big\}(X_t^{(1)})\big\>_{HS},\\
   N_t := &\big\{(\nn \Theta_t^{\ll, n})\si_t\big\}(X_t^{(1)})-\big\{(\nn \Theta_t^{\ll, n})\si_t\big\}(X_t^{(2)}),\end{split} \end{equation}
\beq\label{POB} \beg{split}B_t^{(1)}:= &\big\<\nn_1  g(\xi_t,\eta_t), \n(X_t^{(1)})\big\>
 + \big\<\nn_2  g(\xi_t,\eta_t), \nn_\n \{\nn (h\circ\rr)\}(X_t^{(1)})\big\>,\\
B_t^{(2)}:= &\big\<\nn_1 g(\xi_t,\eta_t), \n(X_t^{(2)})\big\>, \end{split} \end{equation}
\beq\label{POC} \beg{split} \d M_t=&\,\big\<\nn_1 g(\xi_t,\eta_t), \big[\{(\nn \Theta_t^{\ll,n})\si_t\}(X_t^{(1)})- \{(\nn \Theta_t^{\ll,n})\si_t\}(X_t^{(2)})\big]\d W_t\big\>\\
&\, + \big\<\nn_2 g(\xi_t,\eta_t), \big[\{\nn^2(h\circ\rr)\}\si_t\big](X_t^{(1)}) \d W_t\big\>.\end{split}\end{equation}
In the following we estimate these terms respectively.

Firstly,   \eqref{LCC} implies
$$\<\Theta_t^{\ll,n}(x)-\Theta_t^{\ll,n}(y),\n(x)\>\le \ff{|x-y|^2}{2r_0}+\|\nn u_t^{\ll, n}\|_\infty|x-y|,\ \ x\in\pp D, y\in\bar D.$$
Combining this with \eqref{POE3}, we find   constants $\vv_0,\ll_1>0$ such that for any $\ll\ge \ll_1$,
  \beg{align*}& \<\Theta_t^{\ll,n}(x)-\Theta_t^{\ll,n}(y),\n(x)\> \le k_1 |\Theta_t^{\ll,n}(x)-\Theta_t^{\ll,n}(y)|,\\
 & \ x\in\pp D, y\in \bar D, |x-y|\le \vv_0, n\ge 1, t\in [0,T].\end{align*}
So, Lemma \ref{L3} yields
\beq\label{POE4}\beg{split} &\big\<\nn_1 g (\Theta_t^{\ll,n}(x)-\Theta_t^{\ll,n}(y),\n(x)), \n(x)\big\> \le k_2 1_{\{|x-y|>\vv_0\}} |\Theta_t^{\ll,n}(x)-\Theta_t^{\ll,n}(y)|\\
&\le k_2\vv_0^{-1}  |\Theta_t^{\ll,n}(x)-\Theta_t^{\ll,n}(y)|^2,\ \ x\in\pp D, y\in \bar D,   n\ge 1, t\in [0,T].\end{split}\end{equation}

Next,    by the same reason leading to \eqref{POE4},  we find a constant $c_1>0$ such that
\beq\label{POE4'}\beg{split} &\big\<\nn_1 g (\Theta_t^{\ll,n}(x)-\Theta_t^{\ll,n}(y),\nn (h\circ\rr)(x)), \n(y)\big\> \\
& \ge \big\<\nn_1 g (\Theta_t^{\ll,n}(x)-\Theta_t^{\ll,n}(y),\n(y)), \n(y)\big\>\\
& -\big|  \nn_1 g (\Theta_t^{\ll,n}(x)-\Theta_t^{\ll,n}(y),\nn (h\circ\rr)(y))- \nn_1 g (\Theta_t^{\ll,n}(x)-\Theta_t^{\ll,n}(y),\nn (h\circ\rr)(x)) \big|\\
&\ge - 1_{\{|x-y|>\vv_0\}} k_2\vv_0^{-1}|\Theta_t^{\ll,n}(x)-\Theta_t^{\ll,n}(y)|^2\\
&   -  \|h'\|_\infty\|\nn_1\nn_2 g(\Theta_t^{\ll,n}(x)-\Theta_t^{\ll,n}(y),\cdot)\|_\infty |\Theta_t^{\ll,n}(x)-\Theta_t^{\ll,n}(y)|^2\\
&\ge -c_1  |\Theta_t^{\ll,n}(x)-\Theta_t^{\ll,n}(y)|^2,\ \ x\in\bar D, y\in \pp D,  n\ge 1, t\in [0,T].\end{split}\end{equation}

Moreover, by    $(A_2^{\si,b})$ and  $h\circ\rr \in C_b^{2,L}(\bar D)$,   there exists a constant $C>0$ such that
  $$ |L_t  \{\nn (h\circ \rr)\}|\le C (1+|b^{(0)}_t|),\ \ t\in [0,T].$$
Combining this with   Lemma \ref{L3}, Lemma \ref{NN}, \eqref{POE5},  and \eqref{POA}-\eqref{POE4'}, we find a constant $K>0$ such that
\beg{align*} &|A_t|\le   K \big\{|b_t^{(0)}-b_t^{0,n}|^2(X_t^{(1)})+ |b_t^{(0)}-b_t^{0,n}|^2(X_t^{(2)}) \big\}\\
&+ K  |X_t^{(1)}-X_t^{(2)}|^2
\bigg\{1+|b_t^{(0)}|(X_t^{(1)})+\sum_{i=1}^2\scr M_D \big\|\nn\{(\nn \Theta_t^{\ll,n})\si_t\}\big\|^2(X_t^{(i)})  \bigg\}, \\
&\d\<M\>_t\le K |X_t^{(1)}-X_t^{(2)}|^4 \bigg\{1+ \sum_{i=1}^2\scr M_D \big\|\nn\{(\nn \Theta_t^{\ll,n})\si_t\}\big\|^2(X_t^{(i)})  \bigg\},\\
&B_t^{(1)}\le   K  |X_t^{(1)}-X_t^{(2)}|^2,\ \ \ \  -B_t^{(2)}\le K  |X_t^{(1)}-X_t^{(2)}|^2.\end{align*}
Combining these with \eqref{POE5} and \eqref{POE6},  for any $k\ge 1$, we find a constant $c_1>0$ such that
\beq\label{POD} \beg{split} \d H_t^k \le &\,c_1 |X_t^{(1)}-X_t^{(2)}|^{2(k-1)} \big\{|b_t^{(0)}-b_t^{0,n}|^2(X_t^{(1)})+ |b_t^{(0)}-b_t^{0,n}|^2(X_t^{(2)}) \big\}\d t \\
&+ c_1 |X_t^{(1)}-X_t^{(2)}|^{2k} \d \scr L_t + k H_t^{k-1} \d M_t,\end{split}\end{equation}
where
\beq\label{LAT'}   \beg{split}    \scr L_t:=   l_t^{(1)}+  l_t^{(2)}
 +\int_0^t \Big\{1+|b_s^{(0)}|(X_s^{(1)})+\sum_{i=1}^2\scr M_D \big\|\nn \{(\nn \Theta_s^{\ll,n})\si_s\}\big\|^2(X_s^{(i)})  \Big\} \d s.\end{split}\end{equation}
For any $j\ge 1$, let
$$\tau_j:= \inf\big\{t\ge 0: |X_t^{(1)}-X_t^{(2)}|\ge j\big\}.$$
By \eqref{POE5} and \eqref{POD}, we find a constant $c_2>0$ such that
\beq\label{LAT}  |X_{t\land \tau_j}^{(1)}- X_{t\land \tau_j}^{(2)}|^{2k}\le   G_j(t) +
   c_2\int_0^{t\land\tau_j} |X_{s}^{(1)}- X_{s}^{(2)}|^{2k}\d \scr L_s +\tt M_t  \end{equation}
holds for some local martingale $\tt M_t$ and
$$G_j(t):=  c_2 |x^{(1)}-x^{(2)}|^{2k}+c_2 j^{2(k-1)} \int_0^{t\land\tau_j} \big\{|b_s^{(0)}-b_s^{0,n}|^2(X_s^{(1)})+ |b_s^{(0)}-b_s^{0,n}|^2(X_s^{(2)})\big\}\d s. $$
 Since   $(A_2^{\si,b})$ and \eqref{POE3} imply
$$\sup_{n\ge 1} \big\| \nn \{(\nn \Theta^{\ll,n})\si\}\big\| \le \sum_{i=0}^l F_i $$ for some $0\le F_i\in \tt L_{q_i}^{p_1}(T), 0\le i\le l, $
  by  \eqref{KR-2}, \eqref{KAS-2},   the stochastic Gronwall lemma,   and  Lemma \ref{NN},    for
 any $p\in (\ff 1 2,1)$ there exist   constants $c_3,c_4>0$ such that
\beg{align*} &\Big(\E  \Big[\sup_{s\in [0,t\land\tau_j]} |X_s^{(1)}-X_s^{(2)}|^{k}\Big]\Big)^{2} \le c_3(\E \e^{\ff{c_2 p}{1-p} \scr L_t})^{\ff {1-p}p} \E G_j(t)\\
&\le c_4 \big(  |x^{(1)}-x^{(2)}|^{2k}+j^{2(k-1)}\| b^{(0)}-b^{0,n}\|_{\tt L_{q_0}^{p_0}(T)}\big),\ \ n,j\ge 1.\end{align*}
By first letting $n \to\infty$ then $j\to\infty$ and applying \eqref{APP-1}, we prove \eqref{INN} for some constant $c>0$.
\end{proof}

\subsection{Functional inequalities: proof of Theorem \ref{T2.3}}

   Let $\{P_{s,t}\}_{t\ge s\ge 0}$ be the Markov semigroup associated with
\eqref{E01}, i.e.
$P_{s,t}f(x):= \E f(X_{s,t}^x)$ for $t\ge s, f\in \B_b(\bar D),$ 
where $(X_{s,t}^x)_{t\ge s}$ is the unique solution of \eqref{E01} starting from $x$ at time $s$. We have
\beq\label{SMG} P_t f(x)= \E (P_{s,t}f)(X_s^x),\ \ s\in [0,t], f\in C_b^1(\bar D),\end{equation} where $X_s^x:=X_{0,s}^x$.
By \eqref{GRD} for \eqref{E01} from time $s$,  for any  $p>1$,   we have
\beq\label{GRD2} |\nn P_{s,t}f| \le c(p) (P_{s,t}|\nn f|^p)^{\ff 1 p},\ \ 0\le s\le t\le T, f\in C_b^1(\bar D).\end{equation}
If $P_{\cdot,t} f\in C^{1,2}([0,t]\times\bar D)$ for  $f\in C_N^2(\bar D)$ such that
\beq\label{OPP} (\pp_s+L_s)P_{s,t}f=0,\ \ f\in C_N^2(\bar D),\nn_\n P_{s,t}f|_{\pp D}=0,\end{equation}
then the desired inequalities  follow  from \eqref{GRD2} by taking derivative in $s$  to the following reference functions respectively:
 $$P_s\{P_{s,t}(\vv+f)\}^p,\ \ P_s\{P_{s,t}(\vv+f)\}^2,\ \ P_s\{\log P_{s,t}(\vv+f)\}(x+s(y-s)/t),\ \ s\in [0,t],$$
 see for instance the proof of \cite[Theorem 3.1]{14WZ}. However, in the present singular setting it is not clear whether \eqref{OPP} holds or not.
 So, below we make an approximation argument.

(a) Proof of \eqref{GRD'}. Let $\{b^{0,n}\}_{n\ge 1}$ be the mollifying approximations of $b^{(0)}$. By $(A_2^{\si,b})$, for any $f\in C_N^2(\bar D)$ and $t\in (0,T]$, the equation
 $$u_{s,t}^n= P_{s,t}^{\si,b^{(1)}}f+ \int_s^t P_{s,r}^{\si,b^{(1)}}(\nn_{b_r^{0,n}}u_{s,t}^n)\d r,\ \ s\in [0,t].$$
has a unique solution in $C^{1,2}([0,t]\times\bar D)$, and $P_{s,t}^nf:= u_{s,t}^n$ satisfies
\beq\label{LNNM} (\pp_s+L_s^{\si,b^{(1)}} +\nn_{b_s^{0,n}})P_{s,t}^n f=0,\ \  s\in [0,t], f\in C_N^2(\bar D).\end{equation}
By this and It\^o's formula for the SDE
$$\d X_{s,t}^{x,n}= (b_t^{(1)}+b_t^{0,n})(X_{s,t}^{x,n}) \d t+ \si_t(X_{s,t}^{x,n})\d W_t,\ \ t\ge s, X_{s,s}^{x,n}=x,$$ we obtain
$P_{s,t}^n f(x)= \E  f(X_{s,t}^{x,n})$ for $0\le s\le t.$ 
Let $X_t$ solve $\eqref{E01}$ from time $s$ with $X_s=x$, and define
$$R_s:=\e^{\int_0^s \<\xi_r^n, \d W_r\> -\ff 1 2\int_0^s |\xi_r^n|^2\d r},\ \ \xi_s^n:= \big\{\si_s^*(\si_s\si_s^*)^{-1}(b_s^{(0)} - b_t^{0,n})\big\}(X_s),\ \ s\in [0,t].$$
By  Girsanov's theorem, we obtain
\beg{align*} &|P_{s,t}f-P_{s,t}^nf|(x)= |\E[f(X_t)- R_t f(X_t) ]| \\
&\le \|f\|_\infty \big(\E\e^{c\int_0^t|b_s^{(0)}-b_s^{0,n}|^2(X_s)}-1)=:\|f\|_\infty\vv_n,\ \ 0\le s\le t\le T,\end{align*}
where $c>0$ is a constant and due to \eqref{KAS-2}, $\vv_n\to 0$ as $n\to\infty$. Consequently,
\beq\label{UNF} \|P_{s,t}f-P_{s,t}^n f\|_\infty\le \vv_n \|f\|_\infty,\ \ n\ge 1, 0\le s\le t\le T.\end{equation}
Moreover,   the proof of   \eqref{GRD2} implies that it holds for $P_{s,t}^n$ replacing $P_{s,t}$   uniformly in $n\ge 1$, since the constant   is increasing in  $\|b^{(0)}\|_{\tt L_{q_0}^{p_0}(T)}$, which is not less that $   \|b^{0,n}\|_{\tt L_{q_0}^{p_0}}(T).$  Thus,
\beq\label{GRD2'}  |\nn P_{s,t}^nf| \le c(p) (P_{s,t}^n|\nn f|^p)^{\ff 1 p},\ \ 0\le s\le t\le T, f\in C_b^1(\bar D), n\ge 1.\end{equation}
Now, let $0\le f\in C_N^2(\bar D)$ and $t\in (0,T].$ For any $\vv>0$ and $p\in (1,2]$,
by \eqref{GRD2'}, \eqref{LNNM}, \eqref{UNF},  $(A_2^{\si,b})$  and It\^o's formula, we find   constants   $c_1,c_2>0$ such that
\beg{align*} &\d (\vv+ P_{s,t}^nf)^p(X_s) = \big\{p (\vv+P_{s,t}^nf)^{p-1} \<b_t^{(0)}-b_t^{0,n}, \nn P_{s,t}^n f\> \\
 &\qquad\qquad\qquad \qquad + p(p-1) (\vv+P_{s,t}^nf)^{p-2}
|\si_s^* \nn P_{s,t}^n f|^2\big\}(X_s)\d s+\d M_s\\
&\ge\big\{ c_2  (\vv+P_{s,t}^nf)^{p-2}|\nn P_{s,t}^n f|^2  -c_1\|\nn f\|_\infty  |b_t^{(0)}-b_t^{0,n}|\big\}(X_s)\d s+\d M_s,\ \ s\in [0,t],\vv>0\end{align*}
holds for some martingale $M_s$. By \eqref{KR}, H\"older's inequality, and $\|b^{(0)}-b^{0,n}\|_{\tt L_{q_0}^{p_0}(T)}\to 0$ as $n\to\infty$,
we find  a constant $c_3>0$ and  sequence $\vv_n\to 0$ as $n\to\infty$ such that
\beg{align*} &\vv_n+P_t(\vv+f)^p- (P_t^n f+\vv)^p\ge  c_2 \int_0^t P_s\big\{(\vv+P_{s,t}^nf)^{p-2}|\nn P_{s,t}^n f|^2\big\}\d s\\
&\ge c_2 \int_0^t \ff{(P_s|\nn P_{s,t}^n f|^p)^{\ff 2 p}}{\{P_s(\vv+P_{s,t}^n f)^p\}^{\ff{2-p}p}} \d s\ge c_3 \int_0^t \ff{|\nn P_s  P_{s,t}^n f|^2 }{\{P_s(\vv+P_{s,t}^n f)^p\}^{\ff{2-p}p}} \d s,\ \ \vv\in (0,1).\end{align*}
Thus, for any $x\in D$ and $x\ne y\in B(x,\dd)\subset D$ for small $\dd>0$ such that
$$x_r:= x+r(y-x)\in D,\ \ r\in [0,1],$$ this implies
\beg{align*} &\ff{|\int_0^t (P_sP_{s,t}^n f(x)- P_sP_{s,t}^n (y))\d s|}{|x-y|}\le \int_0^{1} \d r\int_0^t |\nn P_sP_{s,t}^n f|(x_r)) \d s\\
&\le \int_0^{1}  \bigg(\int_0^t \ff{|\nn P_sP_{s,t}^n f|^2}{\{P_s(\vv+P_{s,t}^n f)^p\}^{\ff {2-p}p}}(x_r)\d s\bigg)^{\ff 1 2}\bigg(\int_0^t \{P_s(\vv+P_{s,t}^nf)^p)\}^{\ff{2-p}p}(x_r)\d s\bigg)^{\ff 1 2}\d r\\
&\le \int_0^{1}   c_3^{-1/2}\big\{\vv_n+ P_t(\vv+ f)^p\big\}^{\ff 1 2} (x+r(y-x)) \bigg(\int_0^t (\vv+P_sP_{s,t}^nf^p)^{\ff{2-p}p})(x_r)\d s\bigg)^{\ff 1 2}\d r.\end{align*}
Combining this with \eqref{UNF} and letting $n\to\infty,\vv\to 0$, we obtain
$$\ff{|P_tf(x)-P_tf(y)|}{|x-y|} \le \ff{1}{t} \int_0^{1} (c_3^{-1}P_tf^p)^{\ff 1 2}(x_r) \bigg(\int_0^t (P_tf^p)^{\ff{2-p}p}(x_r)\d s\bigg)^{\ff 1 2}\d r.$$ Letting $y\to x$ we prove \eqref{GRD'} for some constant $c$ depending on $p$, for $p\in (1,2]$ and all $f\in C_N^2(\bar D)$.
By Jensen's inequality the estimate also holds for $p>2$, and by approximation argument, it holds for all $f\in \B_b(\bar D)$.

(b)  Proof of \eqref{PCC}. By  \eqref{GRD2'},
  It\^o's formula and $(A_2^{\si,b})$, we find a constant $c_4>0$ and a martingale $M_s$ such that
\beg{align*}  &\d (P_{s,t}^n f)^2(X_s) = 2 \big\{\<\nn P_{s,t}^n f, b_s^{(0)}-b_s^{0,n}\> + |\si_s^*\nn P_{s,t}^n f|^2\big\}(X_s) \d s+\d M_s\\
&\le c_4\big\{ \|\nn f\|_\infty |b_s^{(0)}-b_s^{0,n}|+ P_{s,t}^n |\nn f|^2\big\}(X_s) \d s +\d M_s,\ \ s\in [0,t].\end{align*}
Integrating both sides over $s\in [0,t]$, taking expectations and letting $n\to\infty$, and combining with  \eqref{KR} and \eqref{UNF}, we prove \eqref{PCC}.

(c) Proof of \eqref{LHA}.  Let $0<f\in C_N^2(\bar D)$. By taking It\^o's formula to $P_{s,t}^n(\vv+ f)(X_s)$ for $\vv>0$ and taking expectation, we derive
$$\ff{\d}{\d s} P_s \log P_{s,t}^n\{\vv+f\}= -P_s|\si_s^*\nn \log P_{s,t}^n f|^2 + P_s\<b_s^{(0)}-b_s^{0,n}, \nn\log P_{s,t}^n(\vv+f)\>.$$
For any $x,y\in \bar D$, let $\gg: [0,1]\to \bar D$ be a curve linking $x$ and $y$ such that $|\dot \gg_r|\le c |x-y|$ for some constant $c>0$ independent of
$x,y$. Combining these with $(A_2^{\si,b})$ and \eqref{GRD} for $p=2$ we find a constant $c_5>0$ such that
 \beg{align*} &P_t\log \{\vv +f\}(x)- \log P_t^n\{\vv+ f\}(y) =\int_0^t \ff{\d}{\d s} P_s\log P_{s,t}^n f(\gg_{s/t})\d s \\
 &\le \int_0^t \big\{ct^{-1}  |x-y| |\nn P_s \log P_{s,t}^n f(\gg_{s/t})| -P_s  |\si_s^*\nn\log P_{s,t}^n f|^2\big\}(\gg_{s/t}) \d s\\
 &\le c_5 \int_0^t  \ff{|x-y|^2}{ t^2}\d s=\ff{c_5|x-y|^2}t,\ \ t\in (0,T].\end{align*} Therefore, \eqref{LHA} holds.

\section{Well-posedness for    DDRSDEs}

To characterize  the dependence on the distribution, we will use different  probability distances. For a measurable function
$$\psi: \bar D\times \bar D\to [0,\infty)\ \text{with}\   \psi(x,y)=0  \ \text{if\ and\ only\ if\ }x=y,$$
we introduce the associated  Wasserstein $``$distance'' (also called transportation cost)
$$\W_\psi (\mu,\nu):=   \inf_{\pi\in\C(\mu,\nu)}  \int_{\bar D\times\bar D} \psi(x,y) \pi(\d x,\d y),\ \ \mu,\nu\in \scr P (\bar D),$$
 where $\C(\mu,\nu)$ is the set of all couplings for $\mu$ and $\nu$. In general, $\W_\psi$ is not necessarily a distance as it may be infinite and the triangle inequality may not hold. In particular,
 when $\psi(x,y)=|x-y|^k$ for some constant $k> 0$, the $L^k$-Wasserstein distance $\W_k:= (\W_\psi)^{\ff 1 {1\lor k}}$ is a complete metric on the space
 $$\scr P_k(\bar D):=  \big\{\mu\in \scr P(\bar D):\ \|\mu\|_k:=\mu(|\cdot|^k)^{\ff 1 k}<\infty\big\},$$ 
 where $\mu(f):=\int f\d\mu$ for $f\in L^1(\mu)$. 
 When $k=0$ we set $\|\mu\|_0=1$ such that $\scr P_2(\bar D)= \scr P(\bar D)$ and $\W_0$ reduces to the total variation norm
 $$\W_0(\mu,\nu)=\ff 1 2 \|\mu-\nu\|_{var}:=\ff 1 2 \sup_{|f|\le 1} |\mu(f)-\nu(f)|=\sup_{A\in\scr B(\bar D)}|\mu(A)-\nu(A)|,$$
 where $\B(\bar D)$ is the Borel $\si$-algebra of $\bar D$.
 We will also use the weighted variation norm for $k>0$:
 $$\|\mu-\nu\|_{k,var} := \sup_{|f|\le 1 + |\cdot|^k} |\mu(f)-\nu(f)|,\ \ \mu,\nu\in \scr P_k(\bar D).$$
  According to \cite[Theorem 6.15]{VV},   there exists a constant $c>0$ such that
\beq\label{001} \|\mu-\nu\|_{var}+\W_k(\mu,\nu)^{1\lor k} \le c \|\mu-\nu\|_{k,var},\ \ \mu,\nu\in \scr P_k(\bar D).\end{equation}
However, when $k>1$, for any constant $c>0$,  $\W_k(\mu,\nu)\le c \|\mu-\nu\|_{k,var}$ does not hold. Indeed, by taking
$$\mu=\dd_0, \ \ \nu= (1-n^{-1-k})\dd_0 + n^{-1-k} \dd_{n e},\ \ n\ge 1, e\in\R^d \ \text{with}\ |e|=1,$$
we have $\W_k(\mu,\nu)= n^{-\ff 1 k}$, while
$$\|\mu-\nu\|_{k,var} =n^{-1-k} \|\dd_0-\dd_{ne}\|_{k,var}\le n^{-1-k} \big\{\dd_0(1+|\cdot|^k) + \dd_{ne}(1+|\cdot|^k)\big\}\le \ff 3 n,\ \ n\ge 1,$$
so that $\lim_{n\to\infty} \ff{\W_k(\mu,\nu)}{\|\mu-\nu\|_{k,var}}=\infty$ for $k>1$.

In Theorem \ref{T1} below, we use the enlarged probability distance $\|\cdot\|_{k,var}+ \W_k$ to measure the distribution dependence of the  DDRSDE \eqref{E1}.
For any subspace $\hat {\scr P}$ of $\scr P(\bar D)$ and any $T\in (0,\infty]$, let $C([0,T]; \hat {\scr P})$ be the set of
all continuous maps from $[0,T]\cap [0,\infty)$ to $\hat {\scr P}$ under the weak topology.
For any  $\mu\in C([0,\infty);\scr P(\bar D)),$ let $\si^\mu$ and $b^\mu$ be in \eqref{SBB}.

\subsection{Singular case}

  We make the following assumption. Recall that $b^\mu_t:=b_t(\cdot,\mu_t)$ for $\mu\in C([0,\infty);\scr P(\bar D))$.

\emph{ \beg{enumerate} \item[{\bf (A1)}]  Let $T>0$ and $k\ge 0$.     $\si^\mu=\si$ does not depend on $\mu$, and there exists $\hat\mu\in\scr P_k(\bar D)$ such that
at least {\bf one of the following two conditions} holds.
\item[$(1)$]   $(A_2^{\si,\hat b})$ holds for  $\hat b:= b(\cdot,\hat\mu)$, and  there exist  a constant $\aa\ge 0 $ and      $1\le f_i\in \tt L_{q_i}^{p_i}(T,D), 0\le i\le l,$   such that  for any $t\in [0,T],$
 $ x\in \bar D,$ and $\mu,\nu\in \scr P_k(\bar D),$
 \beq\label{XI11} |b_t^\mu (x)-\hat b_t^{(1)}(x)|\le f_0(t,x)  +\aa\|\mu\|_k,  \end{equation}
\beq\label{XI1}   |b_t^\mu(x)-b_t^\nu(x)|
 \le      \big\{\|\mu-\nu\|_{k,var}+ \W_k(\mu,\nu)\big\}\sum_{i=0}^lf_i(t,x).\end{equation}
 \item[$(2)$] $(A_1^{\si,\hat b})$ holds,  and    \eqref{XI11}-\eqref{XI1}  holds for  
    $|f_i|^2\in \sup_{(p,q)\in\scr K} \tt L_q^p(T,D), 0\le i\le l$.
      \end{enumerate}}

Since $\hat b_t^{(1)} $ is regular,  \eqref{XI11}  gives a control for the singular term  of $b^\mu$. Moreover, \eqref{XI1} is a Lipschitz condition on $b_t(x,\cdot)$ in
$\|\cdot \|_{k,var}+ \W_k$ with a singular Lipschitz coefficient.

 \beg{thm}\label{T1}  Assume   {\bf (A1)}.
  \beg{enumerate} \item[$(1)$]    $\eqref{E1}$ is weak well-posed up to time $T$ for distributions in $\scr P_k(\bar D)$. Moreover,   for any $\gg\in \scr P_k(\bar D)$, and any  $ n>0$,   there exists a constant $c>0$, such that
  \beq\label{ESS} \E\Big[ \sup_{t\in [0,T]}  |X_t|^{n }\Big|X_0 \Big]\le c(1+  |X_0|^{n }),\ \ \E   \e^{n l_T}\le c\end{equation}
  holds for the solution with $\L_{X_0}=\gg$.
  \item[$(2)$]  $\eqref{E1}$ is  well-posed up to time $T$ for distributions in $\scr P_k(\bar D)$ in each of the following situations:

  $(i)$   $d=1$ and {\bf (A1)}$(2) $  holds.

   $(ii)$    {\bf (A1)}$(1) $ holds with $p_1>2$ in  $(A_2^{\si,\hat b})$.   \end{enumerate}
  \end{thm}
To prove Theorem \ref{T1}, we first present a general  result on the well-posedness of the DDRSDE \eqref{E1}  by  using  that of  the  reflecting SDE  \eqref{E01}.

For any $k\ge 0, \gg\in \scr P_k, N\ge 2$, let
   $$\scr P_{k,\gg}^{T,N} = \Big\{\mu\in C([0,T]; \scr P_k(\bar D)): \mu_0=\gg, \sup_{t\in [0,T]} \e^{-Nt} (1+\mu_t(|\cdot|^k)) \le N\Big\}. $$
   Then as $N\uparrow\infty$,
  \beq\label{NPP} \scr P_{k,\gg}^{T,N}\uparrow \scr P_{k,\gg}^T:= \big\{\mu\in C([0,T];\scr P_k(\bar D)): \mu_0=\gg \big\}.\end{equation}
   For any    $\mu\in \scr P_{k,\gg}^{T},$ we will assume that the reflecting SDE
    \beq\label{XM} \d X_t^{\mu,\gg}= b_t (X_t^{\mu,\gg},\mu_t)\d t+ \si_t (X_t^{\mu,\gg}) \d W_t +\n(X_t^{\mu,\gg})\d l_t^{\mu,\gg},\ \ t\in [0,T], \L_{X_0^{\mu,\gg}}=\gg \end{equation}
      has a unique weak solution with
   $$ H^\gg_t(\mu):= \L_{X_t^{\mu,\gg}}\in \scr P_k(\bar D),\ \ t\in [0,T].$$

   \beg{enumerate} \item[{\bf (H2)}] Let $k\ge 0, T>0$. For any $\gg\in \scr P_k(\bar D)$ and   $\mu\in \scr P_{k,\gg}^{T},$
   \eqref{XM} has a unique weak solution,  and   there exist    constants $\{(p_i',q_i')>1\}_{ 0\le i\le l}, N_0\ge 2$ and  increasing maps
   $C: [N_0,\infty)\to (0,\infty)$ and $ F:[N_0,\infty)\times [0,\infty)\to (0,\infty)  $
    such that for any $N\ge N_0$ and $\mu\in \scr P_{k,\gg}^{T,N},$ the (weak) solution  satisfies
              \beq\label{FFG} H^\gg(\mu):=\L_{(X^{\mu,\gg}_t)_{t\in [0,T]}}\in \scr P_{k,\gg}^{T, N},\end{equation}
       \beq\label{PFF} \big(\E \big[(1+|X_t^{\mu,\gg}|^k)^2\big|X_0^{\mu,\gg}\big]\big)^{\ff 1 2} \le C(N) (1+|X_0^{\mu,\gg}|^k),\ \ t\in [0,T], \end{equation}
           \beq\label{RPP}\beg{split} & \E \bigg(\int_{0}^{t} g_s(X_s^{\mu,\gg})\d s\bigg)^2  \le C(N) \|g\|_{\tt L_{q_i'}^{p_i'}(t_0,t_1)}^2,\\
            &\E \e^{\int_{0}^{t} g_s(X_s^{\mu,\gg})\d s}  \le F(N, \|g\|_{\tt L_{q_i}^{p_i}(t,D)}),\ \ t\in [0,T], g\in \tt L_{q_i'}^{p_i'}(t,D),  0\le i\le l.\end{split} \end{equation}
          \end{enumerate}
Obviously,  when $k=0$,   conditions \eqref{FFG}  and \eqref{PFF}      hold for   $N_0=2.$

    \beg{thm}\label{T3.1}  Assume  {\bf (H2)} and let $\si^\mu=\si$ do  not depend on $\mu$.
   Assume that     there exist a    measurable map
    $ \GG: [0,T]\times  \bar D\times\scr P(\bar D)\to\R^m$ such that
    \beq\label{XI0}    b_t(x,\nu)-b_t(x,\mu)=\si_t(x)\GG_t(x,\nu,\mu),\  \     x\in\bar D, t\in [0,T],  \nu,\mu\in \scr P_k(\bar D).  \end{equation}
    Let     $f:=\big(\sum_{i=0}^l \tt f_i \big)^{\ff 1 2} $ for some  $ 1\le \tt f_i \in \tt L_{q_i'}^{p_i'}(T),  0\le i\le l.$   
     \beg{enumerate} \item[$(1)$] If  
\beq\label{XI0-1}   |\GG_t(x,\nu,\mu)| \le  f_t(x)  \|\nu-\mu\|_{k,var}, \ \  x\in\bar D, t\in [0,T],  \nu,\mu\in \scr P_k(\bar D),\end{equation}
    Then $\eqref{E1}$ is  weak    well-posed up to time $T$ for distributions in $\scr P_k(\bar D)$.  If, furthermore, in {\bf (H2)} the SDE \eqref{XM} is strongly well-posed for any $\gg\in \scr P_k(\bar D)$ and $\mu\in \scr P_{k,\gg}^T$, so is $\eqref{E1}$ up to time $T$ for distributions in $\scr P_k(\bar D)$.
     \item[$(2)$] Let $k>1$ and for any  $\mu,\nu\in \scr P_k(\bar D),$
\beq\label{XI0-2}   |\GG_t(x,\nu,\mu)| \le   f_t(x)  \big\{\|\nu-\mu\|_{k,var} +\W_k(\mu,\nu)\big\}, \ \  (t, x)\in [0,T]\times\bar D.\end{equation}
If for any $\gg\in \scr P_k(\bar D)$ and $N\ge N_0$, there exists a constant $C(N)>0$ such that for any $\mu,\nu\in \scr P_{k,\gg}^{T,N},$
\beq\label{XI0-3} \W_k(H^\gg_t(\mu), H_t^\gg(\nu))^{2k}\le C(N) \int_0^t \big\{\|\mu_s-\nu_s\|_{k,var}^{2k}+\W_k(\mu_s,\nu_s)^{2k}\big\}\d s,\ \ t\in [0,T],\end{equation}
then assertions in $(1)$ holds. \end{enumerate}
   \end{thm}

 \beg{proof}  Let $\gg\in \scr P_k(\bar D)$. Then the weak solution to \eqref{XM} is a weak solution to \eqref{E1} if and only if $\mu$ is a fixed point of the map $H^\gg$ in $\scr P_{k, \gg}^T$. So, if $H^\gg$ on $\scr P_{k, \gg}^T$ has a unique fixed point in $\scr P_{k, \gg}^T$,
 then the (weak) well-posedness of  \eqref{XM} implies that of   \eqref{E1}.
 Thus, by \eqref{NPP}, it suffices to show that for any $N\ge N_0$, $H^\gg$ has a unique fixed point in $\scr P_{k,\gg}^{T,N}$. By \eqref{FFG} and the fixed point theorem, we only need to prove that for any $N\ge N_0$, $H^\gg$ is contractive with respect to a complete metric on $ \scr P_{k,\gg}^{T,N}.$

 (1)   For any $\ll>0$, consider the metric
 $$\W_{k,\ll,var}(\mu,\nu):= \sup_{t\in [0,T]} \e^{-\ll t} \|\mu_t-\nu_t\|_{k,var},\ \ \mu,\nu\in \scr P_{k,\gg}^{T,N}.$$
Let    $ (X_t^{\mu,\gg}, l_t^{\mu,\gg})$ solve \eqref{XM} for some Brownian motion $W_t$ on a  complete probability filtration space $(\OO,\{\F_t\},\P)$. 
By \eqref{RPP}, \eqref{XI0-1} or \eqref{XI0-2},  we find a constant $c_1>0$ depending on $N$  such that
\beq\label{RPP-1} \beg{split} &\sup_{\mu,\nu\in \scr P_{k,\gg}^{T,N}  } \E \big(\e^{2\int_0^T |\GG_s(X_s^{\mu,\gg},\nu_s,\mu_s)|^2 \d s}|\F_0\big)\le c_1^2,\\
& \sup_{\mu\in \scr P_{k,\gg}^{T,N}  } \E \bigg(\bigg(\int_0^T g_s(X_s^{\mu,\gg})  \d s\bigg)^2\bigg|\F_0\bigg)
\le c_1^2\|g\|_{\tt L_{q_i}^{p_i}(T)}^2,\ \ g\in \tt L_{q_i}^{p_i}(T),\  0\le i\le l.\end{split}\end{equation}
 Then by Girsanov's theorem,
$$\tt W_t:= W_t- \int_0^t \GG_s(X_s^{\mu,\gg},\nu_s,\mu_s)\d s,\ \ t\in [0,T]$$ is a Brownian motion under the probability $\Q:=R_T \P$, where
$$R_t:=\e^{\int_0^t \<\GG_s(X_s^{\mu,\gg},\nu_s,\mu_s),\d W_s\>-\ff 1 2 \int_0^t |\GG_s(X_s^{\mu,\gg},\nu_s,\mu_s)|^2\d s  },\ \ t\in [0,T]$$ is a $\P$-martingale. By \eqref{XI0}, we may
formulate  \eqref{XM} as
$$\d X_t^{\mu,\gg}= b_t (X_t^{\mu,\gg},\nu_t)\d t+ \si_t (X_t^{\mu,\gg}) \d \tt W_t +\n(X_t^{\mu,\gg})\d l_t^{\mu,\gg},\ \ t\in [0,T], \L_{X_0^{\mu,\gg}}=\gg.$$
By the weak uniqueness due to {\bf (H2)},  the definition of $\|\cdot\|_{k,var}, $    \eqref{PFF} and \eqref{XI0},   we obtain
\beq\label{RPP-2} \beg{split} &\|H_t^\gg(\mu)-H_t^\gg(\nu)\|_{k,var}
=\sup_{|\tt f|\le 1+|\cdot|^k}\big|\E\big[(R_t-1) \tt f(X_t^{\mu,\gg})\big]\big|\\
&\le \E \big[(1+|X_t^{\mu,\gg}|^k)|R_t-1|\big]
 \le \E \Big[\big\{\E\big((1+|X_t^{\mu,\gg}|^k)^2|\F_0\big)\big\}^{\ff 1 2} \big\{\E\big(|R_t-1|^2|\F_0\big)\big\}^{\ff 1 2} \Big]\\
& \le C(N) \E\Big[(1+|X_0^{\mu,\gg}|^k) \big\{\E(\e^{\int_0^t |\GG_s(X_s^{\mu,\gg},\nu_s,\mu_s)|^2\d s}-1|\F_0)\big\}^{\ff 1 2}\Big].\end{split}\end{equation}
Moreover,  \eqref{RPP-1} implies
\beg{align*} &\E(\e^{\int_0^t |\GG_s(X_s^{\mu,\gg},\nu_s,\mu_s)|^2\d s}-1|\F_0)\\
&\le \E\bigg(\e^{\int_0^t |\GG_s(X_s^{\mu,\gg},\nu_s,\mu_s)|^2\d s}\int_0^t |\GG_s(X_s^{\mu,\gg},\nu_s,\mu_s)|^2\d s\bigg|\F_0\bigg)\\
&\le c_1 \bigg\{\E\bigg( \bigg(\int_0^t| f_s(X_s^{\mu,\gg})|^2 \|\mu_s-\nu_s\|_{k,var}^2\d s\bigg)^2\bigg|\F_0\bigg)\bigg\}^{\ff 1 2}\\
&\le c_1 \e^{2\ll t} \W_{k,\ll, var}(\mu,\nu)^2 \bigg\{\E\bigg( \bigg(\int_0^t |f_s(X_s^{\mu,\gg}) |^2\e^{-2\ll(t-s)}\d s\bigg)^2\bigg|\F_0\bigg)\bigg\}^{\ff 1 2}\\
&\le c_1^2 \e^{2\ll t}\sum_{i=0}^l \|\tt f_i\e^{-2\ll(t-\cdot)}\|_{\tt L_{q_i'}^{p_i'}(t)}\W_{k,\ll, var}(\mu,\nu)^2,\ \ t\in [0,T].\end{align*}
  Combining this with \eqref{RPP-2} and the definition of $\W_{k,\ll,var}$, we obtain
\beq\label{HX10} \W_{k,\ll,var}(H^\gg(\mu),H^\gg(\nu))\le C(N) (1+ \gg(|\cdot|^k))c_1\ss{\vv(\ll)} \W_{k,\ll,var}(\mu,\nu), \ \ \ll>0,\end{equation}
where
$$\vv(\ll):= \sup_{t\in [0,T]} \sum_{i=0}^l \|\tt f_i\e^{-2\ll(t-\cdot)}\|_{\tt L_{q_i'}^{p_i'}(t)}\downarrow 0\ \ \text{as}\ \ \ll\uparrow \infty.$$
 So,   $H^\gg$ is contractive on $(\scr P_{k,\gg}^{T,N}, \W_{k,\ll,var})$ for large enough $\ll>0$.

 (2) Let $k>1$. We consider the metric $\tt \W_{k,\ll,var}:= \W_{k,\ll,var}+ \W_{k,\ll}$, where
 $$\W_{k,\ll}(\mu,\nu):= \sup_{t\in [0,T]} \e^{-\ll t}\W_k(\mu_t,\nu_t),\ \ \mu,\nu\in \scr P_{k,\gg}^{T,N}.$$
 By using \eqref{XI0-2} replacing \eqref{XI0-1},  instead of \eqref{HX10} we find   constants $\{C(N,\ll)>0\}_{\ll>0}$ with $C(N,\ll)\to 0$ as $\ll\to\infty$ such that
 \beq\label{HX1} \W_{k,\ll,var}(H^\gg(\mu),H^\gg(\nu))\le C(N,\ll)  \tt \W_{k,\ll,var}(\mu,\nu), \ \ \ll>0, \mu,\nu\in \scr P_{k,\gg}^{T,N}.\end{equation}
 On the other hand,  \eqref{XI0-3} yields
\beg{align*}  & \W_{k,\ll} (H^\gg(\mu), H^\gg(\nu))\le \sup_{t\in [0,T]} \bigg(C(N) \e^{-\ll k t} \int_0^t \big\{\|\mu_s-\nu_s\|_{k,var}^{2k}+\W_k(\mu_s,\nu_s)^{2k}\big\}\d s\bigg)^{\ff 1 {2k}} \\
 &\le \tt \W_{k,\ll,var} (\mu,\nu) \sup_{t\in [0,T]}  \bigg(C(N)\int_0^t \e^{-2\ll k(t-s)}\d s \bigg)^{\ff 1 {2k}} \le  \ff{C(N)^{\ff 1 {2k}}}{(2\ll k)^{\ff 1 {2k}}} \tt \W_{k,\ll,var} (\mu,\nu),\ \ \ll>0.\end{align*}
 Combining this with \eqref{HX1}, we concluded that $H^\gg$ is contractive in $ \scr P_{k,\gg}^{T,N}$ under the metric $\tt \W_{k,\ll,var}$ when $\ll$ is large enough, and hence finish the proof.
  \end{proof}

\beg{proof}[Proof of Theorem \ref{T1}]  Let $\gg\in \scr P_k(\bar D)$ be fixed. By \eqref{XI11}, for any $i=1,2$,  condition  $(A_i^{\si,\hat b})$ implies $(A_i^{\si,b^\mu})$ for any $\mu\in C([0,\infty);\scr P_k(\bar D))$. So,
by Theorem \ref{T2.1}, {\bf (A1)} implies the weak well-posedness of \eqref{XM} for distributions in $\scr P_k(\bar D)$ with
\beq\label{*BC} H_t^\gg(\mu)\in \scr P_k(\bar D),\ \ \ \E \e^{\ll l_T^{\mu,\gg}}<\infty,\ \ \ll>0, \gg\in \scr P_k(\bar D), \mu\in C([0,\infty);\scr P_k(\bar D)),\end{equation}
 and also implies  the  strong well-posedness of    \eqref{XM} in each situation of Theorem \ref{T1}(2). Moreover,  by  Lemma \ref{L1} and Lemma \ref{L1'},   {\bf (A1)} implies that   \eqref{RPP} holds for any $(p,q) \in \scr K$,
 as well as for $(p,q)= (p_0/2, q_0/2)$ under   $(A_2^{\si,\hat b}),$
 \eqref{XI0} with \eqref{XI0-1} holds  for $k\le 1$ due to \eqref{001},  and    \eqref{XI0} with \eqref{XI0-2} holds for $k>1$.
Therefore, by Theorem \ref{T3.1}, it remains to verify \eqref{ESS}, \eqref{FFG}, \eqref{PFF}, and \eqref{XI0-3} for $k>1$.
 Since \eqref{PFF} and \eqref{FFG} are trivial for $k=0$,
  we only need to prove:     \eqref{ESS}, 
  \eqref{PFF} and \eqref{FFG} for $k>0$, 
   \eqref{XI0-3} for $k>1$ for case $(i)$, 
and   \eqref{XI0-3} for $k>1$ for case $(ii)$.
 

(a)  Simply denote 
 $$f_t(x):=\sum_{i=0}^l f_i(t,x).$$  We first prove that under {\bf (A1)}, there exits a constant $c>0$ and an increasing function $c: [1,\infty)\to(0,\infty)$ such that for any $j\ge 1$ and $\mu\in \scr P_{k,\gg}^T$,
\beq\label{WY1} \beg{split} &\E\bigg(\int_0^t |f_s(X_s^{\mu,\gg})|^2\d s\bigg)^{j} \le c(j)+c(j) \bigg(\int_0^t \|\mu_s\|_k^2\d s\bigg)^j,\\
&\E \exp\bigg[ j\int_0^t |f_s(X_s^{\mu,\gg})|^2\d s\bigg]  \le c(j) \exp\bigg[ c  \int_0^t \|\mu_s\|_k^2\d s\bigg],\ \ t\in [0,T],\end{split}\end{equation}
where $X_t^{\mu,\gg}$ solves \eqref{XM}.
We will prove these estimates by Lemmas \ref{L1} and \ref{L1'} for the following   reflecting SDE:
$$\d \hat X_s= \hat b_s(\hat X_s)\d s +\si_s(\hat X_s)\d W_s +\n(\hat X_s)\d \hat l_s,\ \ \hat X_{0}= X_{0}^{\mu,\gg}, s\in [0,t].$$
By \eqref{KAS-2} under {\bf (A1)}(1), and  \eqref{KAS} under {\bf (A1)}(2), for any $j\ge 1$  we find a constant $c_1(j)>0$ such that
\beq\label{GP1}  \E \e^{ j\int_0^t (|\hat b_s^{(0)}|^2+|f_s|^2)(\hat X_s) \d s}   \le c_1(j),\ \ \ t\in [0,T].\end{equation}
 Let $\gg_s= \big\{[\si_s^*(\si_s\si_s^*)^{-1}](b^\mu_s-\hat b_s)\big\}(\hat X_s),$ and
\beg{align*} R_t:= \e^{\int_{0}^{t}\<\gg_s,\d W_s\>-\ff 1 2 \int_{0}^{t}|\gg_s|^2\d s},\ \  \tt W_s:= W_s-\int_{0}^s \gg_r\d r,\ \ s\in [0,t].
 \end{align*}
By Girsanov's theorem, $(\tt W_s)_{s\in [0,t]}$ is a Brownian motion under $R_t\P$, and the SDE for $\hat X_s$ becomes
$$\d \hat X_s=  b_s^\mu (\hat X_s)\d s +\si_s(\hat X_s)\d \tt W_s +\n(\hat X_s)\d \hat l_s,\ \ \hat X_{0}= X_{0}^{\mu,\gg}, s\in [0,t].$$
So, by \eqref{XI11}, \eqref{GP1} and H\"older's inequality,
we find constants $c_1, c, c(j)>0$ such that
\beg{align*} &\E \e^{j\int_{0}^{t} |f_s(X_s^{\mu,\gg})|^2 \d s }= \E\big[R_t \e^{j\int_{0}^{t} |f_s(\hat X_s)|^2 \d s }\big]
\le \big( \E\e^{2j\int_{0}^{t} |f_s(\hat X_s)|^2 \d s}  \big)^{\ff 1 2}\big(\E[R_t^2]\big)^{\ff 1 2}  \\
&\le\ss{ c_1(2j)} \big(\E\e^{c_1 \int_{0}^{t} \{|\hat b^{(0)}_s|^2 + (f_s+\aa\|\mu_s\|_k)^2\}(\hat X_s) \d s}\big)^{\ff 1 2}
\le c(j) \e^{c\int_0^t \|\mu_s\|_k^2\d s}. \end{align*}
Next, taking $c_2(j)>0$ large enough such that the function
$r\mapsto [\log (r+c_2(j))]^j$ is concave for $r\ge 0$, so that this and Jensen's inequality imply
\beg{align*} &\E\bigg(\int_0^t |f_s(X_s^{\mu,\gg})|^2\d s\bigg)^{j} \le \E \big( \big[\log (c_2(j)+  \e^{\int_{0}^{t} |f_s(X_s^{\mu,\gg})|^2 \d s })\big]^j\big)\\
&\le \big[\log (c_2(j)+\E \e^{\int_{0}^{t} |f_s(X_s^{\mu,\gg})|^2 \d s })\big]^j
\le c(j) + c(j) \bigg(\int_0^t\|\mu_s\|_k^2\d s\bigg)^j\end{align*}
holds for some constant $c(j)>0$. Therefore, \eqref{WY1} holds.

(b) Proof of  \eqref{FFG}. Simply denote $X_t=X_t^{\mu,\gg}$. By \eqref{XI11},  the boundedness of $\si$ and
the condition on $\hat b^{(1)} $ in $(A_0^{\si,\hat b})$ which follows from $(A_2^{\si,b})$ due to Lemma \ref{LNN}, we find a constant $c_1>0$ such that
$$L_{t,\mu}:= \ff 1 2 {\rm tr}\{\si_t\si_t^*\nn^2\} +\nn_{b_t^\mu},\ \ L^{\si,\hat b^{(1)}}  :=\ff 1 2 {\rm tr}\{\si_t\si_t^*\nn^2\} + \nn_{\hat b^{(1)}_t}$$ satisfy
$$L_{t,\mu}\tt\rr\ge  L_{t}^{\si,\hat b^{(1)}}  \tt\rr -|b_t^\mu   - \hat b_t^{(1)} |\cdot|\nn\tt\rr| \ge -c_1 (f_t+  \|\mu_t\|_k).$$
Since  $\<\n,\tt\rr\>|_{\pp D}\ge 1$,  by It\^o's formula  we obtain
\beq\label{*BC2} \d \tt\rr (X_t)\ge -c_1 \big\{f_t(X_t)  +\|\mu_t\|_k \big\}\d t +\d M_t +\d l_t\end{equation}
 for some martingale $M_t$ with $\<M\>_t\le c t$ for some constant $c>0$.
This together with \eqref{WY1} yields that for some constant  $k_0>0$,
$$\E l_t^k\le k_0 +k_0 \E\bigg(\int_0^t\{f_s(X_s)+\|\mu_s\|_k\}\d s \bigg)^k.$$
Combining this with    \eqref{KR},    \eqref{XI1}, \eqref{WY1}  and $\|\si\|_\infty<\infty$, and using the formula
 $$X_t=X_0+\int_0^t b_s^\mu(X_s)\d s + \int_0^t \si_s(X_s)\d W_s+\n(X_t)\d l_t,\ \ \L_{X_0}=\gg,$$
 we find constants $k_1,k_2>0$ such that
\beq\label{*Y2}\beg{split} & \E (1+|X_t|^k)  \le  k_1(1+\|\gg\|_k^k) + k_1\E\bigg(\int_0^t\big\{ |X_s|+ |f_s(X_s)|+  \|\mu_s\|_k \big\}\d s\bigg)^k \\
&\le k_2    + k_2   \E\bigg(\int_0^t  \big\{|X_s|^2+ \|\mu_s\|_k^2\big\}\d s\bigg)^{\ff k 2},\  \ t\in [0,T]. \end{split}\end{equation}

(b1) When $k\ge 2$, by \eqref{*Y2} we find a constant $k_3>0$ such that
$$\E (1+|X_t|^k)\le  k_2 + k_3 \int_0^t \big\{\E |X_s|^k+ \|\mu_s\|_k^k\big\}\d s,\ \ t\in [0,T].$$
By Gronwall's lemma, and noting that $\mu\in \scr P_{k,\gg}^{T,N}$, we find constant $k_4>0$ such that
\beg{align*}  \E (1+|X_t|^k) \le k_4 + k_4 \int_0^t (1+\|\mu_s\|_k^k)\d s
 \le k_4 + k_4 N \e^{Nt} \int_0^t \e^{-N(t-s)}\d s
\le 2k_4 \e^{Nt}, \ \ t\in [0,T].\end{align*}
Taking $N_0= 2 k_4$ we prove
$$\sup_{t\in [0,T]} \e^{-Nt} (1+ \|H_t(\mu)\|_k^k)=\sup_{t\in [0,T]} \e^{-Nt} \E(1+|X_t|^k)\le N_0\le N,\ \ N\ge N_0, \mu\in \scr P_{k,\gg}^{T,N},$$
so that \eqref{FFG} holds.

(b2) When $k\in (0,2)$, by BDG's inequality, and by the same reason    leading to \eqref{*Y2}, we find constants $k_5,k_6, k_7>0$ such that
\beg{align*} & U_t:=\E\Big[\sup_{s\in [0,t]} (1+|X_s|^k)\Big] \le k_5+ k_5 \E\bigg(\int_0^t \big\{|X_s|^2 +   \|\mu_s\|_k^2\big\}\d s\bigg)^{\ff k 2}\\
&\le k_6 + k_6 \E\bigg\{ \Big[\sup_{s\in [0,t]} |X_s|^k\Big]^{1-\ff k 2} \bigg(\int_0^t  |X_s|^k \d s\bigg)^{\ff k 2}\bigg\}+ k_6 \bigg(\int_0^t \|\mu_s\|_k^2\d s\bigg)^{\ff k 2}\\
&\le k_6 + \ff 1 2 U_t + k_7\int_0^t  U_s \d s+ k_6 \bigg(\int_0^t \|\mu_s\|_k^2\d s\bigg)^{\ff k 2},\ \ t\in [0,T].\end{align*}
By Gronwall's lemma, we find   constants $k_8, k_9>0$ such that for any $\mu\in \scr P_{k,\gg}^{T,N}$,
\beg{align*}  \E (1+|X_t|^k) &\le  U_t\le k_8+ k_8 \bigg(\int_0^t  \|\mu_s\|_k^2\d s\bigg)^{\ff k 2}\\
&\le k_8 +k_8N \e^{Nt} \bigg(\int_0^t \e^{-2N(t-s)/k}\d s\bigg)^{\ff k 2} \le k_8+k_9 N^{1-\ff k 2}\e^{Nt},\ \ t\in [0,T].\end{align*}
Thus, there exists  $N_0>0$  such that for any $N\ge N_0$,
$$\sup_{t\in [0,T]} \e^{-Nt} (1+ \|H_t(\mu)\|_k^k)=\sup_{t\in [0,T]} \e^{-Nt} \E(1+|X_t|^k)\le k_8+k_9 N^{1-\ff k 2} \le N,\ \   \mu\in \scr P_{k,\gg}^{T,N},$$
which implies \eqref{FFG}.

(c) Proofs of  \eqref{PFF} and  \eqref{ESS}.  Simply denote $(\hat X_t, \hat l_t)= (X_t^{\mu,\gg}, l_t^{\mu,\gg})$ in \eqref{XM} for $\mu_t=\hat\mu, t\in [0,T];$ that is,
\beq\label{XX} \d\hat X_t =\hat  b_t(\hat X_t)\d t +\si(\hat X_t) \d W_t +\n(\hat X_t) \d\hat l_t,\ \ \L_{\hat X_0}=\gg.\end{equation}
 By {\bf (A1)} and Theorem \ref{T2.1}, this SDE has a unique weak solution, and for any $n\ge 1$ there exists a constant $c>0$ such that
\beq\label{ESS'} \E\Big[\sup_{t\in [0,T]} |\hat X_t|^n\Big|\hat X_0\Big] \le c (1+|\hat X_0|^n),\ \ \E \e^{n \hat l_T} \le c.\end{equation}
So, by \eqref{XI1}, Lemma \ref{L1}, Lemma \ref{L1'} under $(A_2^{\si, \hat b})$,  and Girsanov's theorem,
$$\tt W_t:= W_t-\int_0^t\{\si_s^*(\si_s\si_s^*)^{-1} \}(\hat X_s)  \big\{b^\mu_s(\hat X_s)- \hat b_s(\hat X_s)\big\}\d s,\ \ t\in [0,T]$$
is a $\Q$-Brownian motion for $\Q:= R_T\P$, where
$$R_T:=\e^{\int_0^T \<\{\si_s^*(\si_s\si_s^*)^{-1} \}(\hat X_s)  \{b^\mu_s(\hat X_s)- \hat b_s(\hat X_s)\}, \d W_s\>
-\ff 1 2\int_0^T|\{\si_s^*(\si_s\si_s^*)^{-1} \}(\hat X_s)  \{b^\mu_s(\hat X_s)- \hat b_s(\hat X_s)\}|^2\d s}.$$
By {\bf (A1)}, \eqref{ESS'},  Lemma \ref{L1} when $|f_i|^2\in\cup_{(p,q)\in \scr K}  \tt L_q^p(T)$, and Lemma \ref{L1'}  when   $(A_2^{\si, \hat b})$ holds,    we find an   increasing function $F$ such that
$$\E(|R_T|^2|\F_0)\le \E (\e^{  \int_0^T |f_s(\hat X_s)|^2\{\|\mu_s-\hat\mu\|_{k,var}+\W_k(\mu_s,\hat\mu)\}^2\d s}|\F_0) \le F(\|\mu\|_{k,T}),$$
where $\|\mu\|_{k,T}:= \sup_{t\in [0,T]} \mu_t(|\cdot|^k).$
Reformulating \eqref{XX}  as
$$\d \hat X_t= b_t^\mu(\hat X_t)\d t + \si_t(\hat X_t) \d\tt W_t+\n(\hat X_t)\d \hat l_t,\ \ \L_{\hat X_0}=\gg,$$
by the weak uniqueness  we have  $\L_{\hat X|\Q}= \L_{X^{\mu,\gg}}$, so that \eqref{ESS'}  with $2n$ replacing $n$ implies
\beg{align*} &\E\Big[\sup_{t\in [0,T]}|X_t^{\mu,\gg}|^n\Big| \F_0\Big] =  \E_\Q\Big[\sup_{t\in [0,T]}|\hat X_t|^n\Big| \F_0\Big]\\
&\le  \bigg(\E \Big[\sup_{t\in [0,T]}|\hat X_t|^{2n}\Big| \F_0\Big]\bigg)^{\ff 1 2} (\E R_T^2|\F_0)^{\ff 1 2}\le c(1+|\hat X_0|^n) F(\|\mu\|_{k,T}).\end{align*}
Since $\sup_{\mu\in \scr P_{k,\gg}^{T,N}} \|\mu\|_{k,T}$ is a finite increasing function of $N$, this   implies \eqref{PFF}.

Finally, since $X_t:=X_t^{\mu,\gg}$ solves \eqref{E1} with initial distribution $\gg$ and $\mu_t=\L_{X_t}$ (i.e. $\mu$ is the fixed point of $H^\gg$),  and since $H^\gg$ has a unique fixed point in
  $\scr P_{k,\gg}^{T,N}$ for some  $N>0$ depending on $\gg$ as proved  in the proof of Theorem \ref{T2.1} using  \eqref{RPP} and \eqref{FFG}, we have
  $\L_{X_\cdot}\in \scr P_{k,\gg}^{T,N}$, and hence  \eqref{ESS} follows from \eqref{EPP}.

  (d)  Proof of \eqref{XI0-3} for $k>1$ in case $(i)$.
   Let $u_t^\ll $ and $ \Theta_t^\ll$ be constructed for $b^\mu$ replacing $b$ in   the proof of Theorem \ref{T2.2} under $(A_1^{\si,b})$ for $d=1$.   Let $X_0^{(1)}=X_0^2$ be $\F_0$-measurable with $\L_{X_0^{(i)}}=\gg, i=1,2$.
   As explained in the beginning in the present proof,   the  following reflecting SDEs are well-posed:
  \beg{align*} & \d X_t^{(1)}= b_t(X_t^{(1)},\mu_t)\d t+\si_t(X_t^{(1)})\d W_t+\n(X_t^{(1)})\d l_t^{(1)},\\
  &\d X_t^{(2)}= b_t(X_t^{(2)},\nu_t)\d t+\si_t(X_t^{(2)})\d W_t+\n(X_t^{(2)})\d l_t^{(2)},\ \ t\in [0,T].\end{align*}
  Then instead of \eqref{YYE},  the processes
   $$Y_t^{(i)}:= \Theta_t^\ll (X_t^{(i)}),\ \ i=1,2$$ satisfy
 \beg{align*} & \d Y_t^{(1)}= B_t (Y_t^{(1)}) \d t+ \Sigma_t (Y_t^{(1)}) \d W_t+ \{1+\nn u_t^{\ll}(X_t^{(1)})\}\n(X_t^{(1)}) \d l_t^{(1)},\\
 &\d Y_t^{(2)}=  B_t (Y_t^{(2)})  \d t+ \Sigma_t (Y_t^{(2)}) \d W_t+ \{1+\nn u_t^{\ll}(X_t^{(2)})\}\n(X_t^{(2)}) \d l_t^{(2)}\\
 &\qquad \qquad + \big\{b_t(X_t^{(2)},\nu_t)-b_t(X_t^{(2)}, \mu_t)\big\}\d t.\end{align*}
By     \eqref{XI1},    $Y_0^{(1)}=Y_0^{(2)}$,       It\^o's formula to $|Y_t^{(1)}-Y_t^{(2)}|^{2k}$ with     this formula replacing \eqref{YYE},
the calculations in the proof of Theorem \ref{T2.2} under $(A_1^{\si,b})$ for $d=1$    yield that when $\ll$ is large enough,
\beg{align*}  & |Y_t^{(1)}-Y_t^{(2)}|^{2k} \le c_1 \int_0^t |Y_s^{(1)}-Y_s^{(2)}|^{2k} \d \scr L_s + M_t\\
&\qquad +c_1\int_0^t |Y_s^{(1)}-Y_s^{(2)}|^{2k-1}  f_s(X_s^{(2)}) \big\{\|\mu_s-\nu_s\|_{k,var}+\W_k(\mu_s,\nu_s)   \big\}\d s \\
&\le c_1 \int_0^t   |Y_s^{(1)}-Y_s^{(2)}|^{2k} \d  \tt {\scr L}_s + c_1  \int_0^t  \big\{\|\mu_s-\nu_s\|_{k,var}+ \W_k(\mu_s,\nu_s)\big\}^{2k} \d s + M_t ,\ \ t\in [0,T]\end{align*}
  holds for some constant $c_1>0$ depending on $N$ uniformly in $\mu\in \scr P_{k,\gg}^{T,N}$,   some martingale $M_t$,  $\scr L_t$  in \eqref{LLT},  and
  $$\tt {\scr L}_t:= \scr L_t + \int_0^t | f_s(X_s^{(2)})|^{\ff {2k}{2k-1}}\d s\le \scr L_t + \int_0^t |  f_s(X_s^{(2)})|^2 \d s.$$
  By the stochastic Gronwall  lemma,   Lemma \ref{L1}, we find a constant $c_2>0$ depending on $N$ such that
  $$\Big(\E\Big[\sup_{s\in [0,t]} |Y_s^{(1)}-Y_s^{(2)}|^{k}\Big]\Big) ^{2}\le c_2 \int_0^t  \big\{\|\mu_s-\nu_s\|_{k,var}+ \W_k(\mu_s,\nu_s)\big\}^{2k}  \d s,$$
  which implies  \eqref{XI0-3} since by \eqref{HHH} and the definition of $H^\gg$, there exists a constant $c>0$ depending on $N$ such that
  $$(\E|Y_t^{(1)}-Y_t^{(2)}|^{k}) ^{2}\ge c (\E |X_t^{(1)}-X_t^{(2)}|^k)^2\ge c \W_k(H_t^\gg(\mu), H_t^\gg(\nu))^{2k}.$$

  (e)    Proof of \eqref{XI0-3} for $k>1$ in case $(ii)$. Let $u_t^{\ll,n}$ solve \eqref{POE2} for $L_t=L_{t,\nu}, b^{(0)}= b_t^{(0)}(\cdot,\nu_t)$ and the mollifying approximation $b^{0,n}= b_t^{0,n}(\cdot,\nu_t).$   Then in \eqref{POEE} the equation for $\xi_t$ becomes
  \beg{align*} &\d  \xi_t=  \Big\{\ll   u_t^{\ll, n}(X_t^{(1)})-\ll    u_t^{\ll, n}(X_t^{(2)})+ (b_t^{(0)}-b_t^{0,n})(X_t^{(1)})\\
  &\qquad\qquad -(b_t^{(0)}-b_t^{0,n})(X_t^{(2)})+ b(X_t^{(2)},\mu_t)-b_t(X_t^{(2)},\nu_t)\Big\}\d t \\
 &\qquad +\big\{[(\nn \Theta_t^{\ll, n})\si_t](X_t^{(1)})- [(\nn \Theta_t^{\ll, n})\si_t](X_t^{(2)})\big\}\d W_t+ \n(X_t^{(1)})\d l_t^X-\n(X_t^{(2)})\d l_t^{(2)}.\end{align*}
 So, as shown in step (d)   by \eqref{XI1}, instead of \eqref{LAT},   we have
$$ |X_{t\land \tau_m}^{(1)}- X_{t\land \tau_m}^{(2)}|^{2k}\le   G_m(t) +
   c_2\int_0^{t\land\tau_m} |X_{s\land \tau_m}^{(1)}- X_{s\land \tau_m}^2|^{2k}\d {\tt {\scr L}}_s +\tt M_t  $$ for some local martingale $\tt M_t$,
 $$\tt{\scr L}_t:= \scr L_t+ \int_0^t |f_s(X_s^{(2)})|^2\d s,\ \ t\in [0,T]$$ for $\scr L_t$ in \eqref{LAT'}, and due to $X_0^{(1)}=X_0^{(2)}=X_0$ in the present setting,
$$G_m(t):= \int_0^t \Big\{c_2m^{2(k-1)}\sum_{i=1}^2|b_s^{(0)}-b_s^{0,n}|^2(X_s^{(i)}) +\big(\|\mu_s-\nu_s\|_{k,var} +\W_k(\mu_s,\nu_s)\big)^{2k}\Big\}\d s.$$
By the stochastic Gronwall inequality, Lemma \ref{L1'} and \eqref{WY1},  we find a constant $c>0$ such that
\beq\label{EEE} \beg{split} &\W_k(H_t^\gg(\mu), H_t^\gg(\nu))^{2k}\le (\E|X_t^{(1)}-X_t^{(2)}|^{k})^{2}\\
&\le c \liminf_{m\to\infty}\liminf_{n\to\infty} \E G_m(t)= c\int_0^t \big\{\|\mu_s-\nu_s\|_{k,var}^{2k}+\W_k(\mu_s,\nu_s)^{2k}\big\}\d s.\end{split} \end{equation}
Thus,  \eqref{XI0-3} holds.
\end{proof}

\subsection{Monotone case}

For any $k\ge 0$,
$\scr P_k(\bar D)$ is a complete metric space under the $L^k$-Wasserstein distance $\W_k$, where $\W_0(\mu,\nu):=\ff 1 2 \|\mu-\nu\|_{var}$ and
$$\W_k(\mu,\nu):=   \inf_{\pi\in\C(\mu,\nu)} \bigg(\int_{\bar D\times\bar D} |x-y|^k \pi(\d x,\d y)\bigg)^{\ff 1 {1\lor k}},\ \ \mu,\nu\in \scr P_k(\bar D),\ \ k>0.$$
In the following, we first study the well-posedness  of \eqref{E1} for distributions in $\scr P_k(\bar D)$ with $k>1$, then extend to a   setting including   $k=1$.   

\emph{\beg{enumerate} \item[{\bf (A2)}] Let $k> 1.$ {\bf (D)} holds,
    $b$ and $ \si$ are bounded on bounded subsets of  $[0,\infty)\times\bar D\times\scr P_k(\bar D)$, and  the following two conditions hold.
\item[$(1)$] For any $T>0$ there exists a constant  $K>0$ such that
  \beg{align*} &  \|\si_t(x,\mu)-\si_t(y,\nu)\|^2_{HS}+2\<x-y,b_t(x,\mu)-b_t(y,\nu)\>^+ \\
  &\le K  \big\{|x-y|^2+|x-y|\W_k(\mu,\nu)+1_{\{k\ge 2\}}\W_k(\mu,\nu)^2\big\},\ \ t\in [0,T], x,y\in \bar D, \mu,\nu\in \scr P_k(\bar D).\end{align*}
\item[$(2)$] There exists a subset $\tt\pp D\subset \pp D$ such that
\beq\label{A21} \<y-x, \n(x)\>\ge 0,\ \ x\in \pp D\setminus \tt\pp D,\  y\in \bar D,\end{equation} and when $\tt\pp D\ne \emptyset,$ there exists $\tt\rr\in C_b^2(\bar D)$ such that $\tt\rr|_{\pp D}=0$,
$\<\nn\tt\rr, \n\>|_{\pp D}\ge 1_{\tt\pp D}$ and
\beq\label{A22}  \sup_{(t,x)\in   [0,T]\times \bar D} \big\{\|(\si_t^\mu)^*\nn\tt\rr\|^2(x)+ \<b_t^\mu,\nn \tt\rr\>^-(x)\big\}<\infty,\ \ \mu\in C([0,T]; \scr P_k(\bar D)).\end{equation} \end{enumerate}}

{\bf (A2)}(1) is a monotone condition, when $k\ge 2$ it allows $\si_t(x,\mu)$ depending on $\mu$, but when $k\in [1,2)$ it implies that
$\si_t(x,\mu)=\si_t(x)$ does not depend on $\mu$.

 {\bf (A2)}(2) holds for $\tt\pp D=\emptyset$ when $D$ is convex, and it holds for $\tt\pp D=\pp D$ if $\pp D\in C_b^2$ and for some $r>0$
$$\sup_{(t,x)\in   [0,T]\times \pp_{r_0} D} \big\{\|(\si_t^\mu)^*\nn \rr\|^2(x)+ \<b_t^\mu,\nn \rr\>^-(x)\big\}<\infty,\ \ \mu\in C([0,T]; \scr P_k(\bar D)),$$ where in the second case we may take $\tt\rr=h\circ\rr$ for $0\le h\in C^\infty([0,\infty))$ with $h(r)=r$ for $r\le r_0/2$ and $h(r)=r_0$ for $r\ge r_0$. In general, {\bf (A2)}(2) includes the case where $\pp D$ is partly convex and partly $C_b^2$.

 \beg{thm}\label{T2} Assume {\bf (A2)}. Then   $\eqref{E1}$ is well-posed for distributions in $\scr P_k(\bar D),$ and for any $T>0$, there exist a constant $C>0$ and a map
 $c: [1,\infty)\to (0,\infty)$ such that for any solution $(X_t,l_t)$ of $\eqref{E1}$ with $\L_{X_0}\in \scr P_k(\bar D),$
\beq\label{GPP} \E\Big[\sup_{t\in [0,T] } |X_t|^k \Big] \le C(1+\E|X_0|^k),\end{equation}
\beq\label{LCK} \E \e^{n \tt l_T}\le c(n),\ \ n\ge 1, \tt l_T:= \int_0^T 1_{\tt\pp D}(X_t)\d l_t.\end{equation}
  \end{thm}

  \beg{proof}
  Let $X_0$ be $\F_0$-measurable  with $\gg:=\L_{X_0}\in \scr P_k(\bar D).$   Then
  $$  \scr P_{k,\gg}^T := \big\{\mu\in C([0,T];\scr P_k(\bar D)): \mu_0=\gg\big\}$$
  is a complete   space under the following metric for any $\ll>0$:
   $$\W_{k}^{\ll,T}(\mu,\nu):= \sup_{t\in [0,T]} \e^{-\ll t} \W_k(\mu_t,\nu_t),\ \ \mu,\nu\in \scr P_{k,\gg}^T.$$
  By Lemma \ref{L0}, {\bf (A2)} implies  the well-posedness of   the following reflecting SDE for any $\mu\in \scr P_{k,\gg}^T$:
  \beq\label{GPP1} \d X_t^\mu= b_t(X_t^\mu,\mu_t)\d t +\si_t(X_t^\mu,\mu_t) \d W_t+ \n(X_t^\mu)\d l_{t}^\mu, \ \ X_0^\mu=X_0,\end{equation} and the solution satisfies
\beq\label{MKK}  \E\Big[\sup_{t\in [0,T]}|X_t^\mu|^k\Big]<\infty. \end{equation}
  So, as explained in the proof of Theorem \ref{T3.1}, for the well-posedness of \eqref{E1},  it suffices to prove the contraction of the map
  $$\scr P_{k,\gg}^T\ni\mu\mapsto H(\mu):=\L_{X^\mu}\in \scr P_{k,\gg}^T$$
  under the metric $\W_{k}^{\ll,T}$ for large enough $\ll>0$.

Denote
$$  \tt l_t^\mu:= \int_0^t 1_{\tt\pp D}(X_s^\mu)\d l_s^\mu,\ \ \ \ \tt l_t^\nu:= \int_0^t 1_{\tt\pp D}(X_s^\nu)\d l_s^\nu,\ \ t\ge 0.$$
By   \eqref{LCC},    {\bf (A2)} and It\^o's formula, for any $k\ge 1$ we find a constant $c_1>0$ such that
  \beq\label{GPP2} \d |X_t^\mu-X_t^\nu|^k \le c_1\big\{|X_t^\mu-X_t^\nu|^k+\W_k(\mu_t,\nu_t)^k\big\} \d t + \ff k{r_0} |X_t^\mu-X_t^\nu|^k ( \d \tt l_{t}^\mu+ \d \tt l_t^\nu)+\d M_t\end{equation}
  for some martingale $M_t$ with
  $$\d \<M\>_t\le c_1\big\{|X_t^\mu-X_t^\nu|^{2k}+\W_k(\mu_t,\nu_t)^{2k}\big\}\d t. $$
  To estimate $\int_0^t  |X_s^\mu-X_s^\nu|^k (\d \tt l_s^\mu+\d \tt l_s^\nu),$
  we take
  \beq\label{PHH} 0\le h\in C_b^\infty([0,\infty))\ \text{such\  that\ } h'\le 0,\ h'(0)=-(1+2r_0^{-1}k),\ h(0)=  1,\end{equation}   where $r_0>0$ is in \eqref{LCC}. Let
  $$F(x,y):= |x-y|^k \big\{(h\circ\tt\rr) (x)+(h\circ\tt\rr)(y)\big\},\ \ x,y\in\bar D.$$
  By {\bf (A2)}(2), we have $\tt\rr|_{\pp D}=0$ and $\nn_\n \tt\rr|_{\pp D}\ge 1_{\tt\pp D},$ so that \eqref{PHH} and \eqref{LCC} imply
  $$\nn_{\n} F(\cdot,X_t^\nu)(X_t^\mu) \d l_t^\mu + \nn_{\n} F(X_t^\mu,\cdot)(X_t^\nu)  \d l_t^\nu \le - |X_t^\mu-X_t^\nu|^k (\d \tt l_t^\mu+\d\tt l_t^\nu).$$
   Therefore,  by {\bf (A2)} and   applying It\^o's formula,  we find a constant $c_2>0$ such that
  $$\d F(X_t^\mu,X_t^\nu) \le c_2 \big\{|X_t^\mu-X_t^\nu|^k+\W_k(\mu_t,\nu_t)^k\big\} \d t -   |X_t^\mu-X_t^\nu|^k (\d \tt l_{t}^\mu+\d \tt l_t^\nu)+ \d \tt M_t$$
  for some martingale $\tt M_t$. This  and $F(X_0^\mu,X_0^\nu)=F(X_0,X_0)=0$ imply
 \beq\label{YYP0} \E\int_0^t  |X_s^\mu-X_s^\nu|^k (\d \tt l_s^\mu+\d \tt l_s^\nu)\le   c_2 \int_0^t \big\{\E |X_s^\mu-X_s^\nu|^k+\W_k(\mu_s,\nu_s)^k\big\}  \d s.\end{equation}
  Substituting  \eqref{YYP0}  into \eqref{GPP2} and applying BDG's inequality, we find a constant $c_3>0$ such that
  $$\zeta_t:=\sup_{s\in [0,t]} |X_s^\mu-X_s^\nu|^k,\ \ t\in [0,T]$$ satisfies
  \beq\label{YYP1} \E \zeta_t \le   c_3 \int_0^t \big\{\E\zeta_s + \W_k(\mu_s,\nu_s)^k\big\}  \d s,\ \ t\in [0,T],\end{equation}
  so that for any $\ll>c_3$,
\beq\label{MMK2}   \beg{split} &\E\zeta_t \le c_3 \int_0^t \e^{c_3(t-s)}  \W_k(\mu_s,\nu_s)^k \d s\le  c_3\e^{k\ll t} \W_{k}^{\ll,T}(\mu,\nu)^k \int_0^t\e^{-(k\ll-c_3)(t-s)}\d s\\
  &\le   \ff{c_3\e^{k\ll t}}{k\ll-c_3} \W_{k}^{\ll,T}(\mu,\nu)^k,\ \ t\in [0,T].\end{split} \end{equation}
   Therefore, $H$ is contractive in $\W_{k}^{\ll,T}$ for large $\ll>0$ as desired.

 It remains to prove \eqref{GPP} and \eqref{LCK}. Let $X_t$ be the unique solution to \eqref{E1}. By {\bf (A2)},  for any $k>1$, we find a constant $c(k)>0$ such that
\beq\label{KKL} \d   |X_t|^{k}  \le c(k) \big\{1+|X_t|^k +\E |X_t|^k\big\} \d t + k |X_t|^{k-2} \<X_t,\si_t(X_t,\L_{X_t})\d W_t\>+ k |X_t|^{k-1} \d \tt l_t,\end{equation} where $\d\tt l_t:= 1_{\tt\pp D}(X_t) \d l_t$.
   By applying It\^o's formula to $(1+ |X_t|^k)  (h\circ\tt\rr)(X_t)$,   similarly to      \eqref{YYP0} we obtain
\beq\label{*BU} \E\int_0^t (1+ |X_s|^k) \d \tt l_s  \le \tt c(k) \int_0^t \E\big\{1 + |X_s|^k \big\}\d s\end{equation}
   for some constant $\tt c(k)>0$. Combining \eqref{*BU} with \eqref{KKL} and using Gronwall's lemma, we derive
   $$ \E\Big[\sup_{t\in [0,T]}|X_t|^k\Big]\le c'(1+ \E|X_0|^k)$$ for some constant $c'>0$. Substituting this into \eqref{KKL} and using BDG's inequality, we prove \eqref{GPP} for some constant $c>0$.

 Finally, by {\bf (A1)}(2) and applying It\^o's formula to $ \tt\rr(X_t)$, we prove \eqref{LCK}.
     \end{proof}

We now   solve  \eqref{E1}  for distributions  in
   $$\scr P_{\psi}(\bar D):= \big\{\mu\in \scr P(\bar D): \|\mu\|_\psi:= \mu(\psi(|\cdot|))<\infty\big\},$$
   where  $\psi$ belongs to the following class  for some $\kk>0$:
  \beq\label{PKK} \beg{split}
  \Psi_\kk:= \big\{\psi\in C^2((0,\infty))\cap C^1([0,\infty)):\ & \psi(0)=0, \ \psi'|_{(0,\infty)}>0, \|\psi'\|_\infty<\infty \\
  &\ r\psi'(r)+ r^2\{\psi''\}^+(r) \le \kk \psi(r)\ \text{for}\ r> 0 \big\}.\end{split}\end{equation}
   Let \beq\label{PKK'} \W_\psi (\mu,\nu):=   \inf_{\pi\in\C(\mu,\nu)} \int_{\bar D\times\bar D} \psi(|x-y|) \pi(\d x,\d y),\ \ \mu,\nu\in \scr P_\psi(\bar D).\end{equation}
If $\psi''\le 0$ then $\W_\psi$ is a complete metric on $\scr P_\psi$.  In general, it is only a complete quasi-metric since   the triangle inequality not necessarily holds.

\emph{    \beg{enumerate} \item[{\bf (A3)}]  {\bf (D)} holds, $\si_t(x,\mu)=\si_t(x)$ does not depend on $\mu$,      $b$ and $\si$ are bounded on bounded subsets of
    $[0,\infty)\times\bar D\times\scr P_{\psi}(\bar D)$ for some $\psi\in \Psi_\kk$ and $\kk>0$. Moreover,    for any $T>0$ there exists a constant $ K>0$ such that
  \beg{align*} &  \|\si_t(x)-\si_t(y)\|^2_{HS}+2\<x-y,b_t(x,\mu)-b_t(y,\nu)\>^+ \\
  &\le K  |x-y|\big\{|x-y|+\W_\psi(\mu,\nu)\big\},\ \ t\in [0,T], x,y\in \bar D, \mu,\nu\in \scr P_k(\bar D).\end{align*}
    \end{enumerate}}

  \beg{thm}\label{T3} Assume {\bf (A3)} and {\bf (A2)}$(2)$. Then   $\eqref{E1}$ is well-posed for distributions in $\scr P_\psi(\bar D),$ and
\beq\label{GPP'} \E\Big[\sup_{t\in [0,T] }\psi(|X_t|)\Big] <\infty,\ \ T>0, \L_{X_0}\in \scr P_\psi(\bar D).\end{equation}
  \end{thm}

\beg{proof}   Let $X_0$ be $\F_0$-measurable  with $\E\psi( |X_0|)<\infty$, and consider  the path space
  $$ \scr P_{\psi}^T:= \big\{\mu\in C([0,T];\scr P_\psi(\bar D)): \mu_0=\L_{X_0}\big\}.$$
 For any $\ll>0$, the quasi-metric
   $$\W_{\ll,\psi}(\mu,\nu):= \sup_{t\in [0,T]} \e^{-\ll t} \W_\psi(\mu_t,\nu_t),\ \ \mu,\nu\in \scr P_{\psi}^T$$ is complete.
By Lemma \ref{L0}, {\bf (A3)} implies the well-posedness of  the SDE \eqref{GPP1} for any $\mu\in \scr P_{\psi}^T.$
By {\bf (A2)}(2)  and     It\^o's formula for $\gg_t:= \ss{1+|X_t^\mu-X_0|^2}$, we find a constant $c_1>0$ such that
$$\d \gg_t\le c_1 \{\|\mu_t\|_\psi+\gg_t\}\d t + \gg_t^{-1} \<X_t^\mu-X_0, \si_t(X_t^\mu)\d W_t\>+ \d \tt l_{t}^\mu,$$
where $\d\tt l_t^\mu:= 1_{\tt\pp D}(X_t^\mu)\d l_t^\mu.$
Combining this with $\psi\in \Psi_\kk$ and the linear growth of $\|\si_t\|$ implied by {\bf (A3)}, we find a constant $c_2>0$ such that
\beq\label{GGT} \d \psi(\gg_t) \le c_2 \{\|\mu_t\|_\psi+\psi(\gg_t)\}\d t + \psi'(\gg_t)\gg_t^{-1} \<X_t^\mu-X_0, \si_t(X_t^\mu)\d W_t\>+ \psi'(\gg_t) \d \tt l_{t}^\mu.\end{equation}
Next, by {\bf (A2)}(2), $\psi\in \Psi_\kk$ which implies $\psi'(\gg_t)\le \kk \psi(\gg_t)$ since $\gg_t\ge 1$,
and applying It\^o's formula to $\psi(\gg_t) \{\|\tt\rr\|_\infty-\tt\rr(X_t^\mu)\}$,  we find a constant $c_3>0$ such that similarly to \eqref{YYP0},
\beq\label{LCM} \E\int_0^t \psi'(\gg_s) \d \tt l_s^\mu \le \kk  \E\int_0^t \psi(\gg_s) \d \tt l_s^\mu \le c_3 \E\int_0^t \big\{1+\|\mu_s\|_\psi+ \psi(|X_s^\mu|) \big\} \d s,\ \ t\in [0,T].\end{equation}
Combining this with \eqref{GGT}, $r\psi'(r)\le \kk\psi(r)$, the linear growth of $\si_t$ ensured by {\bf (A3)}, and applying BDG's inequality, we obtain
$$  \E\Big[\sup_{t\in [0,T] }\psi(|X_t^\mu|)\Big] <\infty.$$
Consequently, \eqref{GPP'} holds for solutions of \eqref{E1} with $\L_{X_\cdot}\in \scr P_\psi^T.  $
  So, as explained in the proof of Theorem \ref{T3.1},   it remains  to prove the contraction of the map
  $$\scr P_\psi^T \ni\mu\mapsto H(\mu):=\L_{X^\mu}\in \scr P_\psi^T$$   under the metric $\W_{\ll,\psi}$ for large enough $\ll>0$.

By \eqref{LCC},  {\bf (A2)}(2), $\|\psi'\|_\infty<\infty$  and $r\psi'(r)\le \kk\psi(r)$,  we obtain
\beq\label{*DP} \nn_\n \{\psi(|\cdot-y|)\}(x) \le \ff{\kk}{2r_0}1_{\tt\pp D}(x) \psi(|x-y|),\ \ x\in\pp D, y\in \bar D.\end{equation}
Combining this with    {\bf (A3)}   and It\^o's formula,  we find a constant $c_4>0$ such that
  \beq\label{GPP2B} \d \psi(|X_t^\mu-X_t^\nu|) \le c_4\big\{\psi(|X_t^\mu-X_t^\nu|) +\W_\psi(\mu_t,\nu_t)\big\} \d t + c_4 \psi(|X_t^\mu-X_t^\nu|) (\d \tt l_{t}^\mu+\d \tt l_t^\nu)+\d M_t\end{equation}
  for some martingale $M_t$.

  On the other hand, let $\vv= \ff{r_0}{2\kk}$ and take $h\in C^\infty([0,\infty))$ with $h'\ge 0, h(r)=r$ for $r\le \vv/2$ and $h(r)=\vv$ for $r\ge \vv$. Consider
  $$\eta_t:= \psi(|X_t^\mu-X_t^\nu|) \big\{2\vv -h\circ\tt \rr(X_t^\mu) - h\circ\tt \rr(X_t^\nu)\big\}.$$
  By \eqref{*DP}, {\bf (A2)}(2), $\vv= \ff{r_0}{2\kk}$  and It\^o's formula,   we find a constant $c_5>0$ such that
  \beg{align*} \d\eta_t &\le c_5\big\{\psi(|X_t^\mu-X_t^\nu|)+ \W_\psi(\mu_t,\nu_t)\big\}\d t + \Big(  \ff{2\vv\kk}{2r_0} -1\Big)\psi(|X_t^\mu-X_t^\nu|) (\d \tt l_t^\mu+\d\tt l_t^\nu) +\d\tt M_t\\
  &=  c_5\big\{\psi(|X_t^\mu-X_t^\nu|)+ \W_\psi(\mu_t,\nu_t)\big\}\d t -\ff 1 2 \psi(|X_t^\mu-X_t^\nu|) (\d \tt l_t^\mu+\d\tt l_t^\nu) +\d\tt M_t.\end{align*}
  Since $X_0^\mu=X_0^\nu=X_0$, this implies
   $$  \E\int_0^t \psi( |X_s^\mu-X_s^\nu|) (\d \tt l_{t}^\mu+\d \tt l_t^\nu)\le 2  c_5 \int_0^t \big\{\E\psi( |X_s^\mu-X_s^\nu|)+\W_\psi (\mu_s,\nu_s)\big\}  \d s. $$
  Substituting this into \eqref{GPP2B}, we find a constant $c_6>0$ such that
  $$\W_{\psi}(H_t(\mu), H_t(\nu))\le  \E \psi(|X_t^\mu-X_t^\nu|) \le c_6\int_0^t \W_\psi(\mu_s,\nu_s)\d s,\ \ t\in [0,T],$$
  so that $H$ is contractive in   $\W_{\ll, \psi}$ for large $\ll>0$.  Therefore, the proof is finished.
\end{proof}

 \section{Log-Harnack inequality and applications}

In this section, we  study the log-Harnack inequality introduced in \cite{W10} and applications for  DDRSDEs with singular drift or under monotone conditions.

\subsection{Singular case }

\emph{ \beg{enumerate}\item[{\bf (A4)}] Let $\pp D\in C_b^{2,L}$ and $T>0$.  $\si_t(x,\mu)=\si_t(x)$, and there exists $\hat\mu\in \scr P_2(\bar D)$ such that
  $(A_2^{\si,\hat b})$ holds with $p_1>2$, where $\hat b:= b(\cdot,\hat \mu)$ with regular term $\hat b^{(1)}$.
Moreover, there exist a constant $\aa\ge 0$ and  $1\le f_i  \in \tt L_{q_i}^{p_i}(T), 0\le i\le l$   such that  
\beq\label{BLIP00} |b_t^\mu(x)- \hat b_t^{(1)}(x)|\le  f_0(t,x) +\aa\|\mu\|_2,\ \ \mu\in \scr P_2(\bar D), (t,x)\in [0,T]\times\bar D,\end{equation}
\beq\label{BLIP0}   |b_t^\mu(x)- b_t^\nu(x)|\le     \W_2(\mu,\nu)\sum_{i=1}^l f_i(t,x),\ \ \mu,\nu\in \scr P_2(\bar D), (t,x)\in [0,T]\times\bar D.\end{equation}
\end{enumerate}}

According to Theorem \ref{T1},   {\bf (A4)}  implies the well-posedness of \eqref{E1} up to time $T$ for distributions in $\scr P_2(\bar D).$ Let
$$P_t^*\mu= \L_{X_t}\ \ \text{for \ } X_t\ \text{solving\ \eqref{E1}\ with\ }\L_{X_0}=\mu\in \scr P_2(\bar D), \ \ t\ge 0.$$ We consider
$$P_tf(\mu):= \int_{\bar D} f\d(P_t^*\mu),\ \ t\ge 0, \mu\in \scr P_2(\bar D), f\in \B_b(\bar D),$$
where $\B_b(\bar D)$ is the class of all bounded measurable functions on $\bar D$.

 \beg{thm}  \label{T4}  Assume   {\bf (A4)}.  For any $N>0$, let $\scr P_{2,N}(\bar D):=\{\mu\in \scr P_2(\bar D): \|\mu\|_2\le N\}.$
 \beg{enumerate} \item[$(1)$] For any $N>0$, there exists a constant $C(N)>0$ such that for any $ \nu\in \scr P_{2,N}(\bar D)$   and any $t\in [0,T],$ the following inequalities hold:
 \beq\label{EEE2'} \W_2(P_t^*\mu,P_t^*\nu)^2 \le C(N)\W_2(\mu,\nu)^2,\ \   \mu \in \scr P_2(\bar D),   \end{equation}
 \beq\label{LH'}  P_t\log f(\nu)\le \log P_t f(\mu)+\ff{C(N)}t \W_2(\mu,\nu)^2,
 \ \ 0<f\in \B_b(\bar D),     \mu \in \scr P_{2,N}(\bar D), \end{equation}
 \beq\label{GR'}     \ff 1 2 \|P_t^* \mu-P_t^* \nu\|_{var}^2 \le   {\rm Ent}(P_t^*\nu|P_t^*\mu)
 \le \ff{C(N)}{t} \W_2(\mu,\nu)^2,\ \    \mu  \in \scr P_{2,N}(\bar D), \end{equation}
\beq\label{GR0'}  \|\nn P_tf(\nu)\|_{\W_2}:=\limsup_{\mu\to\nu\ \text{in}\ \W_2} \ff{|P_tf(\nu)-P_tf(\mu)|}{\W_2(\mu,\nu)}
 \le \ff{\ss {2C(N)}}{\ss{t}}\|f\|_\infty,\ \ f\in \B_b(\bar D).
\end{equation}
 \item[$(2)$] Let   $\eqref{BLIP00}$ hold  for $\aa=0$.   Then there exists a constant $C>0$ such that
 \beq\label {EEE2'N} \W_2(P_t^*\mu,P_t^*\nu)^2 \le C\W_2(\mu,\nu)^2,\ \   \mu,\nu \in \scr P_2(\bar D).   \end{equation}
 Moreover, if either $\sup_{1\le i\le l}\|f_i\|_\infty <\infty$ or $D$ is bounded, then $\eqref{LH'}$-$\eqref{GR0'}$ hold for some constant $C$ replacing $C(N)$ and all $\mu,\nu\in \scr P_2(\bar D)$.
 \end{enumerate}
\end{thm}

\beg{proof}  (1) Since the relative entropy of $\mu$ with respect to $\nu$ is given by
$$\Ent(\nu|\mu)= \sup_{g\in \B^+(\bar D), \mu(g)=1}\nu(\log g),$$
\eqref{LH'} is equivalent to
\beq\label{LH''} \Ent(P_t^*\nu|P_t^*\mu)\le \ff{C(N)}t \W_2(\mu,\nu)^2,\ \ t\in (0,T], \mu,\nu\in \scr P_{2,N}(\bar D).\end{equation}
By Pinsker's inequality
$$\ff 1 2 \|\mu-\nu\|_{var}^2\le   \Ent(\nu|\mu),$$
we conclude that \eqref{LH''}   implies  \eqref{GR'}, which further yields \eqref{GR0'}.
So, we only need to prove \eqref{EEE2'} and \eqref{LH''}.

  For any $\mu,\nu\in \scr P_2(\bar D)$, let  $X_t$ solve \eqref{E1}  for $\L_{X_0}=\mu$, and denote
$$\mu_t:=P_t^*\mu=\L_{X_t},\ \ \nu_t:= P_t^*\nu,\ \ \bar\mu_t:=\L_{\bar X_t},\ \ t\in [0,T],$$
where $\bar X_t$ solves
$$\d \bar X_t= b_t(\bar X_t,\nu_t)\d t+\si_t(\bar X_t)\d W_t, \ \ t\in [0,T], \bar X_0=X_0.$$
Let $\si$ and $\hat b:=b(\cdot, \hat\mu)= \hat b^{(1)} + \hat b^{(0)} $ satisfy $(A_2^{\si,\hat b})$.
Consider the decomposition
$$b_t^\nu:=b_t(\cdot,\nu_t)   = \hat b^{(1)}_t +   b^{\nu,0}_t,\ \ b^{\nu,0}_t := b_t^\nu-\hat b_t^{(1)}.$$
By \eqref{ESS} and \eqref{BLIP0}, there exists a constant $K(N)>0$ such that
\beq\label{BON} |b^{\nu,0}_t |\le |\hat b^{(0)}_t |+K(N)f_0(t,\cdot),\ \ \|\nu\|_2\le N,\ \ t\in [0,T].\end{equation}
So, by Theorem \ref{T2.2} and Theorem \ref{T2.3},    the estimate \eqref{INN} and the log-Harnack inequality \eqref{LHA} hold  for solutions of \eqref{E01} with $b^\nu$ replacing $b$ with  a constant   depending on $N$; that is,  there exists a constant $c_1(N)>0$ such that
\beq\label{EEE1} \W_2(\bar\mu_t, \nu_t)^2\le c_1(N) \W_2(\mu,\nu)^2,\ \ t\in [0,T], \mu\in \scr P_2(\bar D),\end{equation}
\beq\label{X1} \Ent(\nu_t|\bar\mu_t)=\sup_{f>0, \bar\mu(f)=1} (P_t f)(\nu)\le \ff{c_1(N)}t \W_2(\mu,\nu)^2,\ \ t\in (0,T],\ \mu\in \scr P_2(\bar D).\end{equation}
Moreover, repeating step (e) in the proof of Theorem \ref{T2.2} for $k=2$ and  $(X_t,\bar X_t)$ replacing $(X_t^{(1)},X_t^{(2)})$, and using
\eqref{BLIP0} replacing \eqref{XI1},   instead of \eqref{EEE} where  $\|\mu_s-\nu_s\|_{k,var}^2$ disappears  in the present case, we derive
$$ \W_2(\mu_t,\bar\mu_t)^4 \le (\E|X_t-\bar X_t|^2)^2 \le   c_2(N)\int_0^t   \W_2(\mu_s,\nu_s)^4 \d s, \ \ t\in [0,T]$$
for some constant $c_2(N)>0.$ This together with \eqref{EEE1} yields
\beg{align*} &\W_2(\mu_t,\nu_t)^4 \le 8 \W_2(\mu_t, \bar\mu_t)^4+ 8\W_2(\bar\mu_t, \nu_t)^2\\
&\le 8c_1(N)^2 \W_2(\mu,\nu)^4+  8 c_2(N) \int_0^t \W_2(\mu_s,\mu_s)^4\d s,\ \ t\in [0,T].\end{align*}
Therefore,  Gronwall's inequality implies \eqref{EEE2'} for some constant $C(N)>0$.

On the other hand,  let $\|\mu\|_2\le N$ and define
\beg{align*} &R_t:=\exp\bigg[ -\int_0^t\<\gg_s,\d W_s\>-\ff 1 2\int_0^t|\gg_s|^2\d s\bigg],\\
&\gg_s:= \big\{\si_s^*(\si_s\si_s^*)^{-1} \big\}(X_s) \big[b_s^\mu(X_s)-b_s^\nu(X_s)\big].
\end{align*} By Girsanov's theorem, we obtain
$$     \int_{\bar D} \Big(\ff{\d \bar\mu_t}{\d\mu_t}\Big)^2 \d\mu_t
 = \E\Big\{ \Big(\ff{\d \bar\mu_t}{\d\mu_t}(X_t)\Big)\Big\}^2= \E\Big\{\Big(\E\big[R_t|X_t]\Big)\Big\}^2\le \E R_t^2.$$
As shown in \cite[p 14-15]{21HS}, by combining this with the Young inequality (see \cite[Lemma 2.4]{ATW09})
\beq\label{Young} \mu(fg)\le \mu(f\log f)+\log \mu(\e^g),\ \ f,g\ge 0, \mu(f)=1,\mu\in \scr P(\bar D),\end{equation}  we derive
\beq\label{X2}  \beg{split} & \Ent(\nu_t|\mu_t)=\int_{\bar D}\log\Big( \ff{\d\nu_t}{\d\mu_t}\Big)\,\d\nu_t
 =\int_{\bar D}\Big\{\log \ff{\d\nu_t}{\d\bar\mu_t}+\log \ff{\d\bar\mu_t}{\d\mu_t}\Big\}\d\nu_t\\
 &=  \Ent(\nu_t|\bar\mu_t) + \int_{\bar D}\Big(\ff{\d\nu_t}{\d\bar\mu_t}\Big) \log \ff{\d\bar\mu_t}{\d\mu_t}\,\d\bar \mu_t
 \le 2 \Ent(\nu_t|\bar\mu_t)+ \log \int_{\bar D} \ff{\d\bar\mu_t}{\d\mu_t}\, \d\bar\mu_t\\
 &= 2 \Ent(\nu_t|\bar\mu_t)+ \log \int_{\bar D} \Big(\ff{\d\bar\mu_t}{\d\mu_t}\Big)^2  \d\mu_t\le 2 \Ent(\nu_t|\bar\mu_t)+\log \E R_t^2.\end{split}\end{equation}
Let $f_s(x):= \sum_{i=1}^l f_i(s,x).$ By \eqref{BLIP0}, \eqref{EEE2'},   $\|\si^*(\si\si^*)^{-1}\|_\infty<\infty$ and \eqref{KAS-2} due to $(A_2^{\si, b^\mu})$, we find   constants $c_3(N),c_4(N)>0$ such that
\beq\label{XPP} \beg{split} & \E[R_t^2]\le \big( \E[R_t^2]\big)^2\le \E\e^{c_3(N)\W_2(\mu,\nu)^2 \int_0^t f_s(X_s)^2 \d s}\\
&\le 1 + \E \bigg[c_3(N)\W_2(\mu,\nu)^2 \bigg(\int_0^t f_s(X_s)^2 \d s\bigg)  \e^{c_3(N)\W_2(\mu,\nu)^2 \int_0^t f_s(X_s)^2\d s}\bigg]\\
&\le 1+ c_3(N) \W_2(\mu,\nu)^2 \bigg[\E\bigg(\int_0^t   f_s(X_s)^2 \d s\bigg)^2\bigg]^{\ff 1 2} \Big[\E \e^{2c_3(N)\W_2(\mu,\nu)^2 \int_0^t f_s(X_s)^2\d s}\Big]^{\ff 1 2}\\
&\le 1+ c_4(N) \W_2(\mu,\nu)^2.\end{split} \end{equation}
Combining this with \eqref{X1} and \eqref{X2},  we prove   \eqref{LH''} for some constant $C(N)>0$.

(2) When $\aa=0$,      \eqref{BON} holds for $K(N)=K$ independent of $N$, so that \eqref{EEE1} and \eqref{X1}  hold for some constant
$C_1(N)=C_1>0$ independent of $N$ and all $\mu,\nu\in \scr P_2(\bar D)$, and in \eqref{XPP} the constant $C_3(N)=C_3$ is independent of $N$ as well. Consequently,  \eqref {EEE2'N} holds and
$$\E[R_t^2]\le \E\e^{C_3\W_2(\mu,\nu)^2 \int_0^t f_s(X_s)^2 \d s}\le \e^{C\W_2(\mu,\nu)^2}$$
if $\sup_{1\le i\le l}\|f_i\|_\infty<\infty$, and when $D$ is bounded we conclude that $C_4(N)=C_4$ in \eqref{XPP} is uniform in $N>0$. Therefore, \eqref{LH'} and hence its consequent inequalities hold for some constant independent of $N$.
\end{proof}


\subsection{Monotone case}

\emph{\beg{enumerate}\item[{\bf (A5)}]  {\bf (D)} and {\bf (A2)}$(2)$ hold,  $\si_t(x,\mu)=\si_t(x)$ does not depend on $\mu$ and is locally bounded on $[0,\infty)\times \bar D$, $\si\si^*$ is invertible,  $b$ is bounded on bounded subsets of $[0,\infty)\times\R^d\times\scr P_2(\bar D),$ and
    for any $T>0$ there exists a  constant $L>0$ such that
 \beg{align*} & \|\si_t(x)-\si_t(y)\|^2_{HS}+2\<x-y,b_t(x,\mu)-b_t(y,\nu)\>^+ \le L  |x-y|^2+L|x-y| \W_2(\mu,\nu),  \\
 & \|\si_t(x)  (\si_t\si_t^*)^{-1}(x)\big\|\le L,\ \  t\in [0,T], x,y\in \bar D, \mu,\nu\in \scr P_2(\bar D).\end{align*}
  \end{enumerate}}
By Theorem \ref{T2},   {\bf (A5)} implies that \eqref{E1} is well-posed for distributions in $\scr P_2(\bar D)$.

 \beg{thm} \label{T5} Assume {\bf (A5)}. Then for any $T>0$, there exists a constant $C>0$ such that the following inequalities hold for all $t\in (0,T]$ and $\nu\in \scr P_2(\bar D)$:
 \beq\label{EEE2}   \W_2(P_t^*\mu,P_t^*\nu)^2 \le C \W_2(\mu,\nu)^2,   \ \ \mu \in \scr P_2(\bar D), \end{equation}
 \beq\label{LH}  P_t\log f(\nu)\le \log P_t f(\mu)+\ff{C }t \W_2(\mu,\nu)^2,
 \ \ 0<f\in \B_b(\bar D),    \mu \in \scr P_2(\bar D), \end{equation}
 \beq\label{GR}     \ff 1 2 \|P_t^* \mu-P_t^* \nu\|_{var}^2 \le   {\rm Ent}(P_t^*\nu|P_t^*\mu)
 \le \ff{C}{t} \W_2(\mu,\nu)^2,\ \  \mu \in \scr P_2(\bar D) \end{equation}
\beq\label{GR0}  \|\nn P_tf(\nu)\|_{\W_2}:=\limsup_{\mu\to\nu\ \text{in}\ \W_2}\ff{|P_tf(\mu)-P_tf(\nu)|}{\W_2(\mu,\nu)} 
  \le \ff{\ss {2C}\|f\|_\infty}{\ss{t}},\ \ f\in \B_b(\bar D).
 \end{equation}
   \end{thm}

\beg{proof}   As explained in the proof of Theorem \ref{T4} that    it suffices to prove \eqref{EEE2} and \eqref{LH}.
To this end, we modify the proof of \cite[Theorem 4.1]{W18} as follows.

Firstly, for $\mu_0,\nu_0\in \scr P_2(\bar D),$  let $(X_0,Y_0)$ be $\F_0$-measurable such that
  \beq\label{INN1} \L_{X_0}=\mu_0,\ \  \L_{Y_0}=\nu_0,\ \ \E|X_0-Y_0|^2=\W_2(\mu_0,\nu_0)^2.\end{equation} Denote
  $$\mu_t:=P_t^*\mu_0,\ \ \nu_t:=P_t^*\nu_0,\ \ t\ge 0.$$
  Let $X_t$ solve \eqref{E1}.  We have
\beq\label{E1'} \d X_t= b_t(X_t, \mu_t)\d t+ \si_t(X_t)\d W_t+ \n(X_t)\d l_t^X,\ \ t\in [0,T],\end{equation}
where $l_t^X$ is the   local time of $X_t$ on $\pp D$. Next,    for any $t_0\in (0,T]$  consider the SDE
\beq\label{CY} \beg{split} \d Y_t = &\Big\{b_t(Y_t,\nu_t)+ \ff {\si_t(Y_t)\{\si_t^*(\si_t\si_t^*)^{-1}\}(X_t)(X_t-Y_t)} {\xi_t} \Big\}\d t\\
&+ \si_t(Y_t) \d W_t +\n(Y_t)\d l_t^Y,\ \ t\in [0,t_0),\end{split}\end{equation}
where   $l_t^Y$ is the local time of $Y_t$ on $\pp D$.  For the constant $L>0$ in {\bf (A5)}, let
\beq\label{Xi} \xi_t:= \ff{1}{L}\Big(1-\e^{L(t-t_0)}\Big),\ \ t\in [0,t_0).\end{equation}
The construction of $Y_t$ goes back to \cite{W11} for the classical SDEs, see also    \cite{W18} for the extension to DDSDEs.
According to Theorem \ref{T01}, {\bf (A5)} implies that \eqref{CY} has a unique solution up to times
$$\tau_{n,m}:=\ff{t_0n}{n+1}\land \inf\big\{t\in [0,t_0): |Y_t|\ge m\},\ \ n,m\ge 1.$$
Let $h$ be in \eqref{PHH} for $k=2$.
By \eqref{LCC} and  {\bf (A2)}(2),   we have
$$\<\nn \big\{(1+h\circ \tt\rr) |\cdot-x_0|^2\big\}(Y_t),\n(Y_t)\>\d l_t^Y\le 0,\ \ x_0\in\bar D,$$   so that {\bf (A5)}, for any $n\ge 1$ we find a constant $c(n)>0$ such that
$$\d \big\{(1+h\circ \tt\rr)(Y_t) |Y_t-x_0|^2\big\}\le c(n) (1+|Y_t|^2) \d t +\d M_t,\ \ t\in [0, \tau_{n,m}], \ n,m\ge 1   $$ holds for some martingale $M_t$. This
 implies   $\lim_{m\to\infty}\tau_{n,m}=\ff{t_0 n}{n+1}$, and hence  \eqref{CY} has a unique solution up to time $t_0$.

Next,   let $\tt Y_t$ solve the SDE
\beq\label{INN'} \d \tt Y_t= b_t(\tt Y_t,\nu_t)\d t + \si_t(\tt Y_t)\d W_t +\n(\tt Y_t) \d l_t^{\tt Y},\ \ \tt Y_0= Y_0, t\in [0,T],\end{equation}
  where $ l_t^{\tt Y}$ is the local time of $\tt Y_t$ on $\pp D$.
  By {\bf (A5)},   \eqref{LCC} and It\^o's formula,  we find a constant $c_2>0$ such that
\beq\label{YYPP} \beg{split}  \E|X_t-\tt Y_t|^2\le &\W_2(\mu_0,\nu_0)^2 + c_2 \int_0^t \big\{\E |X_s-\tt Y_s|^2+\W_2(\mu_s,\nu_s)^2\big\}  \d s\\
&+   \ff 2 {r_0}\E \int_0^t |X_s-\tt Y_s|^2 (\d \tt l_s^X+ \d \tt l_s^{\tt Y}),  \ \ t\in [0,T].\end{split}\end{equation}
For $h$ in \eqref{PHH} with $k=2$, we deduce  from {\bf (A2)}(2) that
\beq\label{*DN} \beg{split} &\big\<\nn \big\{|X_t-\cdot|^2(h\circ\rr(X_t)+ h\circ\rr)\big\}(\tt Y_t),\n(\tt Y_t)\big\>\d \tt l_t^{\tt Y}\le - |X_t-\tt Y_t|^2\d\tt l_t^{\tt Y},\\
&\big\<\nn \big\{|\tt Y_t-\cdot|^2(h\circ\rr(\tt Y_t)+ h\circ\rr)\big\}(X_t),\n(X_t)\big\>\d  \tt l_t^X\le - |X_t-\tt Y_t|^2\d  \tt l_t^X.\end{split}\end{equation}
So,  applying It\^o's formula to
$$\eta_t:= |X_t-\tt Y_t|^2(h\circ\rr(X_t)+ h\circ\rr(\tt Y_t)), $$ and using {\bf (A5)} and \eqref{LCC}, we find a constant $c_3>0$ such that
$$\d \eta_t\le c_3\big\{|X_t-\tt Y_t|^2+\W_2(\mu_t,\nu_t)^2\big\} \d t+\d M_t - |X_t-\tt Y_t|^2 (\d \tt l_t^X+\d\tt l_t^{\tt Y})$$ holds for some martingale $M_t$. This together with \eqref{YYPP} yields
\beg{align*} \E|X_t-\tt Y_t|^2 &\le \W_2(\mu_0,\nu_0)^2 + \E\eta_0 +(c_2+c_3) \int_0^t   \big\{\E |X_s-\tt Y_s|^2+\W_2(\mu_s,\nu_s)^2\big\}  \d s\\
&\le 3\W_2(\mu_0,\nu_0)^2 + 2 (c_2+c_3) \int_0^t    \E |X_s-\tt Y_s|^2 \d s,\ \ t\in [0,T],\end{align*}
where we have used the fact that $\W_2(\mu_s,\nu_s)^2\le \E |X_s-\tt Y_s|^2 $ by definition.
By Gronwall's lemma,   this and $\W_2(\mu_t,\nu_t)^2\le \E|X_t-\tt Y_t|^2$, we find a constant $c_4>0$ such that
 \beq\label{NN0} \W_2(\mu_t,\nu_t)^2\le \E|X_t-\tt Y_t|^2\le c_4\W_2(\mu_0,\nu_0)^2,\ \ t\in [0,T],\end{equation}
so that \eqref{EEE2} holds.

Moreover, for any $n\ge 1$, let
\beq\label{TTA} \tau_n:=\ff{t_0 n}{n+1}\land \inf \{t\in [0,t_0): |X_t-Y_t|\ge n\}.\end{equation}
By Girsanov's theorem,
$$\tt W_t:= W_t+\int_0^{t}   \ff 1 {\xi_s} \{\si_s^*(\si_s\si_s^*)^{-1}\}(X_s)(X_s-Y_s) \d s,\ \ t\in [0,\tau_n]$$
is an  $m$-dimensional Brownian   motion  under the probability $\Q_n:= R_n\P,$ where
\beq\label{TTB} R_n:= \e^{-\int_0^{\tau_n} \ff 1 {\xi_s} \<\{\si_s^*(\si_s\si_s^*)^{-1}\}(X_s)(X_s-Y_s), \d W_s\> - \ff 1 2 \int_0^{\tau_n} \ff{|\{\si_s^*(\si_s\si_s^*)^{-1}\}(X_s)(X_s-Y_s)|^2}{|\xi_s|^2} \d s}.\end{equation}
Then  \eqref{E1'} and \eqref{CY} imply
 \beq\label{INN3} \beg{split} &\d X_t=  \Big\{b_t(X_t,  \mu_t) -\ff {X_t-Y_t} {\xi_t} \Big\}\d t + \si_t(X_t) \d \tt W_t+\n(X_t)\d l_t^X,\\
 &\d Y_t =  b_t(Y_t,\nu_t)\d t + \si_t(Y_t) \d \tt W_t +\n(Y_t)\d l_t^Y,\ \ t\in [0,\tau_n], n\ge 1.\end{split}\end{equation}
Combining this with {\bf (A5)}, \eqref{LCC},   \eqref{NN0} and  It\^o's formula, we obtain
\beq\label{INN4} \beg{split}  & \d \ff{ |X_t-Y_t|^2}{\xi_t}  -  \d M_t \\
&\le   \Big\{\ff{L|X_t-Y_t|^2+L|X_t-Y_t| \W_2(\mu_t,\nu_t))}{\xi_t} -\ff{|X_t-Y_t|^2(2+\xi_t')}{\xi_t^2} \Big\}\d t\\
& \quad   +\ff{|X_t-Y_t|^2}{\xi_t^2} (\d\tt l_t^X+\d\tt l_t^Y)\\
&\le  \Big\{\ff{L^2\W_2(\mu_t,\nu_t)^2}2     -\ff {|X_t-Y_t|^2( 2+\xi_t' -L\xi_t -\ff 1 2)}{\xi_t^2} \Big\}\d t  
    +\ff{|X_t-Y_t|^2}{\xi_t^2} (\d\tt l_t^X+\d\tt l_t^Y)\\
&\le \Big\{\ff{L^2\e^{2Lt}\W_2(\mu_0,\nu_0)^2 } 2   -\ff {|X_t-Y_t|^2} {2\xi_t^2} \Big\}\d t  
   +\ff{|X_t-Y_t|^2}{\xi_t^2} (\d\tt l_t^X+\d\tt l_t^Y),\ \ t\in [0,\tau_n],\end{split}\end{equation}
where $\d M_t:= \ff{2}{\xi_t} \big\<X_t-Y_t, \{\si_t(X_t)-\si_t(Y_t)\}\d \tt W_t\big\>$ is a $\Q_n$-martingale.
By \eqref{*DN} for $(Y_t, \tt l_t^Y)$ replacing $(\tt Y_t, \tt l_t^{\tt Y})$,   and
 applying  It\^o's formula to $\gg_t:= \ff{|X_t-Y_t|^2}{\xi_t} (h\circ\rr(X_t)+h\circ\rr(Y_t))$,
we find a  constant $c_5>0$ such that
$$\d \gg_t\le c_5\gg_t   \d t +\d \tt M_t-\ff{|X_t-Y_t|^2}{\xi_t}(\d \tt l_t^X+\d \tt l_t^Y),\ \ t\in [0,\tau_n],n\ge 1$$
holds for some $\Q_n$-martingale $\tt M_t$. This and  \eqref{INN1} imply that for some constants $c_6,c_7>0$,
\beg{align*} &\E_{\Q_n} \gg_{t\land\tau_n}\le \e^{c_4T}\E\gg_0\le \ff {c_6 }{t_0} \W_2(\mu_0,\nu_0)^2,\\
&\E_{\Q_n} \int_0^{\tau_n}\ff{|X_t-Y_t|^2}{\xi_t}(\d \tt l_t^X+\d \tt l_t^Y)\le \ff{c_7}{t_0}\W_2(\mu_0,\nu_0)^2,\ \ n\ge 1, t\ge 0.\end{align*}
Combining this with \eqref{NN0}, \eqref{INN4} and  {\bf (A5)}, we derive
\beq\label{RR2}\beg{split}  &\E[R_n\log R_n] = \E_{\Q_n} [\log R_n] = \ff 1 2 \E_{\Q_n} \int_0^{\tau_n} \ff{|\{\si_s^*(\si_s\si_s^*)^{-1}\}(X_s)(X_s-Y_s)|^2}{|\xi_s|^2} \d s \\
&\le   \ff{c} {t_0}  \W_2(\mu_0,\nu_0)^2,\ \ n\ge 1\end{split}\end{equation} for some constant $c>0$ uniformly in $t_0\in (0,T].$
Therefore, by the martingale convergence theorem, $R_\infty:=\lim_{n\to\infty} R_n$ exists,
 and
$$N_t:=   \e^{-\int_0^{t} \ff 1 {\xi_s} \<\{\si_s^*(\si_s\si_s^*)^{-1}\}(X_s)(X_s-Y_s), \d W_s\> - \ff 1 2 \int_0^{t} \ff{|\{\si_s^*(\si_s\si_s^*)^{-1}\}(X_s)(X_s-Y_s)|^2}{|\xi_s|^2} \d s},\ \ t\in [0,t_0]$$
is a $\P$-martingale.

Finally, let $\Q:= N_{t_0}\P$. By Girsanov's theorem,
$(\tt W_t)_{t\in [0,t_0]}$ is an $m$-dimensional Brownian motion under the probability $\Q$, and $(X_t)_{t\in [0,t_0]}$ solves the SDE
\beq\label{XXX}\d X_t=  \Big\{b_t(X_t, \mu_t) -\ff {X_t-Y_t} {\xi_t} \Big\}\d t + \si_t(X_t) \d \tt W_t+\n(X_t)\d l_t^X,\ \ t\in [0,t_0].\end{equation}
Let $(Y_t)_{t\in [0,t_0]}$ solve
\beq\label{YYY} \d Y_t =  b_t(Y_t,\nu_t)\d t  + \si_t(Y_t) \d \tt W_t +\n(Y_t)\d l_t^Y,\ \ t\in [0,t_0].\end{equation}
By the well-posedness of \eqref{E1},  this extends the second equation in \eqref{INN3} with   $\L_{Y_{t_0}|\Q}=\nu_{t_0}.$
Moreover,  \eqref{RR2} and Fatou's lemma implies
\beq\label{XDD} \beg{split} &\ff 1 2 \E_{\Q} \int_0^{t_0}  \ff{|\{\si_s^*(\si_s\si_s^*)^{-1}\}(X_s)(X_s-Y_s)|^2}{|\xi_s|^2} \d s \\
&= \E[N_{t_0}\log N_{t_0}]\le \liminf_{n\to\infty} \E[R_n\log R_n]\le  \ff{c} {t_0}   \W_2(\mu_0,\nu_0)^2,\end{split}
\end{equation} which in particular implies
$\Q(X_{t_0}=Y_{t_0})=1$. Indeed, by {\bf (A5)}, if $X_{t_0}(\oo)\ne Y_{t_0}(\oo)$ then there exists a small constant $\vv>0$ such that
$$|\{\si_s^*(\si_s\si_s^*)^{-1}\}(X_s)(X_s-Y_s)|^2(\oo)\ge \vv,\ \ s\in [t_0-\vv,t_0],$$
which implies $\int_0^{t_0}  \ff{|\{\si_s^*(\si_s\si_s^*)^{-1}\}(X_s)(X_s-Y_s)|^2}{|\xi_s|^2} (\oo)\d s=\infty.$ So, \eqref{XDD} implies $\Q(X_{t_0}=Y_{t_0})=1$.
Combining this with the  Young's inequality \eqref{Young}, we  arrive at
\beg{align*} & P_{t_0}\log f(\nu_0) = \E[N_{t_0} \log f(Y_{t_0})] = \E[N_{t_0} \log f(X_{t_0})] \le \E[N_{t_0}\log N_{t_0}] + \log \E[f(X_{t_0})] \\
&\le \log P_{t_0} f(\mu_0) + \ff{c} {t_0}   \W_2(\mu_0,\nu_0)^2,\ \ t_0\in (0,T].\end{align*}
Hence, \eqref{LH} holds.
\end{proof}

\paragraph{Acknowledgement.} The  author would like to thank Bin Xie,  Xing Huang, S. Wang, W. Hong, B. Wu and S. Hu for their helpful comments and corrections to earlier versions of the paper.

  \end{document}